\documentclass[12pt]{article}
\usepackage{amscd,amsthm, amsmath, amsfonts,amssymb,
epsfig}
\usepackage{calc,ifthen,graphics,miniplot,wrapfig,
subfigure, rotfloat
}

\begin{document}

\large
\setlength{\myStandardFigureWidth}{\linewidth}

\title{Catalog of dessins d'enfants with $\le$ 4 edges}

\author{N. M.~Adrianov
\and N. Ya.~Amburg
\and V. A.~Dremov
\and Yu. A.~Levitskaya
\and E. M.~Kreines
\and Yu. Yu.~Kochetkov
\and V. F.~Nasretdinova
\and  G. B.~Shabat}
\date{}

\maketitle
{\it
\begin{center}
ITEP, MSU, Moscow, Russia\\ \ \\
\rm E-mails: \it   kreines@itep.ru, shabat@mccme.ru
\end{center}
}
\begin{abstract}
In this work all the dessins d'enfant with no more than 4 edges are listed and
their Belyi pairs are computed. In order to enumerate all dessins the technique
of matrix model computations was used. 
The total number of dessins is 134; among them 77 are spherical, 53 of genus 1
and 4 of genus 2. The orders of automorphism groups of all the dessins are also
found.

Dessins are listed by the number of edges. Dessins with the same number of
edges are  ordered lexicographically by their lists of 0-valencies. The
corresponding matrix model for any list of 0-valencies is given and computed.
Complex matrix models for dessins with 1 -- 3 edges are used. For the dessins
with 4 edges we use Hermitian matrix model, correlators for which are computed
in~\cite{Alex}.

\end{abstract}
%
%
%

\section{Introduction}

{\em Dessin d'enfant\/} is a compact connected smooth oriented surface $S$
together with a graph $\Gamma$  on it such that the complement $S\setminus
\Gamma$ is homeomorphic to a disjoint union of open discs.  The theory of
Dessins d'enfants was  initiated by A.~Grothendieck in \cite{G,Sketch} and
actively developed thereafter, see \cite{LZ} and references therein.
 Dessins d'enfants became rather popular within the last
decades; they provide a possibility to describe in the easy and visually
effective combinatorial language of graphs on surfaces many difficult and deep
concepts and results of Inverse Galois theory, Teichm\"uller and moduli spaces,
Maps and hypermaps, Matrix models, Quantum gravity, String theory, etc.

Dessins d'enfant 
appear naturally in
different branches of mathematics. 
A smooth irreducible complete complex algebraic curve, defined over the field
$\overline{\mathbb{Q}}$, provides a dessin d'enfant in the following way. On
such a curve $X$ according to the famous Belyi theorem there exists a
nonconstant rational function $\beta$  having at most 3 critical values. Denote
$X_{\mathbb{C}}$ its complexification and $\beta_{\mathbb{C}}$ the natural lift
of $\beta$ to $X_{\mathbb{C}}$. By definition, such $\beta$'s and
$\beta_{\mathbb{C}}$'s are {\em Belyi functions.\/} 
Without loss of generality we assume that the critical values of $\beta$ are in
$\{0,\: 1,\: \infty\}$. Moreover, replacing $\beta$ by $4\beta(1-\beta)$, if
needed, we can assume that $1 - \beta$ {\em has only double zeros\/}; such
Belyi functions are called {\em clean\/}. Then
$\beta_{\mathbb{C}}^{-1\circ}([0,1])$ is a dessin d'enfant on the topological
model of $X_{\mathbb{C}}$ whose edges are
$\{\beta_{\mathbb{C}}^{-1\circ}([0,1])\}$ and vertices are
$\{\beta_{\mathbb{C}}^{-1\circ}(0)\}$. In the main text we omit the
complexification subscripts.

In this work, using the matrix model approach, (see~\cite{Alex,matrix,DIFRAN})
we listed all the dessins d'enfants with no more than 4 edges. There are two 1-edge dessins, both of them are of genus zero, fifteen 2-edge
dessins, among them only one is of genus 1, twenty 3-edge dessins: 14 sperical and
6 of genus 1, and one hundred seven 4-edge dessins: $57$ spherical dessins, $46$ dessins of
genus 1, and 4 dessins of genus 2.
The total number of dessins is 134. The main result
is the calculation of the corresponding Belyi \it pairs \rm (in the case of
positive genus it means finding the curve and the Belyi function on it). This
catalog is the analog of the well-known Betrema-Per\'e-Zvonkine
catalog~\cite{BPZ}, where the trees with no more than 8 edges with their Belyi
functions are collected.

\section*{Acknowledgments}
We would like to thank all participants of the scientific seminar ``Graphs on surfaces and curves over number fields'' for the encouragement and interesting discussions. The 2nd, 5th, 7th, and 8th authors are grateful to Russian Federal Nuclear Energy Agency for financial support. Also the 2nd author would like to thank the grants NWO-RFBR 047.011.2004.026 (RFBR 05-02-89000-NWO-a), NSh-8004.2006.2, RFBR 07-01-00441-a, the 3rd, 4th, 5th and 8th authors would like to thank the grants  RFBR 07-01-00441-a and NSH-1999.2006.1, the 5th author would like to thank the grant MK-2687.2007.1, the 7th author would like to thank the grants  RFBR 07-01-00441-a and NSh-8065.2006.2. The 5th author would like to express her deep gratitude to Max-Plank Institute of Mathematics in Bonn, where some part of this work was done, for worm scientific atmosphere and financial support. 

\section{Comments on Belyi functions}

For genus 0 dessins it is natural to write a Belyi function as a fraction of
two polynomials. For the simplification of checking we factorize these
polynomials.

Since all the curves of genus 1 and 2 are hyperelliptic, we always write their
equations in the form $y^2=F(x)$ and denote by $\tau$ the hyperelliptic
involution $\tau:(x,y)\mapsto(x,-y)$.

For the curves of positive genus there two different cases:
\begin{itemize}
\item
Let $\beta$ be invariant with respect to a hyperelliptic involution
 $\tau$ (i.e. $\beta=\beta \circ \tau$), then $\beta$ can be written as a quatient of two polynomials
depending on the coordinate in a quatient space (by this involution), which is a projective line. In this
case we write $(X:y^2=F(x),\beta=\frac{P(x)}{Q(x)})$.
\item If Belyi function is not invariant with respect to $\tau$, it is convenient to use
symmetric functions in $\beta$ and $\beta^{\tau}$:
$$n_0=\beta \cdot \beta^{\tau},\
n_1=(\beta-1) (\beta^{\tau}-1).
$$
 It is easy to reconstruct the Belyi function
from this pair and the equation of the curve. Indeed,
\begin{equation}
\beta=\frac{n_0-n_1+1}2 +
      y \sqrt{\frac{ (n_0-n_1)^2-2(n_0+n_1)+1}{4F}}
\label{hebeta}
\end{equation}
For the convenience of the reader we give all above objects in these cases, namely:\\
$X:\{y^2=F(x)\}$, $\beta=\frac{P(x)+Q(x)y}{R(x)}$,
$n_0=\frac{P(x)^2-Q(x)^2F(x)}{R(x)^2}$ ,
$n_1=\frac{(P(x)-R(x))^2-Q(x)^2F(x)}{R(x)^2}$.
\end{itemize}

It is also useful to factorize numerators and denominators of functions $n_0$
and $n_1$. Note that the degrees of factors in the numerator of $n_0$ are
related with the valencies of the vertices of the dessin, and degrees of
factors in denominator are related to the valencies of faces.

\section{The list of dessins with 1, 2, 3 edges}
\subsection{1-edge dessins}
\begin{align*}\langle Tr(Z^2)Tr((Z^{+})^{2})\rangle=
\langle\langle Tr(Z^2)Tr((Z^{+})^{2})\rangle\rangle=\\
=2\cdot2\left(\frac{1}{2}N^2\right).\\
\end{align*}
\begin{figure}[h]
\begin{minipage}[b]{.45\linewidth}
\centering\epsfig{figure=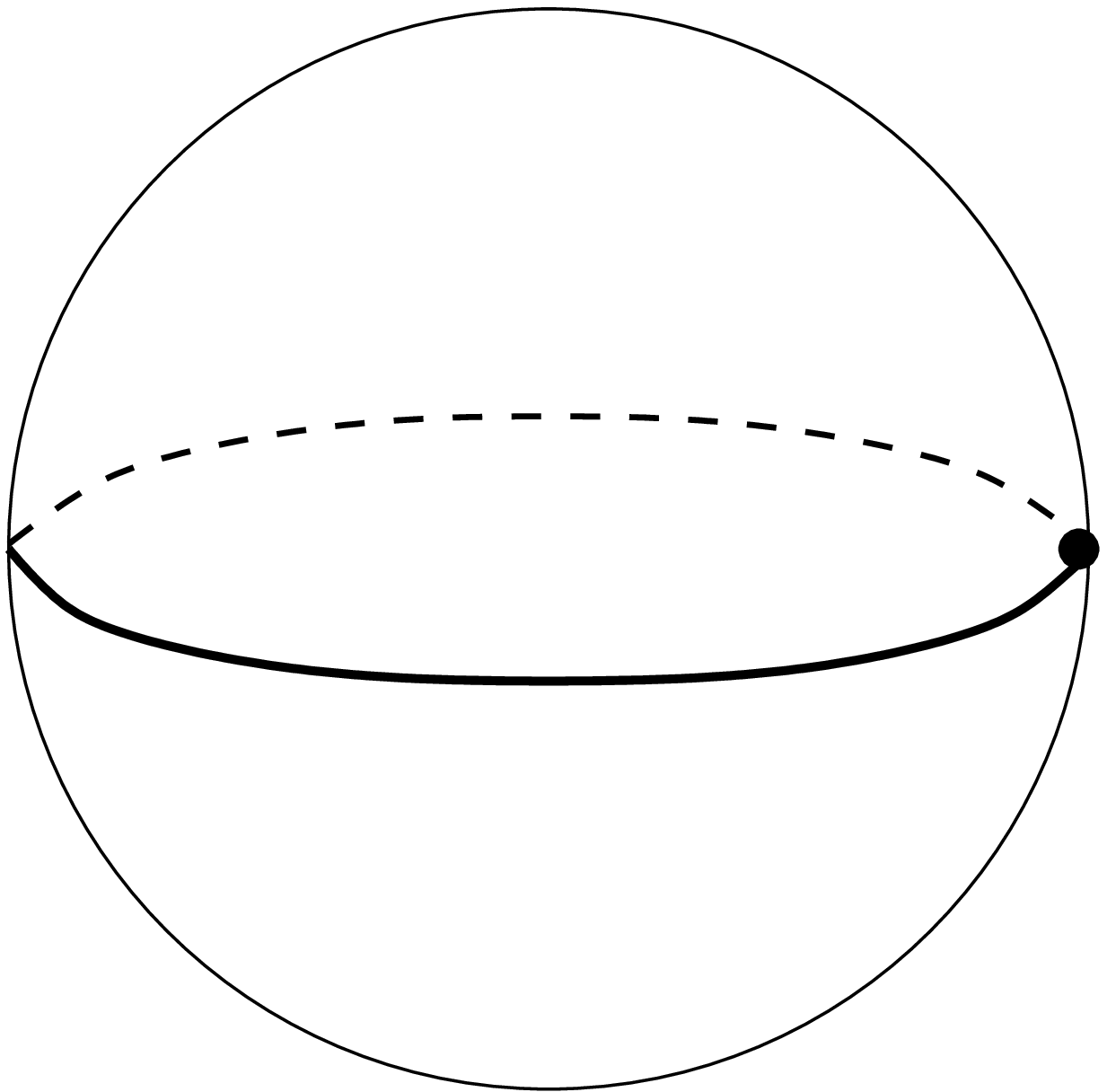,width=\linewidth}
\caption{$S(2|11)$. Valencies $(2|1,1)$. The order of the automorphism group: $2$. Dual dessin
  $S(11|2)$, see Figure \ref{one_v11_f2_sf} on the page \pageref{one_v11_f2_sf}.
Belyi function is $\beta=\frac{z^2}{z^2-1}$.
}\label{one_v2_f11_sf}

\end{minipage}\hfill
\end{figure}
\clearpage

\begin{align*}\langle Tr^2(Z)Tr((Z^{+})^{2})\rangle=
\langle\langle Tr^2(Z)Tr((Z^{+})^{2})\rangle\rangle=\\
=2!\cdot2\left(\frac{1}{2}N\right).\\
\end{align*}
\begin{figure}[h]
\begin{minipage}[b]{.45\linewidth}
\centering\epsfig{figure=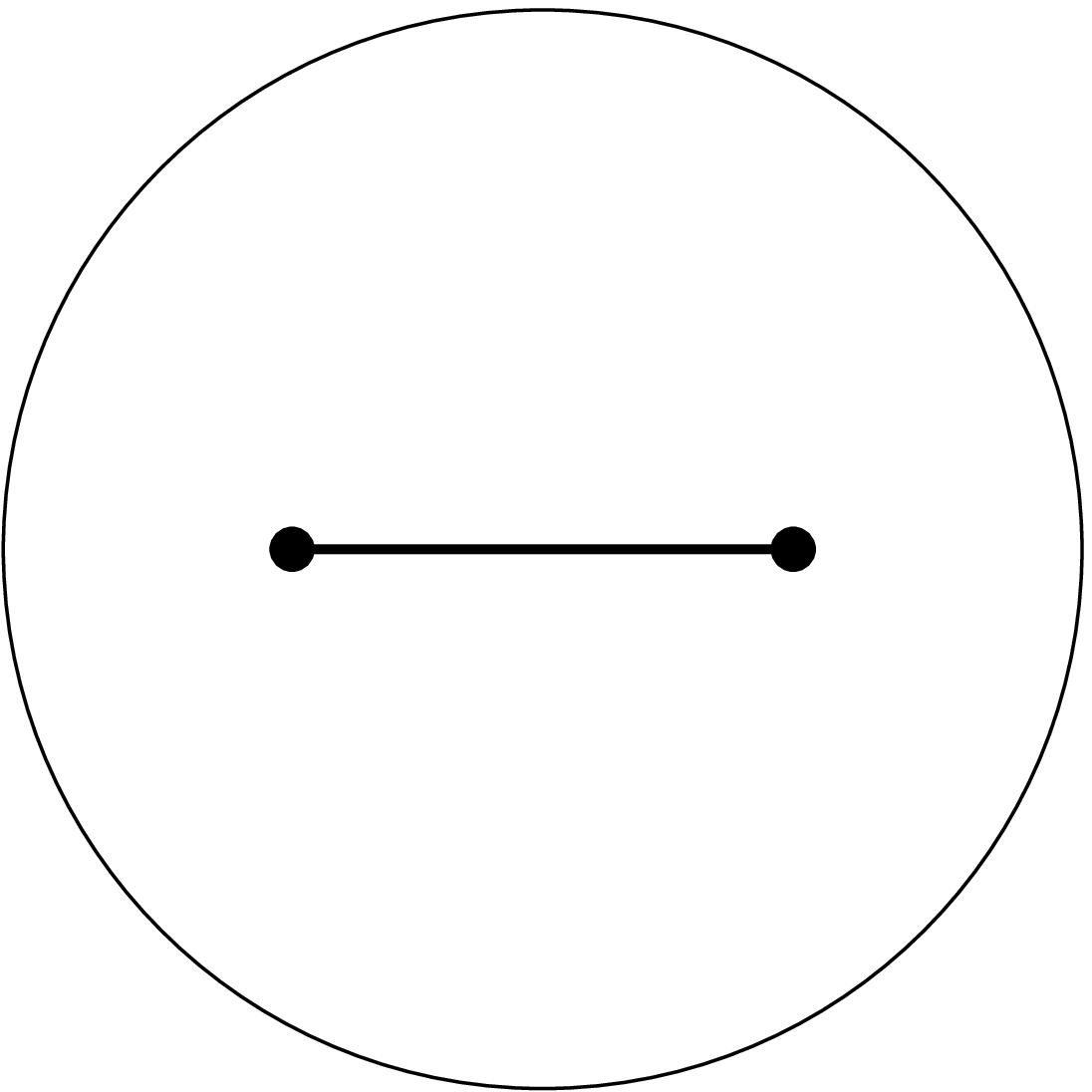,width=\linewidth}
\caption{$S(11|2)$. Valencies $(1,1|2)$. The order of the automorphism group: $2$. Dual dessin
  $S(2|11)$, see Figure \ref{one_v2_f11_sf} on the page \pageref{one_v2_f11_sf}.
Belyi function is $\beta=1-z^2$.
}\label{one_v11_f2_sf}

\end{minipage}\hfill
\end{figure}
\clearpage

\subsection{2-edge dessins}
\begin{align*}\langle Tr(Z^4)Tr^2((Z^{+})^{2})\rangle=
\langle\langle Tr(Z^4)Tr^2((Z^{+})^{2})\rangle\rangle=\\
=4\cdot2!\cdot2^2\left(\frac{1}{2}N^3+\frac{1}{4}N\right).\\
\end{align*}
\begin{figure}[h]
\begin{minipage}[b]{.50\linewidth}
\centering\epsfig{figure=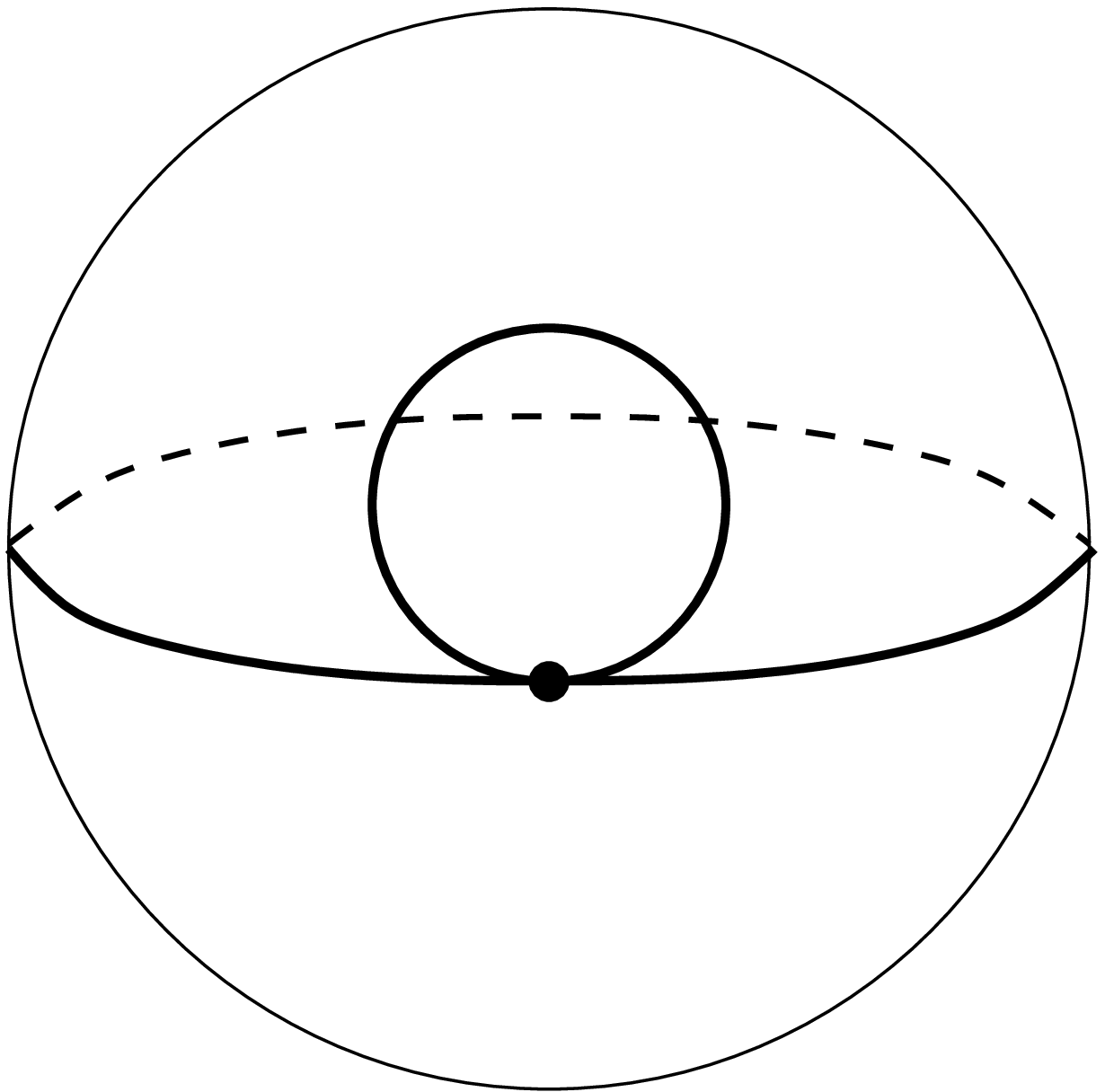,width=.85\linewidth}
\caption{$S(4|211)$. Valencies $(4|2,1,1)$. The order of the automorphism group: $2$.
Dual dessin   $S(211|4)$, see Figure \ref{two_v211_f4_sf} on the page \pageref{two_v211_f4_sf}.
Belyi function is $\beta=\frac{1}{-4z^2(z^2-1)}$.
}\label{two_v4_f211_sf}
\end{minipage}\hfill
\end{figure}
\begin{arrangedFigure}{1}{2}{}{
{$T(4,4)$. Valencies $(4|4)$. The order of the automorphism group: $4$. Dual dessin
$T(4|4)$, see Figure \ref{two_v4_f4_tor} on the page \pageref{two_v4_f4_tor}.
Belyi function is $(X:y^2=x^3-x, \beta = x^2)$.
}\label{two_v4_f4_tor}}
\subFig{two_v4_f4_tor}%
\subFig{two_v4_f4_tor_ucover}%
\end{arrangedFigure}

\clearpage
\begin{align*}
\langle Tr(Z^3)Tr(Z)Tr^2((Z^{+})^{2})\rangle=
\langle\langle Tr(Z^3)Tr(Z)Tr^2((Z^{+})^{2})\rangle\rangle
=3\cdot2!\cdot2^2\left(N^2\right).\\
\end{align*}
\begin{figure}[h]
\begin{minipage}[b]{.45\linewidth}
\centering\epsfig{figure=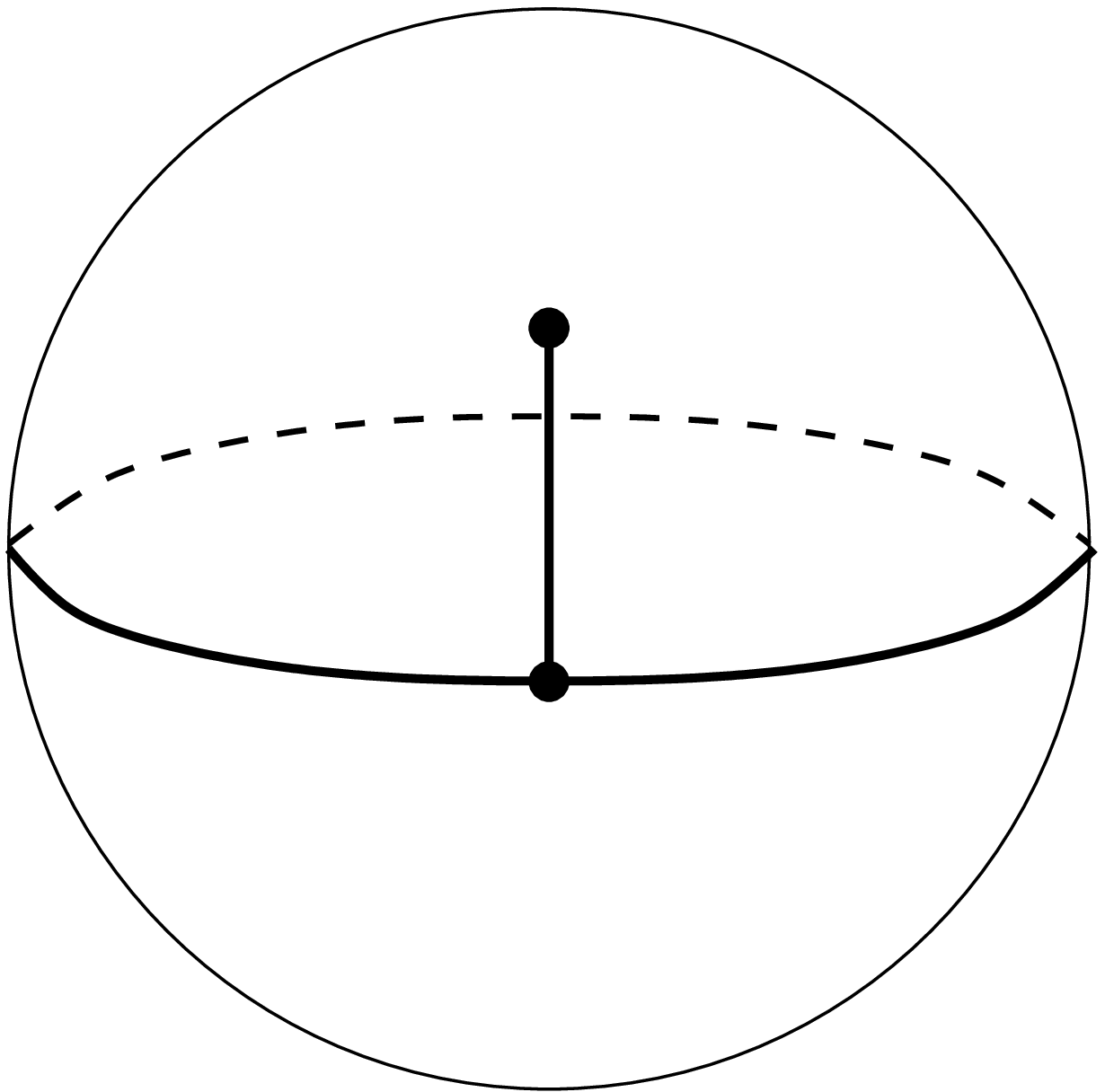,width=.85\linewidth}
\caption{$S(31|31)$. Valencies $(3,1|3,1)$. The order of the automorphism group: $1$.
Dual dessin   $S(31|31)$, see Figure \ref{two_v31_f31_sf} on the page \pageref{two_v31_f31_sf}.
Belyi function is $\beta = -64 \frac{z^3(z-1)}{8z+1}$.
}\label{two_v31_f31_sf}

\end{minipage}\hfill
\end{figure}
\clearpage
\begin{align*}
\langle\langle Tr^2(Z^2)Tr^2((Z^{+})^{2})\rangle\rangle
=2!\cdot2^2\cdot2!\cdot2^2\left(\frac{1}{4}N^2\right).\\
\end{align*}
\begin{figure}[h]
\begin{minipage}[b]{.45\linewidth}
\centering\epsfig{figure=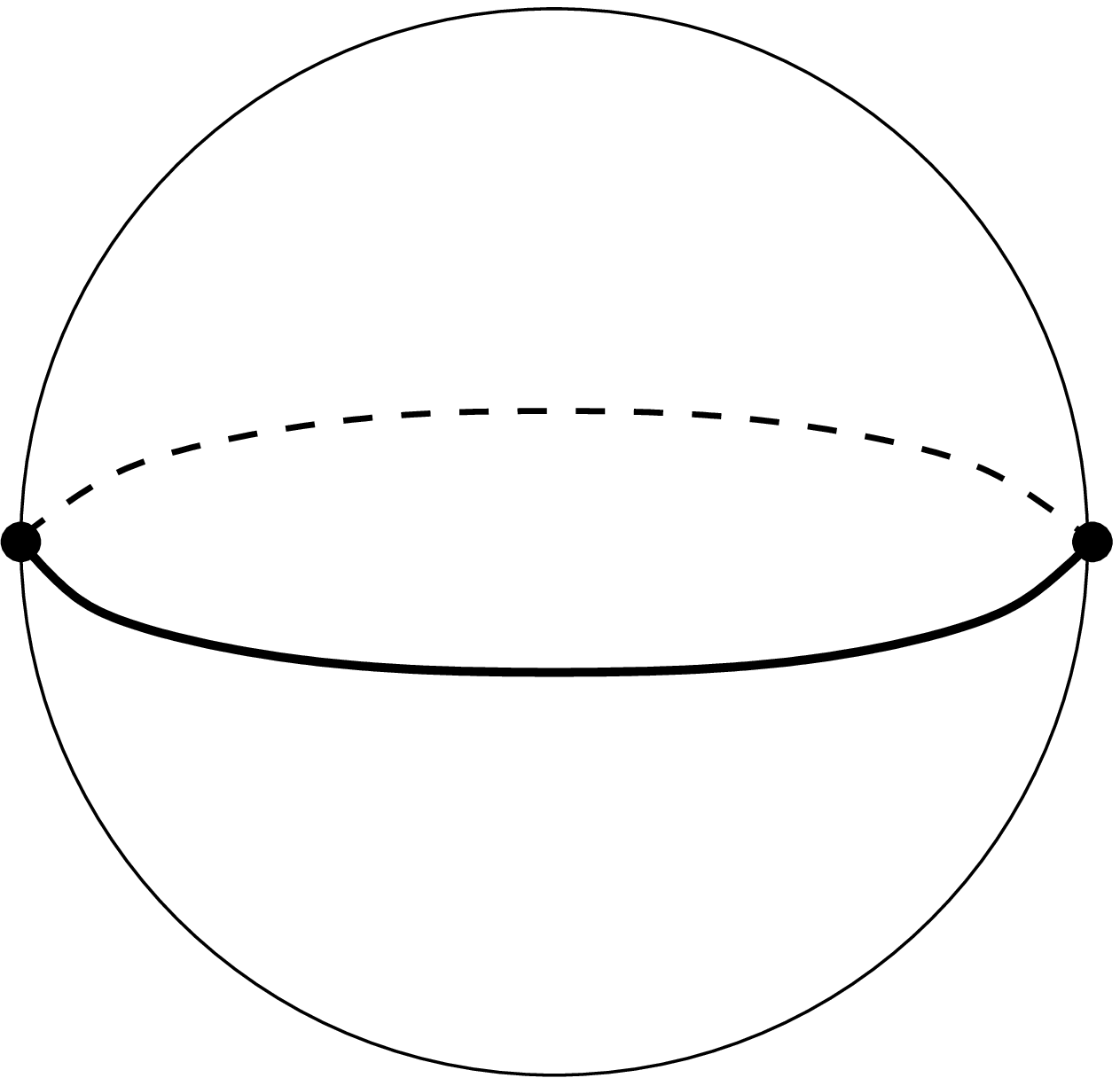,width=.85\linewidth}
\caption{$S(22|22)$. Valencies $(2,2|2,2)$. The order of the automorphism group: $4$.
Dual dessin   $S(22|22)$, see Figure \ref{two_v22_f22_sf} on the page \pageref{two_v22_f22_sf}.
Belyi function is $\beta = \frac{(z^2-1)^2}{-4z^2}$.
}\label{two_v22_f22_sf}

\end{minipage}\hfill
\end{figure}
\clearpage
\begin{align*}
\langle\langle Tr(Z^2)Tr^2(Z)Tr^2((Z^{+})^{2})\rangle\rangle
=2\cdot2!\cdot2!\cdot2^2\left(\frac{1}{2}N\right).\\
\end{align*}
\begin{figure}[h]
\begin{minipage}[b]{.45\linewidth}
\centering\epsfig{figure=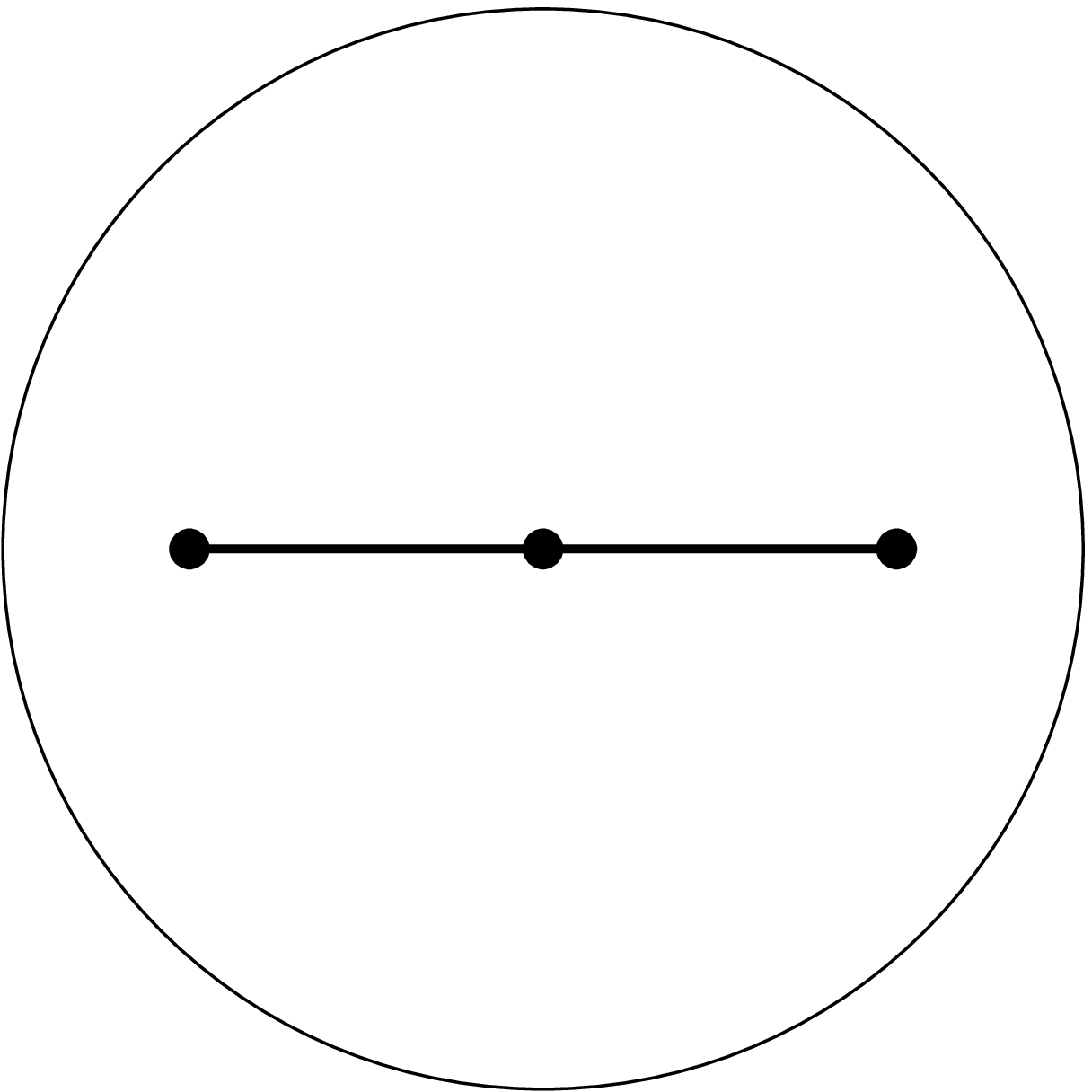,width=.85\linewidth}
\caption{$S(211|4)$. Valencies $(2,1,1|4)$. The order of the automorphism group: $2$.
Dual dessin   $S(4|211)$, see Figure \ref{two_v4_f211_sf} on the page \pageref{two_v4_f211_sf}.
Belyi function is $\beta=-4z^2(z^2-1)$.
}\label{two_v211_f4_sf}
\end{minipage}\hfill
\end{figure}
\clearpage

\subsection{3-edge dessins}
\begin{align*}\langle Tr(Z^6)Tr^3((Z^{+})^{2})\rangle=
\langle\langle Tr(Z^6)Tr^3((Z^{+})^{2})\rangle\rangle=\\
=6\cdot3!\cdot2^3\left(\frac{5}{6}N^4+\frac{5}{3}N^2\right)=\\
=6\cdot3!\cdot2^3\left(\left(\frac{1}{3}+\frac{1}{2}\right)N^4+\left(1+\frac{1}{2}+\frac{1}{6}\right)N^2\right)
\end{align*}

\begin{figure}[h]
\begin{minipage}[b]{.45\linewidth}
\centering\epsfig{figure=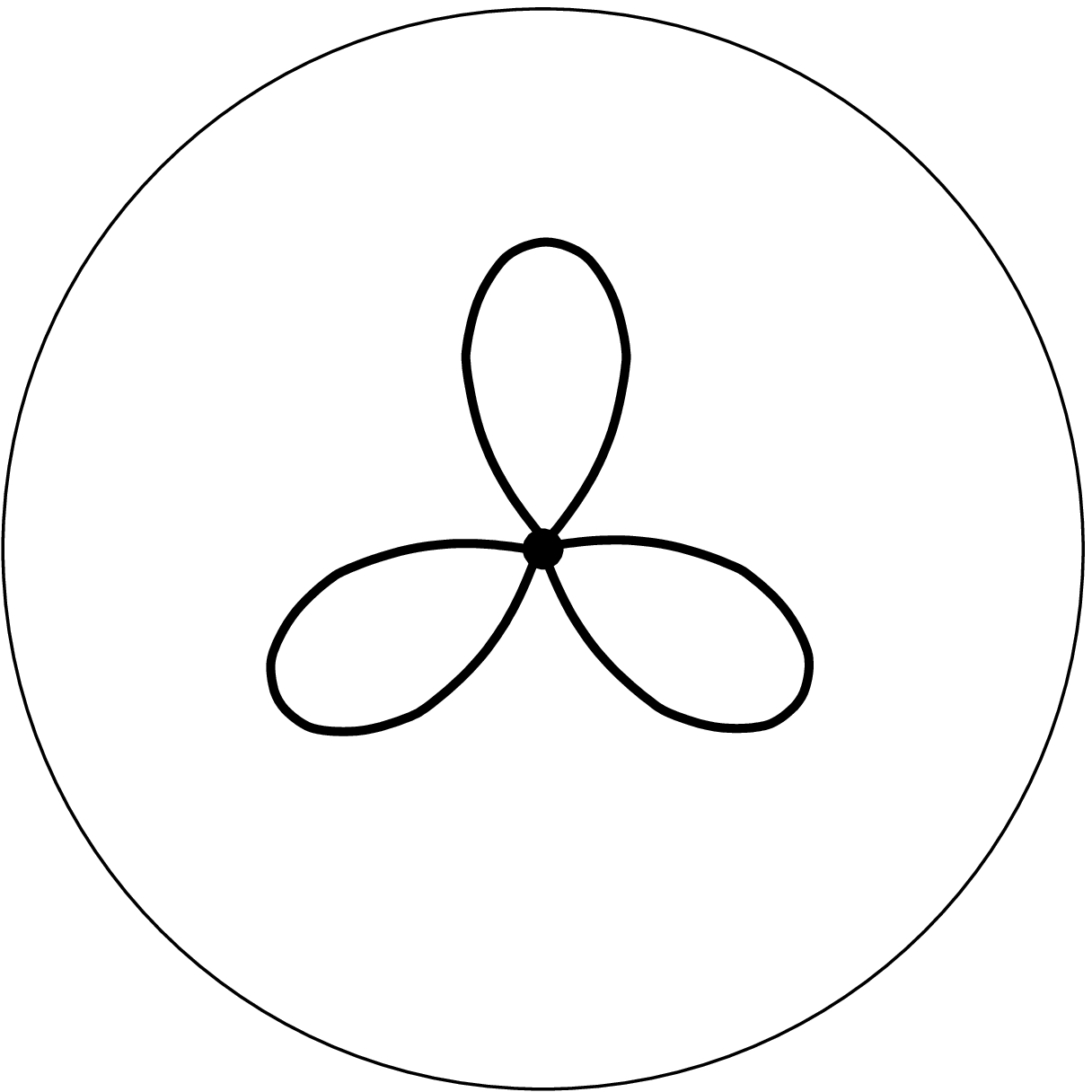,width=\linewidth}
\caption{$S(6|3111)$. Valencies $(6|3,1,1,1)$. The order of the automorphism group: $3$.
Dual dessin   $S(3111|6)$, see  Figure \ref{three_v3111_f6_sf} on the page
\pageref{three_v3111_f6_sf}.
Belyi function is $\beta=\frac{1}{-4z^3(z^3-1)}$.
} \label{three_v6_f3111_sf}
\end{minipage}\hfill
\begin{minipage}[b]{.45\linewidth}
\centering\epsfig{figure=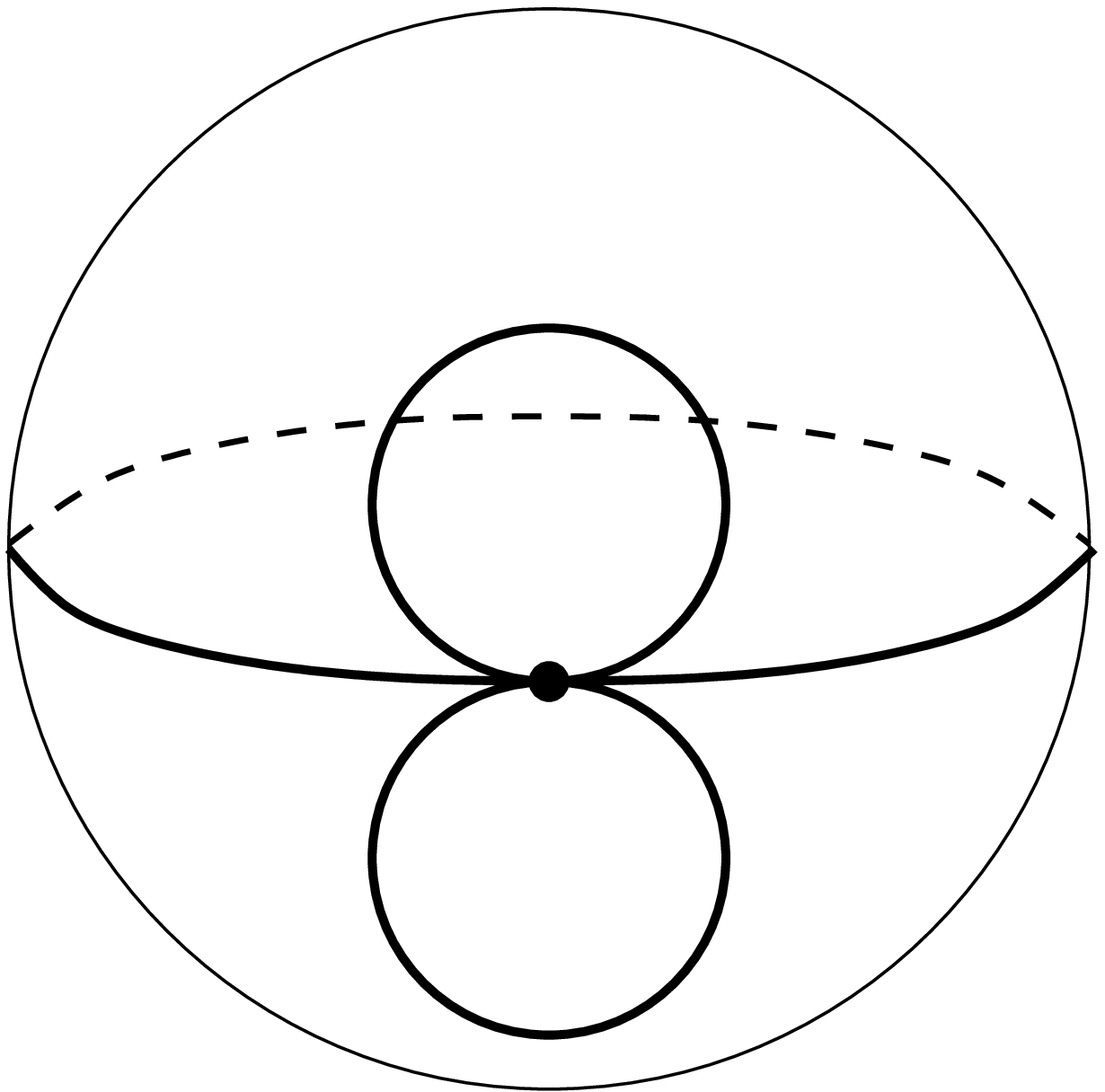,width=.85\linewidth}
\caption{$S(6|2211)$. Valencies $(6|2,2,1,1)$. The order of the automorphism group: $2$.
Dual dessin   $S(2211|6)$, see  Figure \ref{three_v2211_f6_sf}
on the page \pageref{three_v2211_f6_sf}.
Belyi function is $\beta=\frac{-4}{(z^2-1)^2(z^2-4)}$.
} \label{three_v6_f2211_sf}
\end{minipage}\hfill
\end{figure}
\begin{arrangedFigure}{2}{2}{}{
{$T(6|51)$. Valencies $(6|5,1)$. The order of the automorphism group: $1$. Dual dessin
$T(51|6)$, see Figure \ref{three_v51_f6_tor} on the page \pageref{three_v51_f6_tor}.
Belyi function is
$\beta=-{\frac {1}{216}}\,{\frac { \left( 2\,{x}^{3}+2\,yx-29\,{x}^{2}-
15\,y+85\,x-50 \right) {x}^{3}}{2\,x+1}}$   on the curve
$X:{y}^{2}={x}^{4}-14\,{x}^{3}+29\,{x}^{2}-60\,x$.
$n_0={\frac {625}{11664}}\,{\frac {{x}^{6}}{2\,x+1}}$,
$n_1={\frac {1}{11664}}\,{\frac { \left( 29\,{x}^{3}-54\,{x}^{2}+
108\,x+108 \right) ^{2}}{2\,x+1}}$
}\label{three_v6_f51_tor}}
\subFig{three_v6_f51_tor}%
\subFig{three_v6_f51_tor_ucover}%
\newSubFig{}{{$T(6|42)$. Valencies $(6|4,2)$. The order of the automorphism group: $2$.
Dual dessin   $T(42|6)$, see Figure \ref{three_v42_f6_tor} on the page \pageref{three_v42_f6_tor}.
Belyi function is $(X:y^2=x(x+3)(x-1),\beta = \frac{4 x^3}{27 (x-1)})$.
} \label{three_v6_f42_tor}}

\subFig{three_v6_f42_tor}%
\subFig{three_v6_f42_tor_ucover}%
\newSubFig{}{{$T(6|33)$. Valencies $(6|3,3)$. The order of the automorphism group: $6$.
Dual dessin   $T(33|6)$, see Figure \ref{three_v33_f6_tor}
on the page  \pageref{three_v33_f6_tor}.
Belyi function is $(X:y^2=x^4-x,\beta = x^3)$.
} \label{three_v6_f33_tor}}
\subFig{three_v6_f33_tor}%
\subFig{three_v6_f33_tor_ucover}%
\end{arrangedFigure}
\clearpage
\begin{align*}\langle Tr(Z^5)Tr(Z)Tr^3((Z^{+})^{2})\rangle=
\langle\langle Tr(Z^5)Tr(Z)Tr^3((Z^{+})^{2})\rangle\rangle=\\
=5\cdot3!\cdot2^3\left(2N^3+N\right)
=5\cdot3!\cdot2^3\left((1+1)N^3+N\right).
\end{align*}

\begin{figure}[h]
\begin{minipage}[b]{.45\linewidth}
\centering\epsfig{figure=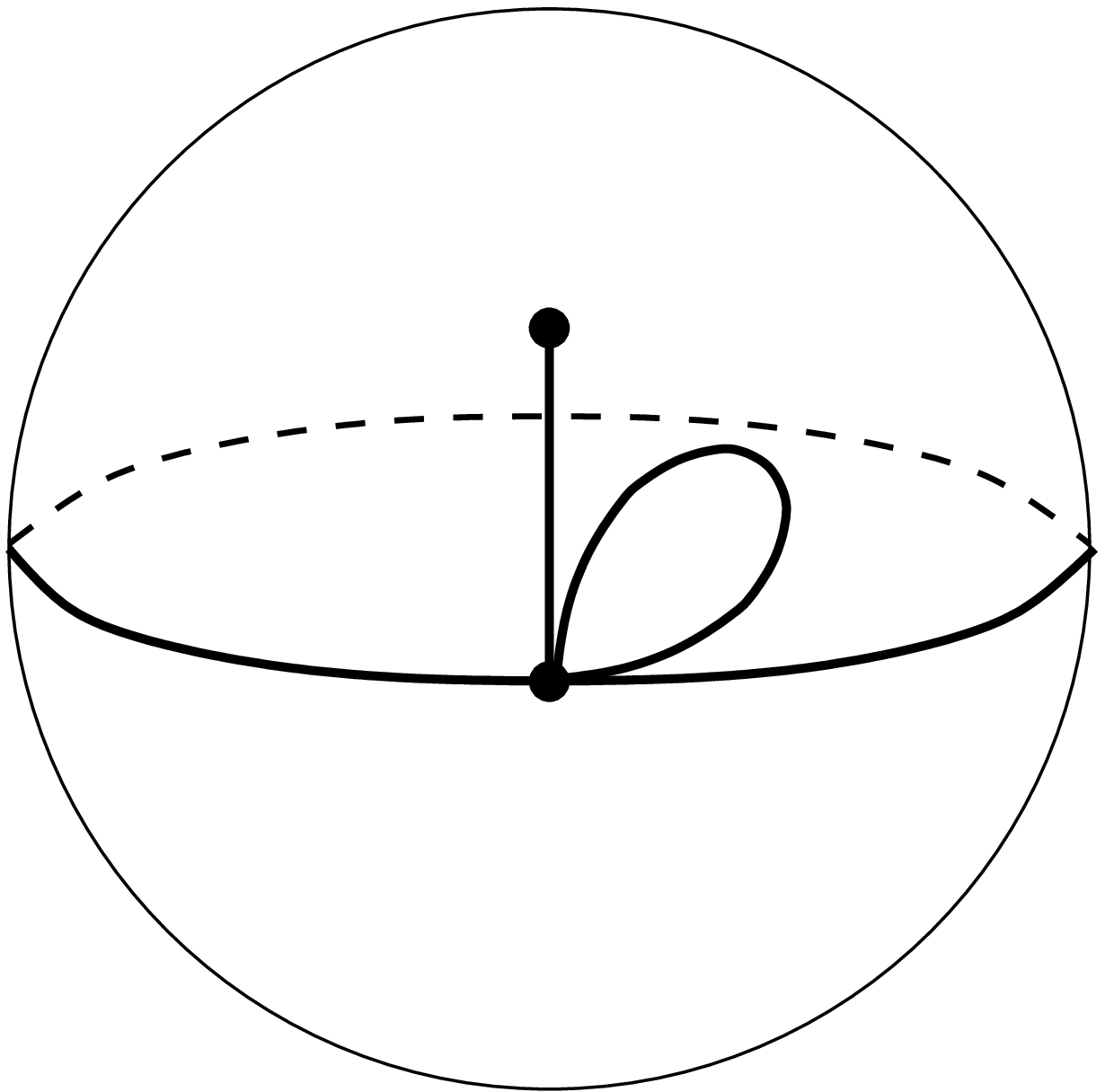,width=.85\linewidth}
\caption{$S(51|411)$. Valencies $(5,1|4,1,1)$. The order of the automorphism group: $1$.
Dual dessin   $S(411|51)$, see Figure \ref{three_v411_f51_sf}
on the page  \pageref{three_v411_f51_sf}.
Belyi function is $\beta = \frac{256(z-1)}{z^4 (z^2+4z+20)}$.
} \label{three_v51_f411_sf}
\end{minipage}\hfill
\begin{minipage}[b]{.45\linewidth}
\centering\epsfig{figure=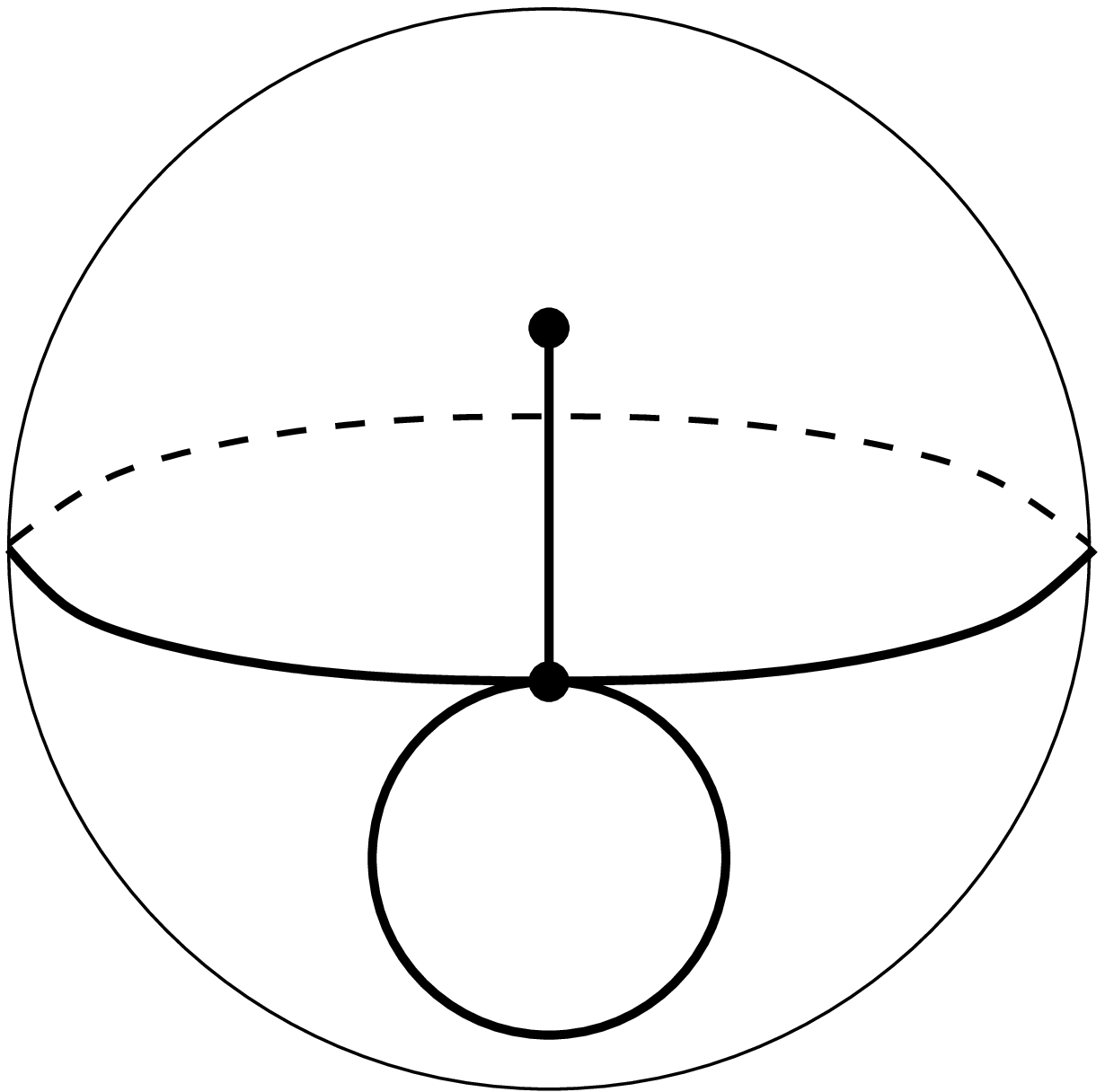,width=.85\linewidth}
\caption{$S(51|321)$. Valencies $(51|3,2,1)$. The order of the automorphism group: $1$.
Dual dessin   $S(321|51)$, see Figure \ref{three_v321_f51_sf}
on the page  \pageref{three_v321_f51_sf}.
Belyi function is $\beta = \frac{-64(15z+1)}{3125 z^3(z-1)^2(5z-8)}$.
} \label{three_v51_f321_sf}
\end{minipage}\hfill
\end{figure}
\begin{arrangedFigure}{1}{2}{}{
{$T(51|6)$. Valencies $(5,1|6)$. The order of the automorphism group: $1$. Dual dessin
$T(6|51)$, see Figure \ref{three_v6_f51_tor} on the page \pageref{three_v6_f51_tor}.
Belyi function is $\beta=
\frac{25\,{x}^{3}-270\,yx-255\,{x}^{2}+216\,y+522\,x-216}{1250}$
  on the curve  $X:{y}^{2}=-{\frac {5}{18}}\,{x}^{3}+{\frac {29}{36}}\,{x}^{2}-\frac{7}{3}x+1$.
$n_0 = \frac{x^5 (12+x)}{2500}$, $n_1=\frac{(x^3+6x^2-18x+58)^2}{2500}$.
}\label{three_v51_f6_tor}}
\subFig{three_v51_f6_tor}%
\subFig{three_v51_f6_tor_ucover}%
\end{arrangedFigure}
\clearpage

\begin{align*}\langle\langle Tr(Z^4)Tr(Z^2)Tr^3((Z^{+})^{2})\rangle\rangle=\\
=4\cdot2\cdot3!\cdot2^3\left(N^3+\frac{1}{2}N\right).
\end{align*}

\begin{figure}[h]
\begin{minipage}[b]{.45\linewidth}
\centering\epsfig{figure=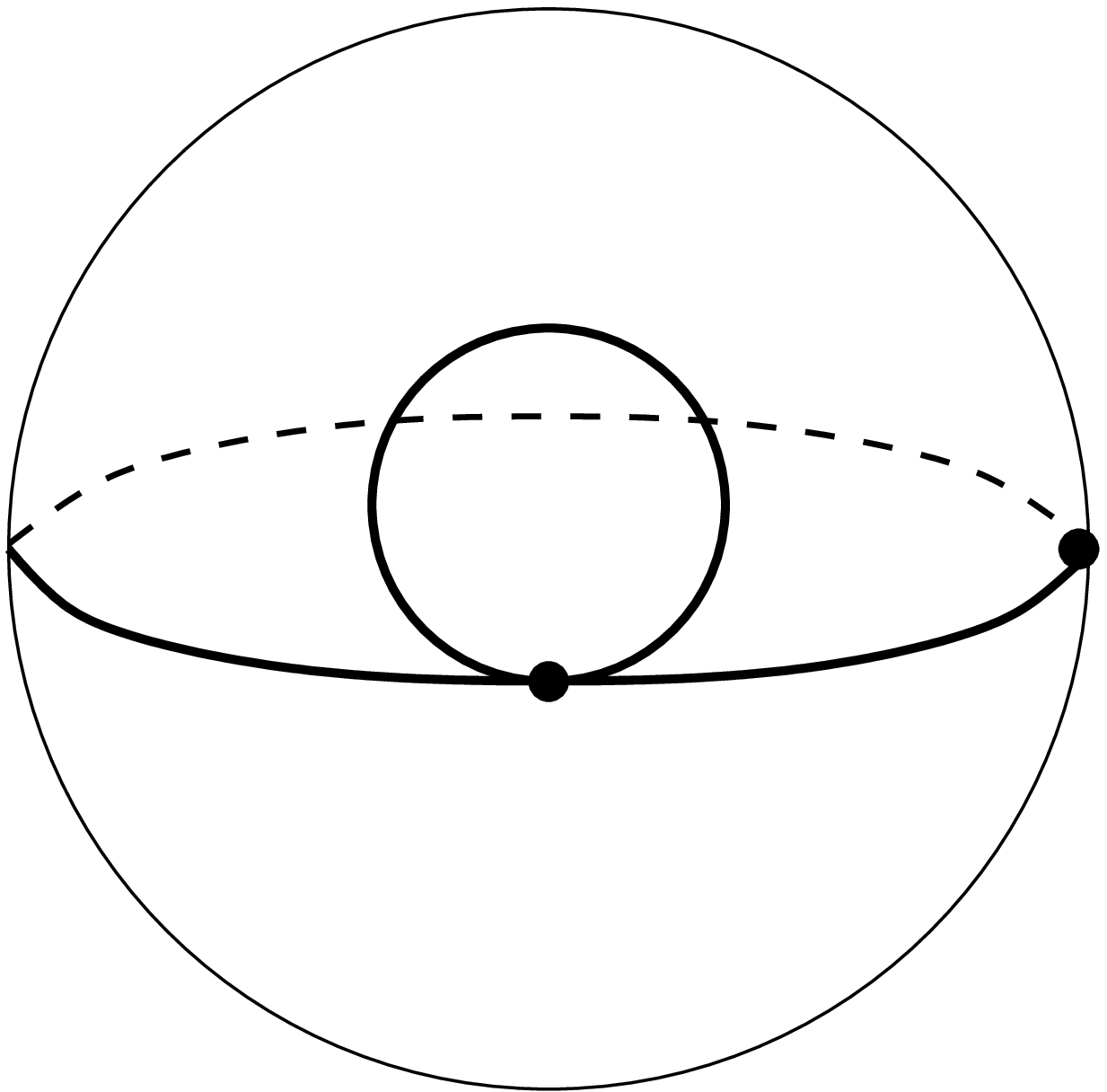,width=.85\linewidth}
\caption{$S(42|321)$. Valencies $(4,2|3,2,1)$. The order of the automorphism group: $1$.
Dual dessin   $S(321|42)$, see Figure \ref{three_v321_f42_sf}
on the page  \pageref{three_v321_f42_sf}.
Belyi function is $\beta = \frac{-(3z-2)^2}{4z^3(z-1)^2(z+2)}$.
} \label{three_v42_f321_sf}
\end{minipage}\hfill
\end{figure}
\begin{arrangedFigure}{1}{2}{}{
{$T(42|6)$. Valencies $(4,2|6)$. The order of the automorphism group: $2$. Dual dessin
$T(6|42)$, see Figure \ref{three_v6_f42_tor} on the page \pageref{three_v6_f42_tor}.
Belyi function is $(X:y^2=4x^3-39x+35, \beta = \frac{(x-1)^2(2x+7)}{27})$.
}\label{three_v42_f6_tor}}
\subFig{three_v42_f6_tor}%
\subFig{three_v42_f6_tor_ucover}%
\end{arrangedFigure}
\clearpage

\begin{align*}
\langle\langle Tr(Z^4)Tr(Z)Tr(Z)Tr^3((Z^{+})^{2})\rangle\rangle=\\
=2!\cdot 4\cdot3!\cdot2^3\left(\frac{3}{2}N^2\right)
=2!\cdot 4\cdot3!\cdot2^3\left(\left(1+\frac{1}{2}\right)N^2\right).
\end{align*}

\begin{figure}[h]
\begin{minipage}[b]{.45\linewidth}
\centering\epsfig{figure=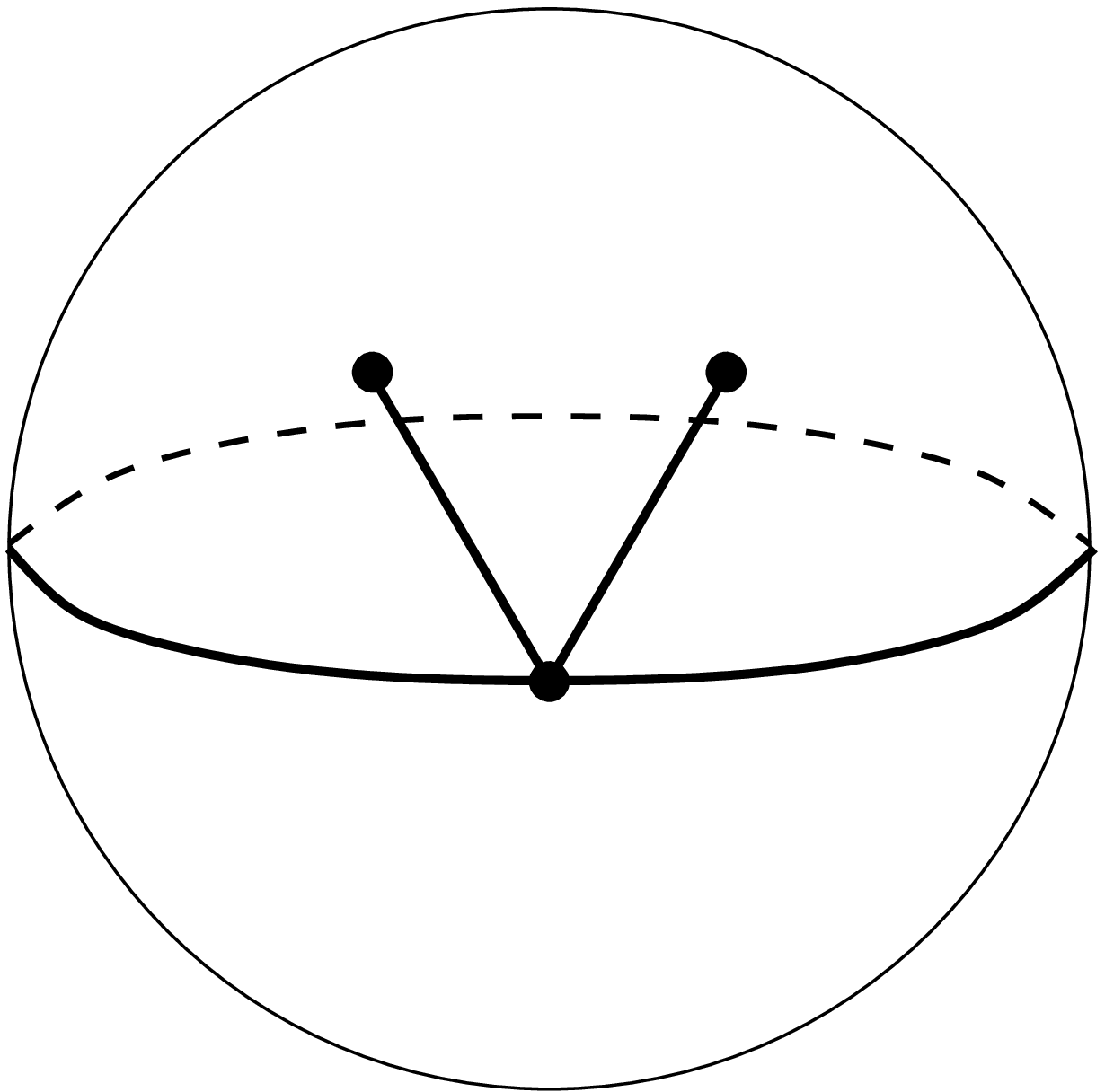,width=.85\linewidth}
\caption{$S(411|51)$. Valencies $(4,1,1|5,1)$. The order of the automorphism group: $1$.
Dual dessin   $S(51|411)$, see Figure \ref{three_v51_f411_sf}
on the page  \pageref{three_v51_f411_sf}.
Belyi function is $\beta = \frac{z^4 (z^2+4z+20)}{256(z-1)}$
} \label{three_v411_f51_sf}
\end{minipage}\hfill
\begin{minipage}[b]{.45\linewidth}
\centering\epsfig{figure=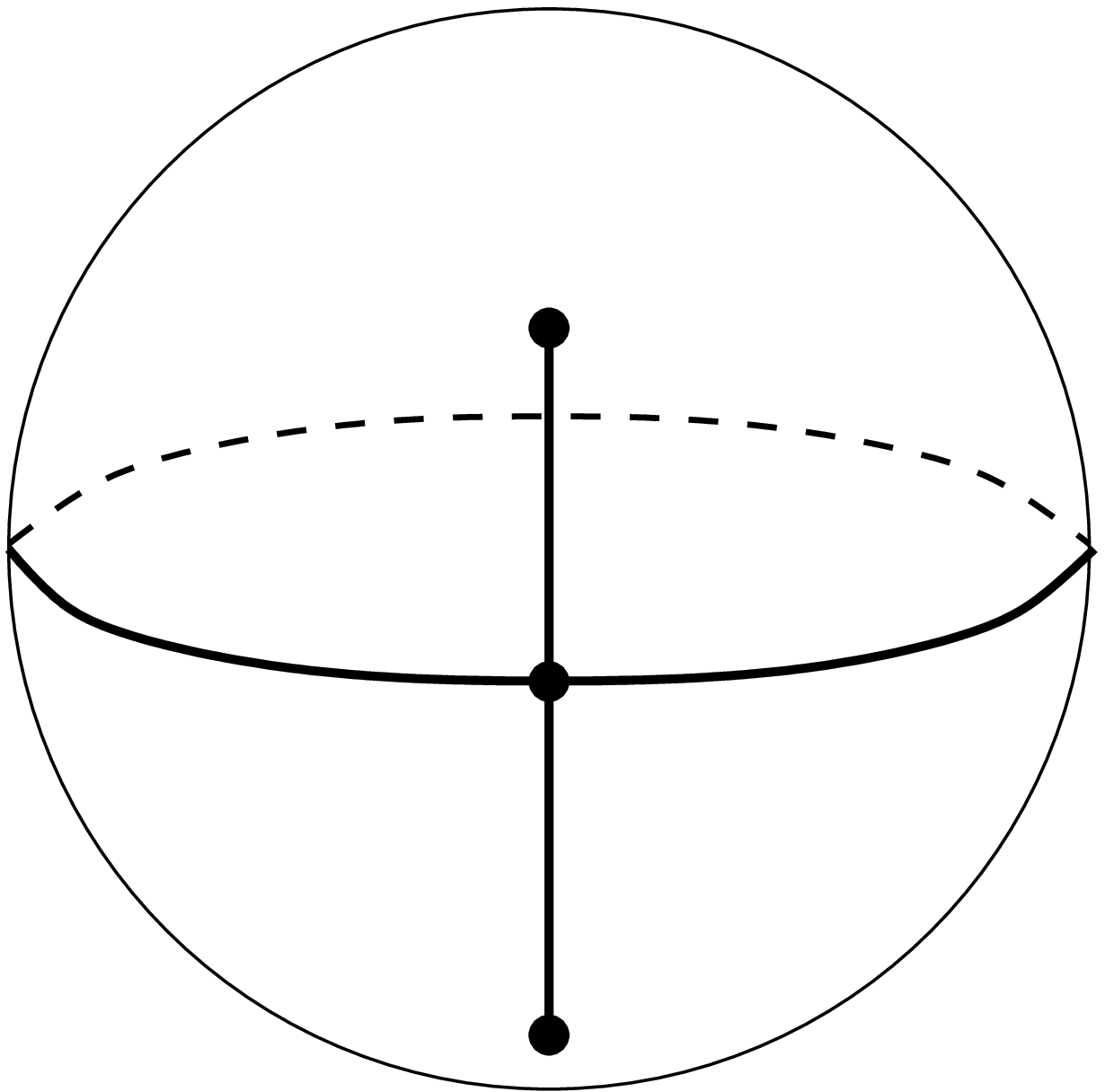,width=.85\linewidth}
\caption{$S(411|33)$. Valencies $(4,1,1|3,3)$. The order of the automorphism group: $2$.
Dual dessin   $S(33|411)$, see Figure \ref{three_v33_f411_sf}
on the page  \pageref{three_v33_f411_sf}.
Belyi function is $\beta=\frac{z^4(4z^2-3)}{4(z^2-1)^3}$.
} \label{three_v411_f33_sf}
\end{minipage}\hfill
\end{figure}

\clearpage
\begin{align*}\langle Tr(Z^3)Tr(Z^3)Tr^3((Z^{+})^{2})\rangle=
\langle\langle Tr(Z^3)Tr(Z^3)Tr^3((Z^{+})^{2})\rangle\rangle=\\
=2!\cdot3^2\cdot3!\cdot2^3\left(\frac{2}{3}N^3+\frac{1}{6}N\right)
=2!\cdot3^2\cdot3!\cdot2^3\left(\left(\frac{1}{2}+\frac{1}{6}\right)N^3+\frac{1}{6}N\right).
\end{align*}

\begin{figure}[h]
\begin{minipage}[b]{.45\linewidth}
\centering\epsfig{figure=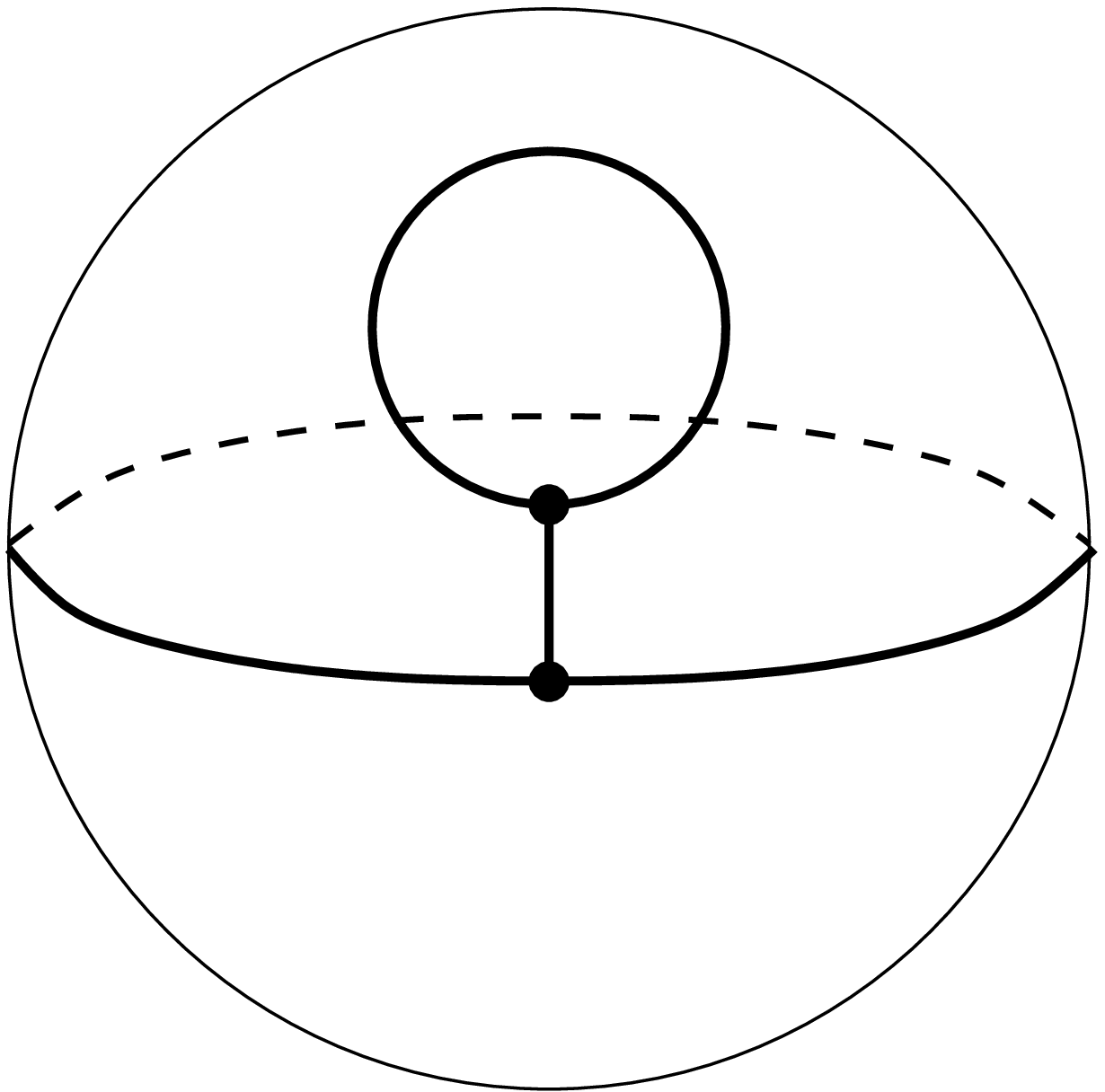,width=.85\linewidth}
\caption{$S(33|411)$. Valencies $(3,3|4,1,1)$. The order of the automorphism group: $2$.
Dual dessin   $S(411|33)$, see Figure \ref{three_v411_f33_sf}
on the page  \pageref{three_v411_f33_sf}.
Belyi function is $\beta=\frac{4(z^2-1)^3}{z^4(4z^2-3)}$.
} \label{three_v33_f411_sf}
\end{minipage}\hfill
\begin{minipage}[b]{.45\linewidth}
\centering\epsfig{figure=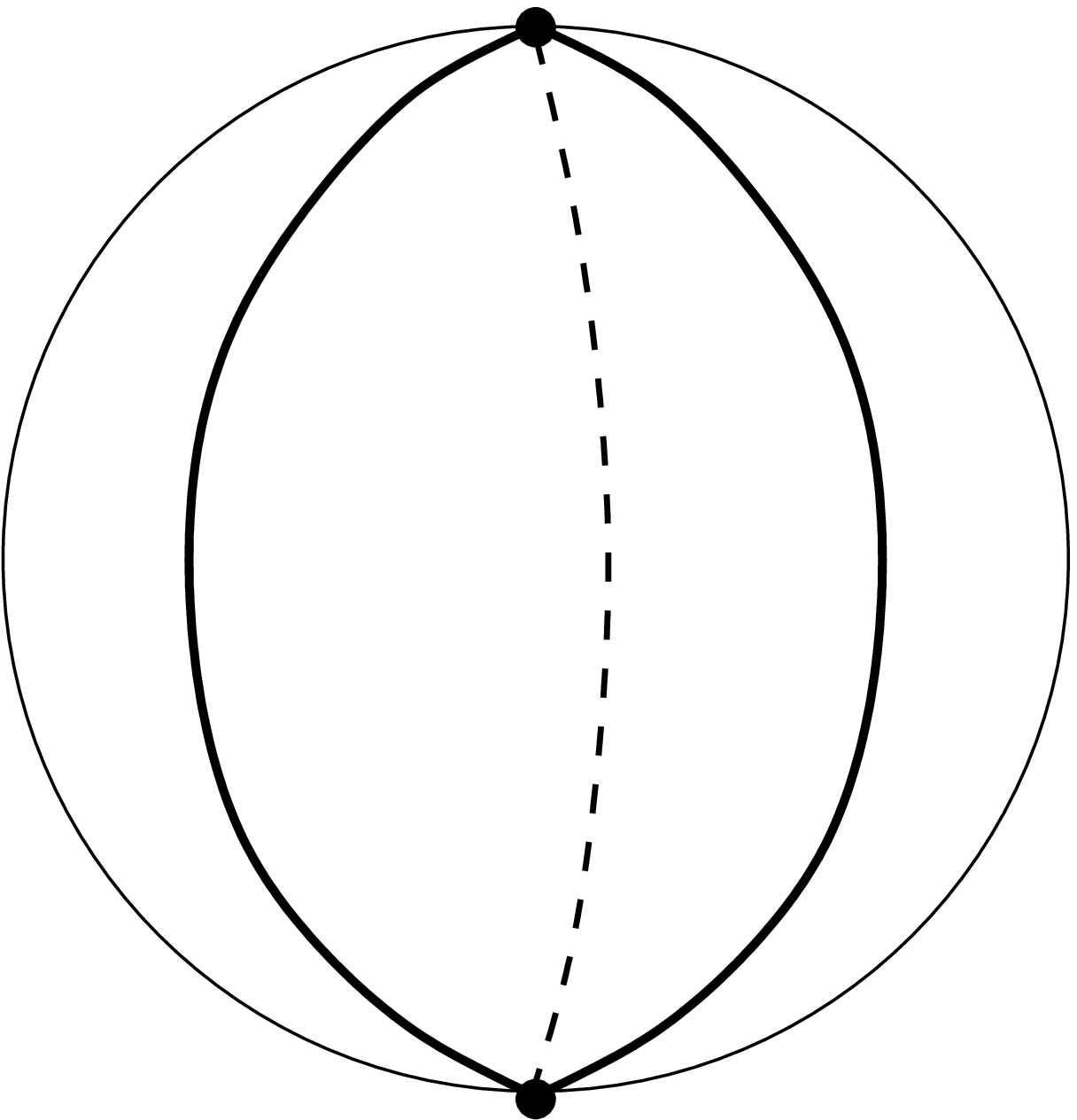,width=.85\linewidth}
\caption{$S(33|222)$. Valencies $(3,3|2,2,2)$. The order of the automorphism group: $6$.
Dual dessin   $S(222|33)$, see Figure \ref{three_v222_f33_sf}
on the page  \pageref{three_v222_f33_sf}.
Belyi function is $\beta = \frac{-4z^3}{(z^3-1)^2}$.
} \label{three_v33_f222_sf}
\end{minipage}\hfill
\end{figure}
\begin{arrangedFigure}{1}{2}{}{
{$T(33|6)$. Valencies $(3,3|6)$. The order of the automorphism group: $6$. Dual dessin
$T(6|33)$, see Figure \ref{three_v6_f33_tor} on the page \pageref{three_v6_f33_tor}.
Belyi function is $(X:y^2=x^3-1,\beta = x^3)$.
}
\label{three_v33_f6_tor}}
\subFig{three_v33_f6_tor}%
\subFig{three_v33_f6_tor_ucover}%
\end{arrangedFigure}

\clearpage
\begin{align*}
\langle\langle Tr(Z^3)Tr(Z^2)Tr(Z)Tr^3((Z^{+})^{2})\rangle\rangle=\\
=3\cdot2\cdot3!\cdot2^3\left(2N^2\right)
=3\cdot2\cdot3!\cdot2^3\left((1+1)N^2\right).
\end{align*}

\begin{figure}[h]
\begin{minipage}[b]{.45\linewidth}
\centering\epsfig{figure=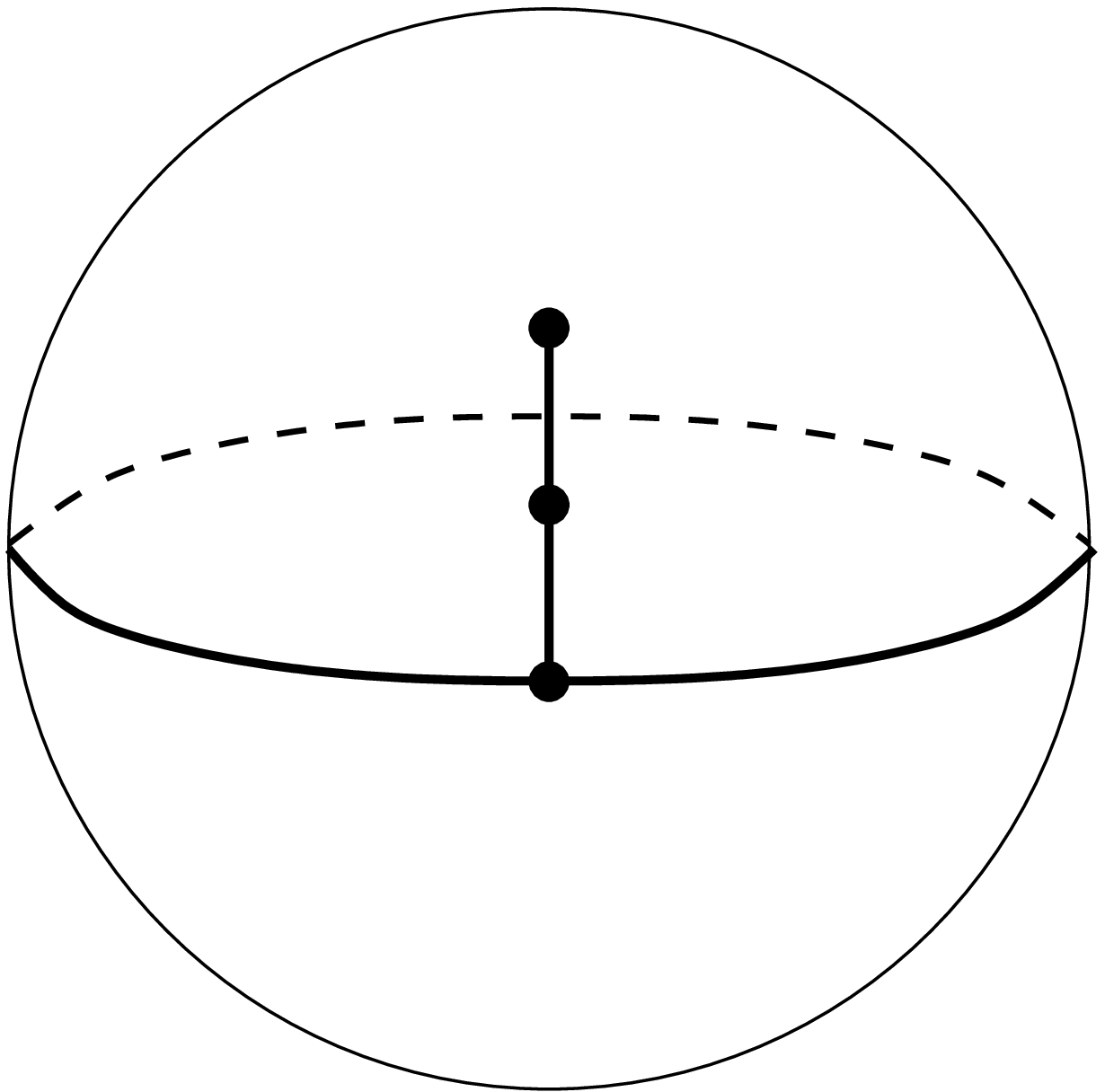,width=.85\linewidth}
\caption{$S(321|51)$. Valencies $(3,2,1|5,1)$. The order of the automorphism group: $1$.
Dual dessin   $S(51|321)$, see Figure \ref{three_v51_f321_sf}
on the page  \pageref{three_v51_f321_sf}.
Belyi function is $\beta = \frac{-z^3(z-5)^2(z-8)}{64(3z+1)}$
} \label{three_v321_f51_sf}
\end{minipage}\hfill
\begin{minipage}[b]{.45\linewidth}
\centering\epsfig{figure=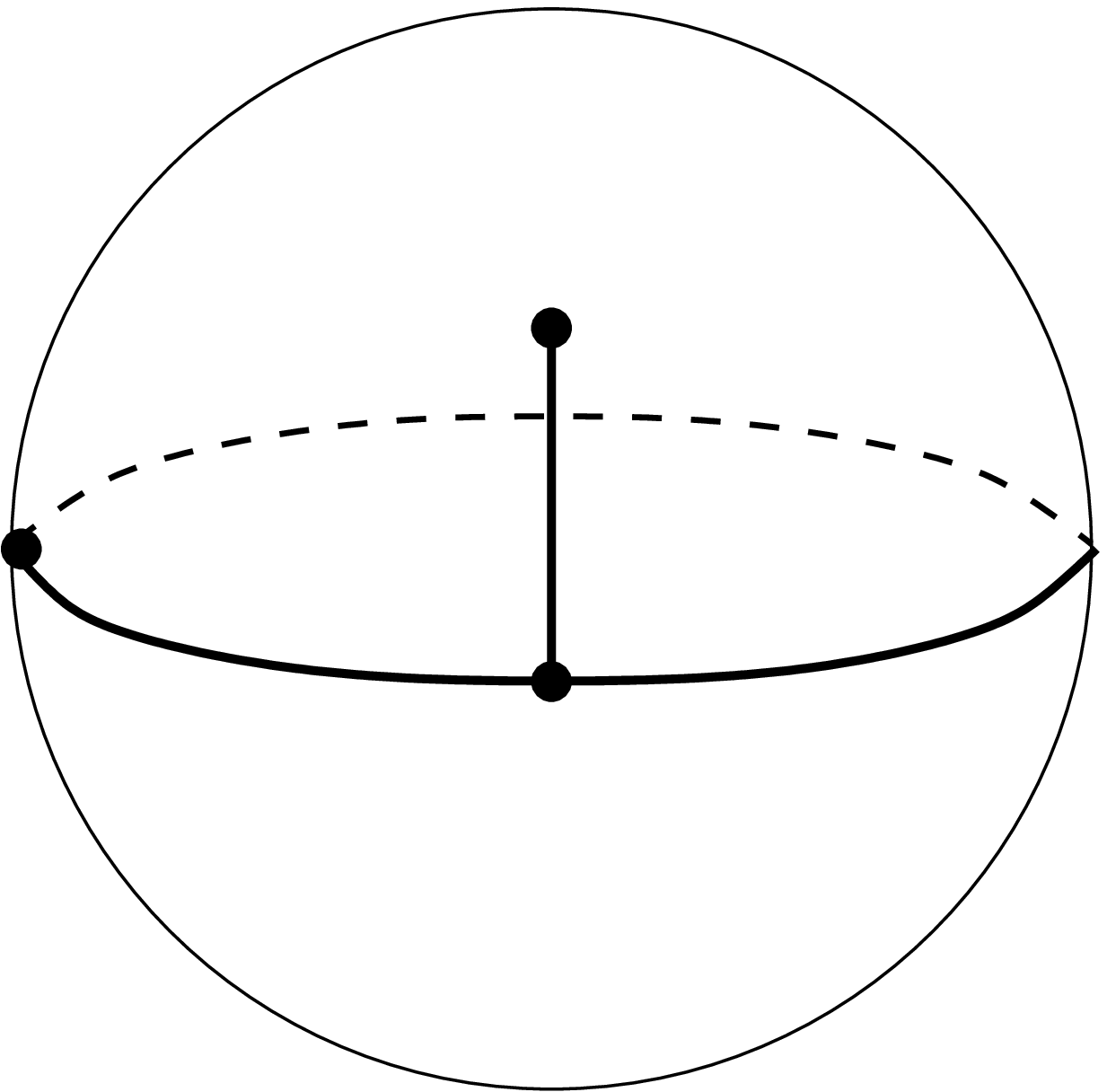,width=.85\linewidth}
\caption{$S(321|42)$. Valencies $(3,2,1|4,2)$. The order of the automorphism group: $1$.
Dual dessin   $S(42|321)$, see Figure \ref{three_v42_f321_sf}
on the page  \pageref{three_v42_f321_sf}.
Belyi function is $\beta = \frac{-4z^3(z-1)^2(z+2)}{(3z-2)^2}$.
} \label{three_v321_f42_sf}
\end{minipage}\hfill
\end{figure}

\clearpage
\begin{align*}
\langle\langle Tr(Z^3)Tr^3(Z)Tr^3((Z^{+})^{2})\rangle\rangle=
=3\cdot3!\cdot2^3\left(\frac{1}{3}N\right).
\end{align*}
\begin{figure}[h]
\begin{minipage}[b]{.45\linewidth}
\centering\epsfig{figure=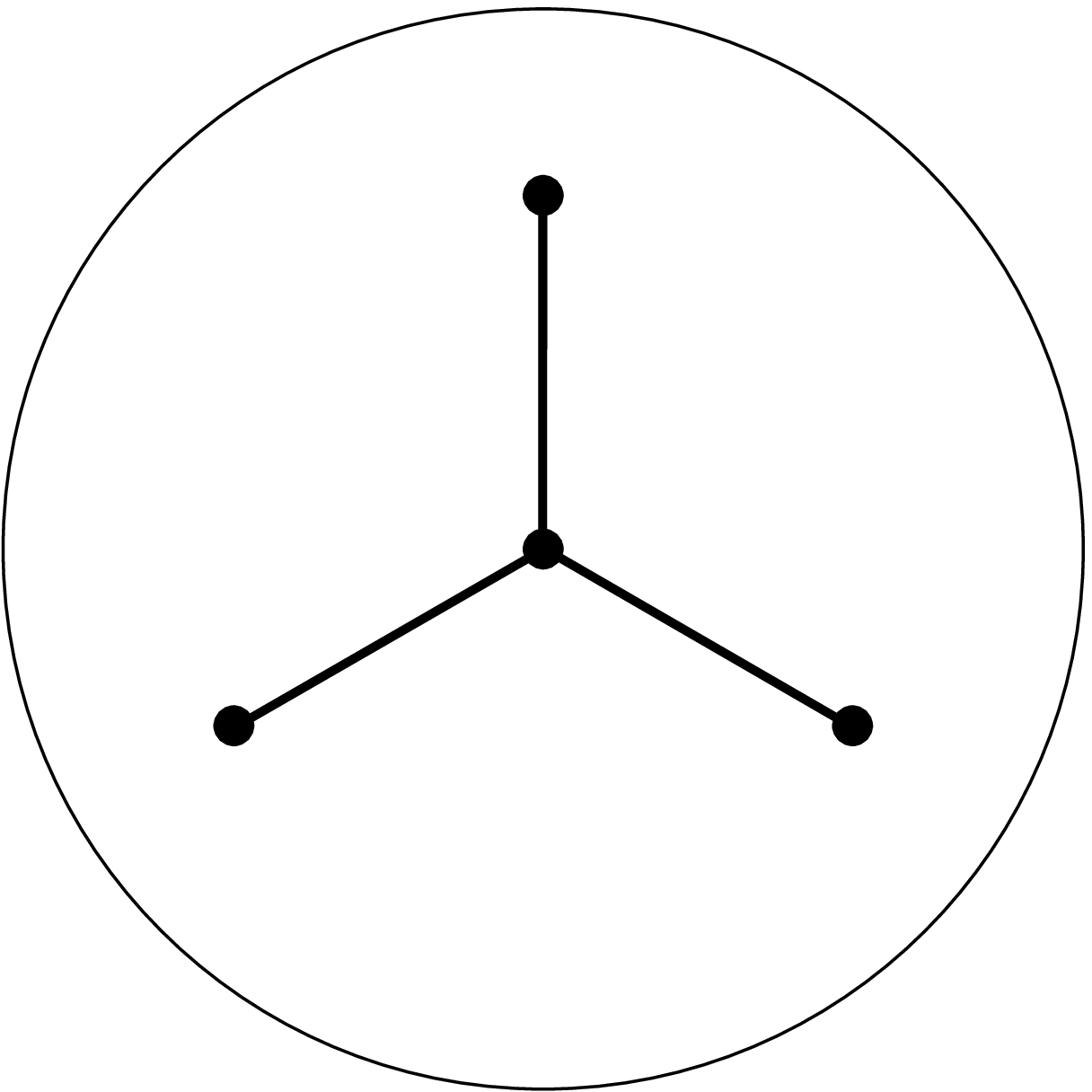,width=.85\linewidth}
\caption{$S(3111|6)$. Valencies $(3,1,1,1|6)$. The order of the automorphism group: $3$.
Dual dessin   $S(6|3111)$, see Figure \ref{three_v6_f3111_sf}
on the page  \pageref{three_v6_f3111_sf}.
Belyi function is $\beta=-4z^3(z^3-1)$.
} \label{three_v3111_f6_sf}
\end{minipage}\hfill
\end{figure}

\clearpage
\begin{align*}
\langle\langle Tr^3(Z^2)Tr^3((Z^{+})^{2})\rangle\rangle=
=3!\cdot2^3\cdot3!\cdot2^3\left(\frac{1}{6}N^2\right).
\end{align*}

\begin{figure}[h]
\begin{minipage}[b]{.45\linewidth}
\centering\epsfig{figure=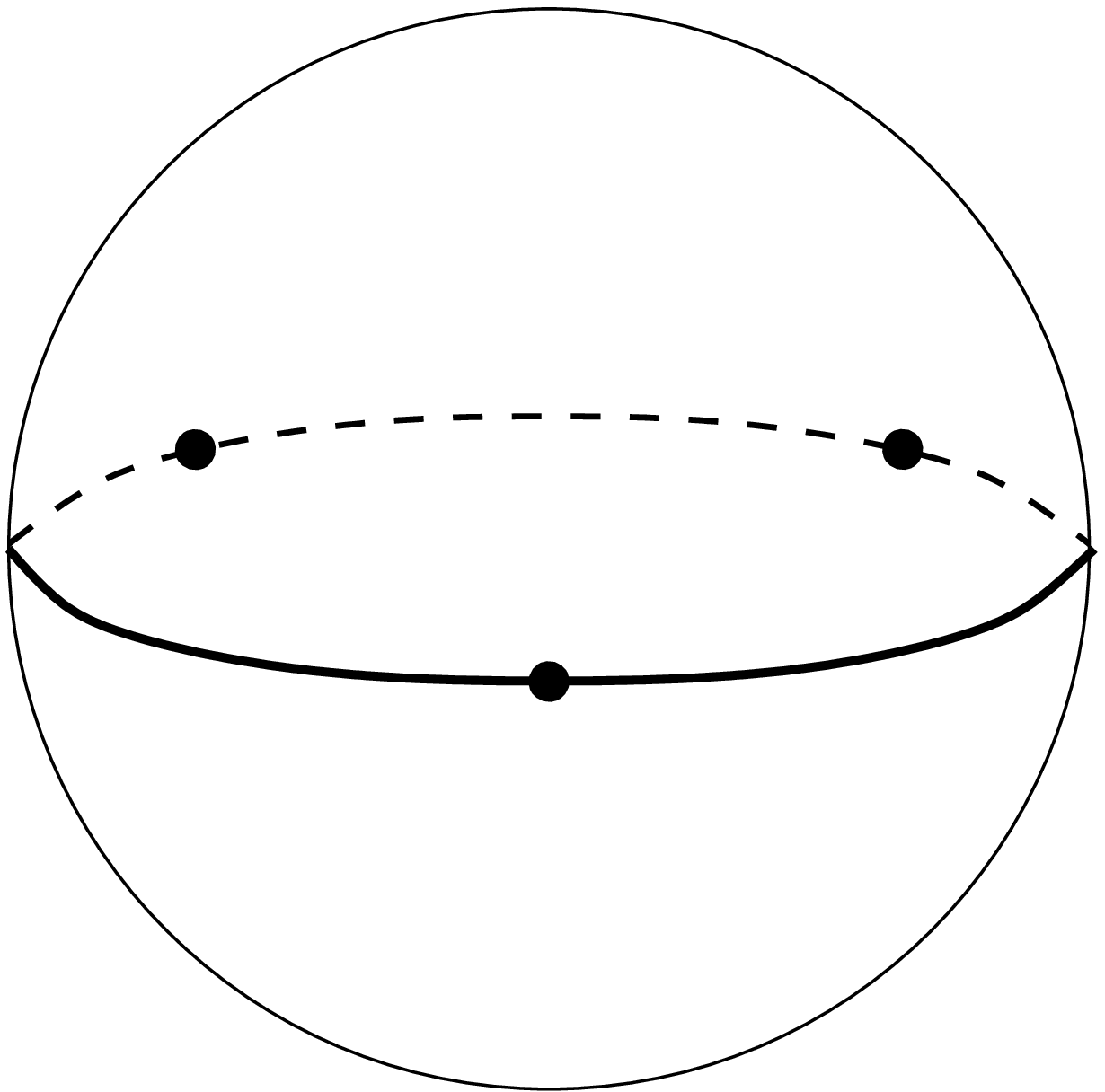,width=.85\linewidth}
\caption{$S(222|33)$. Valencies $(2,2,2|3,3)$. The order of the automorphism group: $6$.
Dual dessin   $S(33|222)$, see Figure \ref{three_v33_f222_sf}
on the page  \pageref{three_v33_f222_sf}.
Belyi function is $\beta = \frac{-(z^3-1)^2}{4z^3}$.
} \label{three_v222_f33_sf}
\end{minipage}\hfill
\end{figure}

\clearpage
\begin{align*}
\langle\langle Tr^2(Z^2)Tr^2(Z)Tr^3((Z^{+})^{2})\rangle\rangle=
=2!\cdot2^2\cdot2!\cdot3!\cdot2^3\left(\frac{1}{2}N\right).
\end{align*}

\begin{figure}[h]
\begin{minipage}[b]{.45\linewidth}
\centering\epsfig{figure=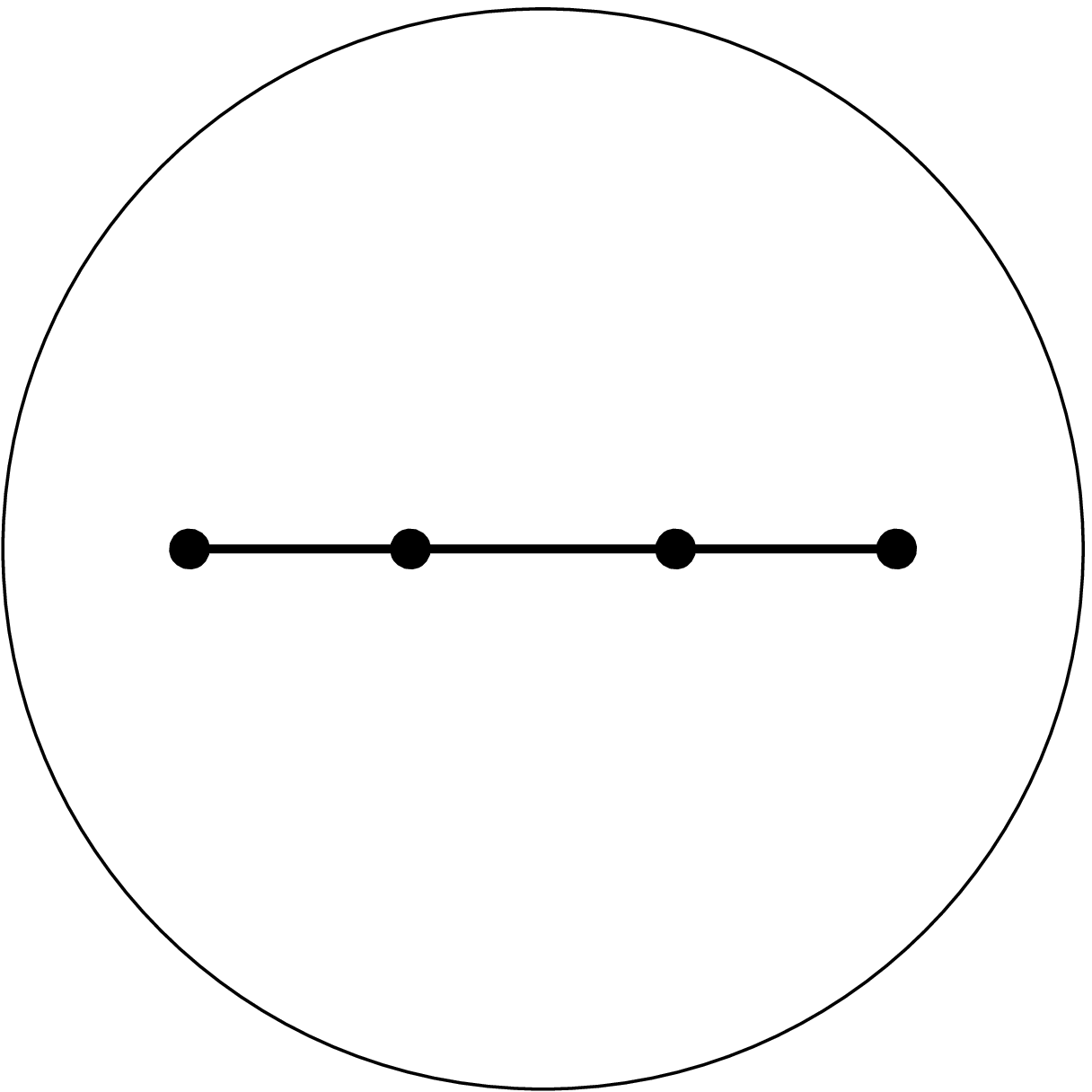,width=.85\linewidth}
\caption{$S(2211|6)$. Valencies $(2,2,1,1|6)$. The order of the automorphism group: $2$.
Dual dessin   $S(6|2211)$, see Figure \ref{three_v6_f2211_sf}
on the page  \pageref{three_v6_f2211_sf}.
Belyi function is $\beta=-(4z^2-1)^2(z^2-1)$.
} \label{three_v2211_f6_sf}
\end{minipage}\hfill
\end{figure}
\section{4-edges dessins}
\centerline{4-edges dessins and Hermitian matrix model:}
\begin{align*}
\langle\langle Tr(H^8)\rangle\rangle
=8\left(\frac{7}{4}N^5+\frac{35}{4}N^3+\frac{21}{8}N\right)=\\
=8\left(\left(1+\frac{1}{2}+\frac{1}{4}\right)N^5+\right(7+3\cdot\frac{1}{2}+\frac{1}{4}\left)N^3
+\left(2+\frac{1}{2}+\frac{1}{8}\right)N\right).
\end{align*}

\begin{arrangedFigure}{3}{1}{}{
$S(8|41111)$. Valencies $(8|4,1,1,1,1)$. The order of the automorphism group: $4$. Dual dessin
  $S(41111|8)$, see Figure \ref{four_v41111_f8_sf} on the page \pageref{four_v41111_f8_sf}.
Belyi function is $\beta=\frac{z^8}{4(z-1)(z+1)(z^2+1)}$
\label{four_v8_f41111_sf}
}
\subFig{four_v8_f41111_sf}%
\newSubFig{}{
{$S(8|32111)$. Valencies $(8|3,2,1,1,1)$. The order of the automorphism group: $1$. Dual dessin
  $S(32111|8)$, see Figure \ref{four_v32111_f8_sf} on the page \pageref{four_v32111_f8_sf}.
Belyi function is $\beta=\frac{729}{1024}\frac{z^8}{(z-1)(3z^2+8z+16)(3z-4)^2}$.
}
\label{four_v8_f32111_sf}
}
\subFig{four_v8_f32111_sf}%
\newSubFig{}{
{$S(8|22211)$. Valencies $(8|2,2,2,1,1)$. The order of the automorphism group: $2$. Dual dessin
  $S(22211|8)$, see Figure \ref{four_v22211_f8_sf} on the page \pageref{four_v22211_f8_sf}.
Belyi function is $\beta=\frac{1}{-4z^2(z^2-2)(z-1)^2(z+1)^2}$
}
\label{four_v8_f22211_sf}
}
\subFig{four_v8_f22211_sf}%
\end{arrangedFigure}

\begin{arrangedFigure}{2}{2}{}{
{$T(8|611)\_2$. Valencies $(8|6,1,1)$. The order of the automorphism group: $2$. Dual dessin
  $T(611|8)\_2$, see Figure \ref{four_v611_f8_tor_2} on the page \pageref{four_v611_f8_tor_2}.
Belyi function is $(X:y^2=3\,{x}^{3}+8\,{x}^{2}+16\,x,
\beta={\frac {27}{256}}\,{\frac {{x}^{4}}{x-1}})
$.
}
\label{four_v8_f611_tor_2}
}
\subFig{four_v8_f611_tor_2}%
\subFig{four_v8_f611_tor_2_ucover}%
\newSubFig{}{
{$T(8|611)\_1A$. Valencies $(8|6,1,1)$. The order of the automorphism group: $1$. Dual dessin
  $T(611|8)\_1A$, see Figure \ref{four_v611_f8_tor_1_a}
on the page  \pageref{four_v611_f8_tor_1_a}.} \label{four_v8_f611_tor_1_a}
}
\subFig{four_v8_f611_tor_1_a}%
\subFig{four_v8_f611_tor_1_a_ucover}%
\newSubFig{}{
{$T(8|611)\_1B$. Valencies $(8|6,1,1)$. The order of the automorphism group: $1$. Dual dessin
  $T(611|8)\_1B$, see Figure \ref{four_v611_f8_tor_1_b}
on the page  \pageref{four_v611_f8_tor_1_b}.} \label{four_v8_f611_tor_1_b}
}
\subFig{four_v8_f611_tor_1_b}%
\subFig{four_v8_f611_tor_1_b_ucover}%
\end{arrangedFigure}

\begin{figure}
For the dessin $T(8|611)\_1A$ Belyi function $\beta$ can be found from the equality
$\frac{1}{\beta}=-{\frac {1}{22769316864}}\, \left( -835-872\,\sqrt {2} \right)
( {x}^{4}-16\,y{x}^{2}+24\,
\sqrt {2}y{x}^{2}-648\,{x}^{3}+272\,\sqrt {2}{x}^{3}-4032\,\sqrt {2}yx
+4224\,yx-19296\,\sqrt {2}{x}^{2}+34200\,{x}^{2}+30240\,\sqrt {2}y-
50112\,y-818208\,x+706752\,\sqrt {2}x-1088640\,\sqrt {2}+1708560)
$
on the curve
$X:{y}^{2}=-{x}^{3}+20\,{x}^{2}-8\,\sqrt {2}{x}^{2}+1104\,\sqrt {2}x-132\,x$.
\\
For the dessin $T(8|611)\_1B$ Belyi function $\beta$ can be found from the equality
$\frac{1}{\beta}=-{\frac {1}{22769316864}}\, \left( -835+872\,\sqrt {2} \right)
( {x}^{4}-16\,y{x}^{2}-24\,\sqrt {2}y{x}^{2}
-648\,{x}^{3}-272\,\sqrt {2}{x}^{3}+4032\,\sqrt {2}yx+4224\,yx
+19296\,\sqrt {2}{x}^{2}+34200\,{x}^{2}-30240\,\sqrt {2}y
-50112\,y-818208\,x-706752\, \sqrt {2}x+1088640\,\sqrt {2}+1708560 )$
  on the curve
$X:{y}^{2}=-{x}^{3}+20\,{x}^{2}+8\,\sqrt {2}{x}^{2}-1104\,\sqrt {2}x-132\,x$.
\end{figure}

\begin{arrangedFigure}{2}{2}{}{
{$T(8|521)A_{+}$. Valencies $(8|5,2,1)$. The order of the automorphism group: $1$. Dual dessin
  $T(521|8)A_{+}$, see Figure \ref{four_v521_f8_tor_a_plus}
on the page  \pageref{four_v521_f8_tor_a_plus}.} \label{four_v8_f521_tor_a_plus}
}
\subFig{four_v8_f521_tor_a_plus}%
\subFig{four_v8_f521_tor_a_plus_ucover}%
\end{arrangedFigure}
\begin{arrangedFigure}{2}{2}{}{
{$T(8|521)A_{-}$. Valencies $(8|5,2,1)$. The order of the automorphism group: $1$. Dual dessin
  $T(521|8)A_{-}$, see Figure \ref{four_v521_f8_tor_a_minus}
on the page  \pageref{four_v521_f8_tor_a_minus}.} \label{four_v8_f521_tor_a_minus}
}
\subFig{four_v8_f521_tor_a_minus}%
\subFig{four_v8_f521_tor_a_minus_ucover}%
\newSubFig{}{
{$T(8|521)B$. Valencies $(8|5,2,1)$. The order of the automorphism group: $1$. Dual dessin
  $T(521|8)B$, see Figure \ref{four_v521_f8_tor_b} on the page \pageref{four_v521_f8_tor_b}.}
\label{four_v8_f521_tor_b}
}
\subFig{four_v8_f521_tor_b}%
\subFig{four_v8_f521_tor_b_ucover}%
\end{arrangedFigure}

\begin{figure}
For the dessins $T(8|521)A_{+}$, $T(8|521)A_{-}$ and $T(8|521)B$
Belyi function $\beta$ can be found from
$\frac{1}{\beta}={\frac {1}{735306250}}\,
\left( 552\,{\nu}^{2}-617\,\nu+68 \right)
( 42875\,{x}^{4}
+1756160\,\nu\,y{x}^{2}-860160\,{\nu}^{2}y{x}^{2}-4543840\,y{x}^{2}
+3959200\,\nu\,{x}^{3}-10346175\,{x}^{3}-1926400\,{\nu}^{2}{x}^{3}
+31782912\,{\nu}^{2}yx-63438592\,\nu\,yx+168996968\,yx
+18916352\,{\nu}^{2}{x}^{2}-37781632\,\nu\,{x}^{2}+100206428\,{x}^{2}
-257512128\,y-48381952\,{\nu}^{2}y+96684032\,\nu\,y
-62101504\,{\nu}^{2}x-330259656\,x+123960064\,\nu\,x
+48381952\,{\nu}^{2}-96684032\,\nu+257512128)
$
on the curve
$X:y^2=-{\frac {\left( 17\,{\nu}^{2}+8-42\,\nu \right)}{960400}}\,
(19600\,{x}^{2}+55552\,\nu\,x-18432\,{\nu}^{2}x-88408\,x+338963-130592
\,\nu+65792\,{\nu}^{2}) (x-1)$.\\
Here $256 \nu^3-544 \nu^2+1427 \nu-172=0$, and real $\nu$ corresponds to the case $B$.
\end{figure}

\begin{arrangedFigure}{2}{2}{}{
{$T(8|431)A$. Valencies $(8|4,3,1)$. The order of the automorphism group: $1$. Dual dessin
  $T(431|8)A$, see Figure \ref{four_v431_f8_tor_a} on the page \pageref{four_v431_f8_tor_a}.
Belyi function is
$\beta=-{\frac {85766121}{256 }}( 5488\,{x}^{4}+14112\,y{x}^{2}-26264
\,{x}^{3}+37548\,yx-202741\,{x}^{2}+3240\,y-73368\,x-3240 )^{-1}
$  on the curve
$X:{y}^{2}={\frac {1}{81}}\, ( 1-x ) ( 448\,{x}^{2}+1872\,x+81 )$.
$n_0={\frac {62523502209}{65536}}\,{\frac {1}{{x}^{4} \left( 4\,x
+45 \right)  \left( 4\,x+21 \right) ^{3}}}$,
$n_1={\frac { \left( 4096\,{x}^{4}+55296\,{x}^{3}+158976\,{x}^{2}+55296\,
x-247617 \right) ^{2}}{65536{x}^{4} \left( 4\,x+45 \right)
\left( 4\,x+21 \right) ^{3}}}$.
}
\label{four_v8_f431_tor_a}
}
\subFig{four_v8_f431_tor_a}%
\subFig{four_v8_f431_tor_a_ucover}%
\newSubFig{}{
{$T(8|431)B$. Valencies $(8|4,3,1)$. The order of the automorphism group: $1$. Dual dessin
  $T(431|8)B$, see Figure \ref{four_v431_f8_tor_b} on the page \pageref{four_v431_f8_tor_b}.
Belyi function is
$\beta=-36\, ( 81\,{x}^{4}-108\,y{x}^{2}+288\,{x}^{3}-96\,yx+308\,{x}^{2
}-24\,y+160\,x+32) ^{-1}$
 on the curve
$X:{y}^{2}=\frac49 \, \left( x+1 \right)  ( 9\,{x}^{2}+4\,x+4 )$.
$n_0=\frac{16}{3}{\frac {1}
{ \left( x-2 \right)  \left( 3\,x+2 \right)^{3}{x}^{4}}}$,
$n_1=\frac{1}{27} \frac { \left( 27\,{x}^{4}-36\,{x}^{2}-32\,x-20 \right) ^{2}}
{\left( x-2 \right)  \left( 3\,x+2 \right) ^{3}{x}^{4}}$.

}
\label{four_v8_f431_tor_b}
}
\subFig{four_v8_f431_tor_b}%
\subFig{four_v8_f431_tor_b_ucover}%
\end{arrangedFigure}

\begin{arrangedFigure}{2}{2}{}{
{$T(8|422)\_4$. Valencies $(8|4,2,2)$. The order of the automorphism group: $4$. Dual dessin
  $T(422|8)\_4$, see Figure \ref{four_v422_f8_tor_4} on the page \pageref{four_v422_f8_tor_4}.
Belyi function is $(X:{y}^{2}={x}^{3}-x,
\beta=\frac14\,{\frac {x^{4}}{ ( x-1)  (x+1) }})
$.
}
\label{four_v8_f422_tor_4}
}
\subFig{four_v8_f422_tor_4}%
\subFig{four_v8_f422_tor_4_ucover}%
\newSubFig{}{
{$T(8|422)\_2$. Valencies $(8|4,2,2)$. The order of the automorphism group: $2$. Dual dessin
  $T(422|8)\_2$, see Figure \ref{four_v422_f8_tor_2} on the page \pageref{four_v422_f8_tor_2}.
Belyi function is $X:y^2=4 x^3 -4x^2-x,
 \beta = \frac{1}{16 x^2 (x-1)^2}$.
}
\label{four_v8_f422_tor_2}
}
\subFig{four_v8_f422_tor_2}%
\subFig{four_v8_f422_tor_2_ucover}%
\newSubFig{}{
{$T(8|332)$. Valencies $(8|3,3,2)$. The order of the automorphism group: $2$. Dual dessin
  $T(332|8)$, see Figure \ref{four_v332_f8_tor} on the page \pageref{four_v332_f8_tor}.
Belyi function is $X:{y}^{2}=x \left( x-1 \right)  \left( 3\,{x}^{2}+8\,x+16 \right) ,
\beta={\frac {27}{256}}\,{\frac {{x}^{4}}{x-1}}
$.
}
\label{four_v8_f332_tor}
}
\subFig{four_v8_f332_tor}%
\subFig{four_v8_f332_tor_ucover}%
\end{arrangedFigure}

\begin{arrangedFigure}{2}{2}{}{
{$P(8|8)\_8$. Valencies $(8|8)$. The order of the automorphism group: $8$. Dual dessin
$P(8|8)\_8$, see Figure \ref{four_v8_f8_pr_8} on the page \pageref{four_v8_f8_pr_8}.
Belyi function is $(X:y^2=x^5-x,\beta=x^4)$.
}
\label{four_v8_f8_pr_8}
}
\subFig{four_v8_f8_pr_8}%
\subFig{four_v8_f8_pr_8_ucover}%
\newSubFig{}{
{$P(8|8)\_2$. Valencies $(8|8)$. The order of the automorphism group: $2$. Dual dessin
$P(8|8)\_2$, see Figure \ref{four_v8_f8_pr_2} on the page \pageref{four_v8_f8_pr_2}.
Belyi function is $\beta=-(-xy+x^4-2x^2+1)^2$   on the curve  $X:y^2=(x^2-2)(x^4-2x^2+2)$.
$n_0 = 1$, $n_1 = 4 (x-1)^4(x+1)^4$.
}
\label{four_v8_f8_pr_2}
}
\subFig{four_v8_f8_pr_2}%
\subFig{four_v8_f8_pr_2_ucover}%
\newSubFig{}{
{$P(8|8)\_1A$. Valencies $(8|8)$. The order of the automorphism group: $1$. Dual dessin
$P(8|8)\_1A$, see Figure \ref{four_v8_f8_pr_1a} on the page \pageref{four_v8_f8_pr_1a}.
}
\label{four_v8_f8_pr_1a}
}
\subFig{four_v8_f8_pr_1_a}%
\subFig{four_v8_f8_pr_1_a_ucover}%
\newSubFig{}{
{$P(8|8)\_1B$. Valencies $(8|8)$. The order of the automorphism group: $1$. Dual dessin
$P(8|8)\_1B$, see Figure \ref{four_v8_f8_pr_1b} on the page \pageref{four_v8_f8_pr_1b}.
}
\label{four_v8_f8_pr_1b}
}
\subFig{four_v8_f8_pr_1_b}%
\subFig{four_v8_f8_pr_1_b_ucover}%
\end{arrangedFigure}

\begin{figure}
Pairwise conjugate Belyi function of dessins $P(8|8)\_1A$ and $P(8|8)\_1B$:\\
$\beta=-1/8\, \left( -5+4\,\mu \right)  \left( x-1-\mu \right)  \left( x+1-
\mu \right)  \left( x+1 \right)  ( {x}^{2}-2\,x-2\,\mu\,x+1)
y+{\frac {1}{56}}\, \left( -5+4\,\mu \right)  ( 7\,{x}^{
8}-8\,{x}^{7}-40\,\mu\,{x}^{7}+140\,{x}^{6}-56\,\mu\,{x}^{5}+168\,{x}^
{5}-14\,{x}^{4}-224\,\mu\,{x}^{4}+168\,{x}^{3}-56\,\mu\,{x}^{3}+140\,{
x}^{2}-8\,x-40\,\mu\,x+7 )
$   on the curve
$X: \{ {y}^{2}={x}^{6}-{\frac {24}{7}}\,\mu\,{x}^{5}-2/7\,{x}^{5}+{
\frac {107}{49}}\,{x}^{4}-{\frac {200}{49}}\,\mu\,{x}^{4}+{\frac {500}
{49}}\,{x}^{3}-{\frac {48}{49}}\,\mu\,{x}^{3}-{\frac {200}{49}}\,\mu\,
{x}^{2}+{\frac {107}{49}}\,{x}^{2}-{\frac {24}{7}}\,\mu\,x-2/7\,x+1
 \} $. ${n_0}={x}^{8}$, ${n_1}=-{\frac {1}{196}}\, ( -9+4\,\mu)
( 7\,{x}^{4}+12\,{x}^{3}-4\,\mu\,{x}^{3}+6\,{x}^{2}-16
\,\mu\,{x}^{2}+12\,x-4\,\mu\,x+7 ) ^{2}$,
where $\mu=\pm \sqrt 2$.
\end{figure}

\clearpage
\begin{align*}
\langle\langle Tr(H^7)Tr(H)\rangle\rangle
=7\left(5N^4+10N^2\right).
\end{align*}
\vspace{-2cm}
\begin{arrangedFigure}{2}{1}{}{
{$S(71|5111)$. Valencies $(7,1|5,1,1,1)$. The order of the automorphism group: $1$. Dual dessin
  $S(5111|71)$, see Figure \ref{four_v5111_f71_sf} on the page \pageref{four_v5111_f71_sf}.
Belyi function is $\beta=\frac{16384z(z-1)^7}{896z^3-2912z^2+3216z-1225}$.
}
\label{four_v71_f5111_sf}}
\subFig{four_v71_f5111_sf}%
\end{arrangedFigure}
\begin{arrangedFigure}{2}{1}{}{
{$S(71|4211)_{+}$. Valencies $(7,1|4,2,1,1)$. The order of the automorphism group:
$1$. Dual dessin   $S(4211|71)_{+}$, see Figure \ref{four_v4211_f71_sf_plus}
on the page  \pageref{four_v4211_f71_sf_plus}.
Belyi function is
$\beta=\frac{-7340032(i \sqrt 7+21) z(z-1)^7}
{(896z^2-1904z+48 i z \sqrt 7+1029-49 i \sqrt 7 ) (112 z-119-5 i \sqrt 7 )^2}$
} \label{four_v71_f4211_sf_plus}}
\subFig{four_v71_f4211_sf_plus}%
\newSubFig{}{{$S(71|4211)_{-}$. Valencies $(7,1|4,2,1,1)$. The order of the automorphism group:
$1$. Dual dessin   $S(4211|71)_{-}$, see Figure \ref{four_v4211_f71_sf_minus}
on the page  \pageref{four_v4211_f71_sf_minus}.
Belyi function is
$\beta=\frac{-7340032(-i \sqrt 7+21) z(z-1)^7}
{(896z^2-1904z-48 i z \sqrt 7+1029+49 i \sqrt 7 ) (112 z-119+5 i \sqrt 7 )^2}$
} \label{four_v71_f4211_sf_minus}}
\subFig{four_v71_f4211_sf_minus}%
\end{arrangedFigure}

\begin{arrangedFigure}{2}{1}{}{
{$S(71|3311)$. Valencies $(7,1|3,3,1,1)$. The order of the automorphism group: $1$.
Dual dessin   $S(3311|71)$, see Figure \ref{four_v3311_f71_sf}
on the page  \pageref{four_v3311_f71_sf}.
Belyi function is
$\beta=-1728\,{\frac {z}{ \left( 1+{z}^{2}-5\,z \right) ^{3} \left( 49-13\,z+
{z}^{2} \right) }}$.

} \label{four_v71_f3311_sf}}
\subFig{four_v71_f3311_sf}%
\newSubFig{}{{$S(71|3221)$. Valencies $(7,1|3,2,2,1)$. The order of the automorphism group: $1$.
Dual dessin   $S(3221|71)$, see Figure \ref{four_v3221_f71_sf}
on the page  \pageref{four_v3221_f71_sf}.
Belyi function is
$\beta=-{\frac {1}{256}}\,{\frac {{z}^{7} \left( -48+z \right) }
{ \left( z+1 \right)  \left( 7\,{z}^{2}+28\,z+24 \right) ^{2}}}
$
} \label{four_v71_f3221_sf}}
\subFig{four_v71_f3221_sf}%
\end{arrangedFigure}

\begin{arrangedFigure}{2}{2}{}{
{$T(71|71)A_{+}$. Valencies $(7,1|7,1)$. The order of the automorphism group: $1$. Dual dessin
  $T(71|71)A_{-}$, see Figure \ref{four_v71_f71_tor_a_minus}
on the page  \pageref{four_v71_f71_tor_a_minus}.} \label{four_v71_f71_tor_a_plus}}
\subFig{four_v71_f71_tor_a_plus}%
\subFig{four_v71_f71_tor_a_plus_ucover}%
\newSubFig{}{{$T(71|71)A_{-}$. Valencies $(7,1|7,1)$. The order of the automorphism group: $1$.
Dual dessin   $T(71|71)A_{+}$, see Figure \ref{four_v71_f71_tor_a_plus}
on the page  \pageref{four_v71_f71_tor_a_plus}.} \label{four_v71_f71_tor_a_minus}}
\subFig{four_v71_f71_tor_a_minus}%
\subFig{four_v71_f71_tor_a_minus_ucover}%
\end{arrangedFigure}
\begin{figure}
-ы  Ёшёєэъют $T(71|71)A_{+}$ ш $T(71|71)A_{-}$ Belyi function is шьххЄ тшф
$\beta=-\frac { \left( \nu+3 \right)}{8(32\,x-51-7\,\nu)}
( -2\,\nu\,y{x}^{3}+10\,y{x}^{3}+56
\,{x}^{4}-64\,y{x}^{2}-14\,\nu\,{x}^{3}-266\,{x}^{3}+116\,yx+12\,\nu\,
yx+504\,{x}^{2}+56\,\nu\,{x}^{2}-11\,\nu\,y-73\,y-454\,x-82\,\nu\,x+
171+41\,\nu )
$
  on the curve  тшфр
$X:{y}^{2}=-\frac{1}{32} \left( 5+\nu \right)
(16\,{x}^{3}+ 4\,\nu\,{x}^{2}-52\,{x}^{2}-4\,\nu\,x+68\,x+\nu-37)$.
│фхё№ $\nu=\pm i\sqrt 7$.
\\
│Ёш ¤Єюь
$n_0=\frac { \left( 13+7\,\nu \right)  \left( x+3 \right)
 \left( x-1 \right) ^{7}}{8(32\,x-51-7\,\nu)}$,
$n_1=\frac { \left( 13+7\,\nu \right)  \left( 8\,{x}^{4}-16\,{x}^{3}-16
\,{x}^{2}+80\,x-59-7\,\nu \right) ^{2}}{512(32\,x-51-7\,\nu)}$.
\end{figure}
\begin{arrangedFigure}{2}{2}{}{
{$T(71|71)B_{+}$. Valencies $(7,1|7,1)$. The order of the automorphism group: $1$. Dual dessin
  $T(71|71)B_{+}$, see Figure \ref{four_v71_f71_tor_b_plus}
on the page  \pageref{four_v71_f71_tor_b_plus}. Belyi function is: ёь.эшцх.
} \label{four_v71_f71_tor_b_plus}}
\subFig{four_v71_f71_tor_b_plus}%
\subFig{four_v71_f71_tor_b_plus_ucover}%
\newSubFig{}{{$T(71|71)B_{-}$. Valencies $(7,1|7,1)$. The order of the automorphism group: $1$.
Dual dessin   $T(71|71)B_{-}$, see Figure \ref{four_v71_f71_tor_b_minus}
on the page  \pageref{four_v71_f71_tor_b_minus}. Belyi function can be seen below.
} \label{four_v71_f71_tor_b_minus}}
\subFig{four_v71_f71_tor_b_minus}%
\subFig{four_v71_f71_tor_b_minus_ucover}%
\end{arrangedFigure}

\begin{figure}
For the dessin $T(71|71)B_{+}$ Belyi function is
$\beta=-{\frac {343}{33554432}}\,{\frac { \left( -91\,i\sqrt {7}+87
 \right)  \left( 98\,{x}^{2}+21\,x+21\,i\sqrt {7}x+2 \right)  \left(
14\,{x}^{2}+7\,x+7\,i\sqrt {7}x-2 \right) ^{2}y}{x}}+
{\frac {1}{8388608x}}
\,
{( -91\,i\sqrt {7}+87 )} ( 1647086
\,{x}^{8}+3294172\,{x}^{7}+3294172\,i\sqrt {7}{x}^{7}-15764966\,{x}^{6
}+4941258\,i\sqrt {7}{x}^{6}-4705960\,i\sqrt {7}{x}^{5}-17882648\,{x}^
{5}-3882417\,i\sqrt {7}{x}^{4}+5260591\,{x}^{4}+672280\,i\sqrt {7}{x}^
{3}+2554664\,{x}^{3}-321734\,{x}^{2}+100842\,i\sqrt {7}{x}^{2}+2044\,i
\sqrt {7}x+1532\,x+686)
$
  on the curve
$X:-{x}^{4}-{\frac {11}{7}}\,i\sqrt {7}{x}^{3}-{\frac {11}{7}}\,{x}^{3}+
{\frac {519}{98}}\,{x}^{2}-{\frac {153}{98}}\,i\sqrt {7}{x}^{2}+
{\frac {51}{49}}\,i\sqrt {7}x+{\frac {869}{343}}\,x+{y}^{2}-1=0$.
$n_0=1$, $n_1=-{\frac {343}{33554432}}\,
{\frac { \left( -91\,i \sqrt {7}+87 \right)
\left( 14\,{x}^{2}+7\,x+7\,i\sqrt {7}x-2\right) ^{4}}{x}}$.
\\
For the dessin $T(71|71)B_{-}$ Belyi function is
$\beta=- \frac{343}{33554432}
\frac{
\left( 91\,i\sqrt {7}+87 \right)
\left( 98\,{x}^{2}+21\,x-21\,i\sqrt {7}x+2 \right)
\left( 14\,{x}^{2}+7\,x-7\,i\sqrt {7}x-2 \right)^{2}y
}{x}
+
\frac {1}{8388608x}
\left( 91\,i\sqrt {7}+87 \right)
( 1647086\,{x}^{8}+3294172\,{x}^{7}-3294172\,i\sqrt {7}{x}^{7}-15764966\,{x}^{6}-
4941258\,i\sqrt {7}{x}^{6}+4705960\,i\sqrt {7}{x}^{5}-17882648\,{x}^{5
}+3882417\,i\sqrt {7}{x}^{4}+5260591\,{x}^{4}-672280\,i\sqrt {7}{x}^{3
}+2554664\,{x}^{3}-321734\,{x}^{2}-100842\,i\sqrt {7}{x}^{2}-2044\,i
\sqrt {7}x+1532\,x+686
)
$
  on the curve
$X:-{x}^{4}+{\frac {11}{7}}\,i\sqrt {7}{x}^{3}-{\frac {11}{7}}\,{x}^{3}
+{\frac {519}{98}}\,{x}^{2}+{\frac {153}{98}}\,i\sqrt {7}{x}^{2}
-{\frac {51}{49}}\,i\sqrt {7}x+{\frac {869}{343}}\,x+{y}^{2}-1=0
$.
$n_0=1$,
$n_1={\frac {343}{33554432}}\,
{\frac { \left( 91\,i\sqrt {7}+87 \right)
\left( 14\,{x}^{2}+7\,x-7\,i\sqrt {7}x-2 \right) ^{4}}{x}}
$.
\end{figure}

\begin{arrangedFigure}{2}{2}{}{
{$T(71|71)C$. Valencies $(7,1,|7,1)$. The order of the automorphism group: $1$. Dual dessin
  $T(71|71)C$, see Figure \ref{four_v71_f71_tor_c} on the page \pageref{four_v71_f71_tor_c}.
Belyi function is
$\beta={\frac {1}{128}}\,{\frac {{x}^{3} \left( 7\,x+y \right) }{x-1}}$
  on the curve  $X:{y}^{2}=4{x}^{3}+13{x}^{2}+32x$.
$n_0=-{\frac {1}{4096}}\,{\frac {{x}^{7} \left( x-8 \right) }{x-1}}$,
$n_1=-{\frac {1}{4096}}\,{\frac { \left( {x}^{4}-4\,{x}^{3}-8\,
{x}^{2}-32\,x+64 \right) ^{2}}{x-1}}$.
}
\label{four_v71_f71_tor_c}}
\subFig{four_v71_f71_tor_c}%
\subFig{four_v71_f71_tor_c_ucover}%
\end{arrangedFigure}
\begin{arrangedFigure}{2}{2}{}{
{$T(71|62)_{+}$. Valencies $(7,1|6,2)$. The order of the automorphism group: $1$. Dual dessin
  $T(62|71)_{+}$, see Figure \ref{four_v62_f71_tor_plus}
on the page  \pageref{four_v62_f71_tor_plus}.
} \label{four_v71_f62_tor_plus}}
\subFig{four_v71_f62_tor_plus}%
\subFig{four_v71_f62_tor_plus_ucover}%
\newSubFig{}{{$T(71|62)_{-}$. Valencies $(7,1|6,2)$. The order of the automorphism group: $1$.
Dual dessin   $T(62|71)_{-}$, see Figure \ref{four_v62_f71_tor_minus}
on the page  \pageref{four_v62_f71_tor_minus}.} \label{four_v71_f62_tor_minus}}
\subFig{four_v71_f62_tor_minus}%
\subFig{four_v71_f62_tor_minus_ucover}%
\end{arrangedFigure}
\begin{figure}
For the dessins $T(71|62)_{+}$ and $T(71|62)_{-}$
Belyi function is
$\beta=-{\frac {2}{49( -112\,x-147+81\,\nu ) ^{2}}}\,(784\,{x}^{5}+19404\,\nu\,y{x}^{3}+12348\,y{x}^{3}+15435\,{x}^{4}+567\,\nu\,{x}^{4}+
54684\,\nu\,y{x}^{2}+188748\,y{x}^{2}+80892\,{x}^{3}-6804\,\nu\,{x}^{3
}+48384\,\nu\,yx+338688\,yx-40824\,\nu\,{x}^{2}+146664\,{x}^{2}+13608
\,\nu\,y+140616\,y+104976\,x-45360\,\nu\,x-13608\,\nu+25272)
$
  on the curve of the form
$X:{y}^{2}=-{\frac {\left( 35+9\,\nu \right)}{87808}}\,
( 2\,x+1 ) ( 14\,{x}^{2}+18\,\nu\,x+42\,x+9\,\nu+189)
$,
where $\nu=\pm i\sqrt 7$.\\
In this case
$n_0={\frac {4}{49}}\,{\frac { \left( x-24 \right) {x}^{7}}{
 \left( -112\,x-147+81\,\nu \right) ^{2}}}
$,
$n_1=\frac{1}{49}\,{\frac {
 \left( 2\,{x}^{4}-24\,{x}^{3}-144\,{x}^{2}-944\,x-1077+567\,\nu
 \right) ^{2}}{ \left( -112\,x-147+81\,\nu \right) ^{2}}}
$.
\end{figure}
\begin{arrangedFigure}{2}{2}{}{
{$T(71|53)A$. Valencies $(7,1|5,3)$. The order of the automorphism group: $1$.
Dual dessin   $T(53|71)A$, see Figure \ref{four_v53_f71_tor_a}
on the page  \pageref{four_v53_f71_tor_a}.} \label{four_v71_f53_tor_a}}
\subFig{four_v71_f53_tor_a}%
\subFig{four_v71_f53_tor_a_ucover}%
\newSubFig{}{{$T(71|53)B$. Valencies $(7,1|5,3)$. The order of the automorphism group: $1$.
Dual dessin   $T(53|71)B$, see Figure \ref{four_v53_f71_tor_b}
on the page  \pageref{four_v53_f71_tor_b}.} \label{four_v71_f53_tor_b}}
\subFig{four_v71_f53_tor_b}%
\subFig{four_v71_f53_tor_b_ucover}%
\end{arrangedFigure}

\begin{figure}
For the dessins $T(71|53)A$, $T(71|53)B$ Belyi function has the form
$\beta={\frac {\left( 78133+7625\,\nu \right)}{8232}}\,
( 169344\,{x}^{5}+17640\,y{x}^{3}-47250\,\nu\,{x}^{4}-489510\,{x}^{4}
-5250\,\nu\,y{x}^{2}-7350\,y{x}^{2}+1269324\,{x}^{3}+94500\,\nu\,{x}^{3}
+12250\,\nu\,yx-54950\,yx-3537716\,{x}^{2}+94500\,\nu\,{x}^{2}
+45885\,y-7125\,\nu\,y-301000\,\nu\,x+4558440\,x
-1978857+160125\,\nu)
$   on the curve of the form
$X: {y}^{2}={\frac {1}{350}}\,
( 2688\,{x}^{3}-10388\,{x}^{2}+100\,\nu\,{x}^{2}
-500\,\nu\,x+15972\,x+1275\,\nu-17247 )
( 12 x-13)
$. Here $\nu=\pm \sqrt{105}$
\\
and
$n_0={\frac {1944\left( 78133+7625\,\nu \right)}{343}}
\left( x+3 \right)  \left( x-1 \right) ^{7}$,
$n_1={\frac {78133+7625\,\nu}{131712}}\,
( 864\,{x}^{4}-1728\,{x}^{3}-1728\,{x}^{2}+5504\,x-11895+875\,\nu) ^{2}$.
\end{figure}

\begin{arrangedFigure}{2}{2}{}{
{$T(71|44)$. Valencies $(7,1|4,4)$. The order of the automorphism group: $1$.
Dual dessin   $T(44|71)$, see Figure \ref{four_v44_f71_tor}
on the page  \pageref{four_v44_f71_tor}.
Belyi function is $\beta={\frac {128}{343}} ( 294\,{x}^{4}-28\,y{x}^{2}+154\,{x}^{3}-7
\,yx+77\,{x}^{2}-2\,y+13\,x+2 ) $   on the curve
$X:y^2=112\,{x}^{4}+56\,{x}^{3}+37\,{x}^{2}+6\,x+1$.
$n_0=-{\frac {65536}{343}}\, \left( x-2 \right) {x}^{7}$,
$n_1=-{\frac {1}{343}}\, ( 256\,{x}^{4}-256\,{x}^{3}-128\,{x}^{2}-128
\,x-13 ) ^{2}$.
} \label{four_v71_f44_tor}}
\subFig{four_v71_f44_tor}%
\subFig{four_v71_f44_tor_ucover}%
\end{arrangedFigure}
\clearpage
\begin{align*}
\langle\langle Tr(H^6)Tr(H^2)\rangle\rangle
=6\cdot 2\left(\frac{5}{2}N^4+5N^2\right)=\\
=6\cdot
2\left(\left(2+\frac{1}{2}\right)N^4+\left(4+2\left(\frac{1}{2}\right)\right)N^2\right).\\
\end{align*}

\begin{arrangedFigure}{3}{1}{}{
{$S(62|4211)$. Valencies $(6,2|4,2,1,1)$. The order of the automorphism group: $1$. Dual dessin
  $S(4211|62)$, see Figure \ref{four_v4211_f62_sf} on the page \pageref{four_v4211_f62_sf}.
Belyi function is
$\beta=-{\frac {1}{108}}\,{\frac {{z}^{6} \left( z-4 \right) ^{2}}{ \left( {z
}^{2}+2\,z+3 \right)  \left( z-3 \right) ^{2}}}$
}
\label{four_v62_f4211_sf}}
\subFig{four_v62_f4211_sf}%
\newSubFig{}{{$S(62|3311)$. Valencies $(6,2|3,3,1,1)$. The order of the automorphism group: $2$.
Dual dessin   $(3311|62)$, see Figure \ref{four_v3311_f62_sf}
on the page  \pageref{four_v3311_f62_sf}.
Belyi function is $\beta=-64\,{\frac {{z}^{6}}{ \left( -1+3\,z \right)  \left( 3\,z+1 \right)
 \left( z-1 \right) ^{3} \left( z+1 \right) ^{3}}}
$.
} \label{four_v62_f3311_sf}}
\subFig{four_v62_f3311_sf}%
\newSubFig{}{{$S(62|3221)$. Valencies $(6,2|3,2,2,1)$. The order of the automorphism group: $1$.
Dual dessin   $S(3221|62)$, see Figure \ref{four_v3221_f62_sf}
on the page  \pageref{four_v3221_f62_sf}.
Belyi function is $\beta=-4\,{\frac {{z}^{6} \left( z-2 \right) ^{2}}{ \left( 4\,z+1 \right)
 \left( -1-2\,z+2\,{z}^{2} \right) ^{2}}}
$.
} \label{four_v62_f3221_sf}}
\subFig{four_v62_f3221_sf}%
\end{arrangedFigure}
\begin{arrangedFigure}{2}{2}{}{
{$T(62|71)_{+}$. Valencies $(6,2|7,1)$. The order of the automorphism group: $1$. Dual dessin
  $T(71|62)_{+}$, see Figure \ref{four_v71_f62_tor_plus}
on the page  \pageref{four_v71_f62_tor_plus}.} \label{four_v62_f71_tor_plus}}
\subFig{four_v62_f71_tor_plus}%
\subFig{four_v62_f71_tor_plus_ucover}%
\newSubFig{}{{$T(62|71)_{-}$. Valencies $(6,2|7,1)$. The order of the automorphism group: $1$.
Dual dessin   $T(71|62)_{-}$, see Figure \ref{four_v71_f62_tor_minus}
on the page  \pageref{four_v71_f62_tor_minus}.} \label{four_v62_f71_tor_minus}}
\subFig{four_v62_f71_tor_minus}%
\subFig{four_v62_f71_tor_minus_ucover}%
\end{arrangedFigure}
\begin{figure}
For the dessins $T(62|71)_{+}$ and $T(62|71)_{-}$ Belyi function has the form
$\beta={\frac{ \left( 11\,\nu+7 \right)}
{ 224\left( x-24 \right) {x}^{7}}}
( -686\,{x}^{5}+1078\,\nu\,{x}^{5}+197568\,y{x}^{3}+
20727\,\nu\,{x}^{4}-18963\,{x}^{4}+691488\,y{x}^{2}-211680\,\nu\,y{x}^{2}
+117180\,\nu\,{x}^{3}-5292\,{x}^{3}+762048\,yx-423360\,\nu\,yx
+237384\,\nu\,{x}^{2}+264600\,{x}^{2}-181440\,\nu\,y+254016\,y
+344736\,x+184032\,\nu\,x+46656\,\nu+108864)$
  on the curve of the form
$X:y^2=-{\frac {\left( 35+9\,\nu \right)}{87808}}\,
( 2\,x+1)  ( 14\,{x}^{2}+18\,\nu\,x+42\,x+9\,\nu+189)
$, where $\nu=\pm i \sqrt{7}$.
\\ For these dessins
$n_0=\frac{49}{4}\,{\frac { \left( -112\,x-147+81\,\nu
 \right) ^{2}}{ \left( x-24 \right) {x}^{7}}}$,\\
$n_1=\frac{1}{4}\,{\frac {
\left( 2\,{x}^{4}-24\,{x}^{3}-144\,{x}^{2}-944\,x-1077+567\,\nu
 \right) ^{2}}{ \left( x-24 \right) {x}^{7}}}$.
\end{figure}
\begin{arrangedFigure}{2}{2}{}{
{$T(62|62)\_2$. Valencies $(6,2|6,2)$. The order of the automorphism group: $2$.
Dual dessin   $T(62|62)\_2$, see Figure \ref{four_v62_f62_tor_2}
on the page  \pageref{four_v62_f62_tor_2}.
Belyi function is $(X:{y}^{2}=x \left( x+8 \right)  \left( x-1 \right),
\beta={\frac {1}{64}}\,{\frac {{x}^{3} \left( x+8 \right) }{x-1}})
$.
} \label{four_v62_f62_tor_2}}
\subFig{four_v62_f62_tor_2}%
\subFig{four_v62_f62_tor_2_ucover}%
\newSubFig{}{{$T(62|62)\_1$. Valencies $(6,2|6,2)$. The order of the automorphism group: $1$.
Dual dessin   $T(62|62)\_1$, see Figure \ref{four_v62_f62_tor_1}
on the page  \pageref{four_v62_f62_tor_1}.
Belyi function is
$\beta=-\frac {1}{432}\,{\frac { \left( 18\,{x}^{3}-810\,{x}^{2}+9342\,x-5382 \right) y}
{ \left( 4\,x-89 \right) ^{2}}}-
\frac {1}{432}\,\frac {4\,{x}^{5}-103\,{x}^{4}-1172\,{x}^{3}+28030\,{x}^{2}+
126536\,x-295247}{ \left( 4\,x-89 \right) ^{2}}$
  on the curve
$X:{y}^{2}=4\, ( x-2 ) ( {x}^{2}+2\,x+73)$
.
$n_0={\frac {1}{186624}}\,{\frac { \left( x-29 \right) ^{2}
 \left( x-5 \right) ^{6}}{ \left( 4\,x-89 \right) ^{2}}}
$,
$n_1={\frac {1}{186624}}\,{\frac { \left( {x}^{4}-44\,{x}^{3}+510\,{x}^{2}-
572\,x+35161 \right) ^{2}}{ \left( 4\,x-89 \right) ^{2}}}$.
} \label{four_v62_f62_tor_1}}
\subFig{four_v62_f62_tor_1}%
\subFig{four_v62_f62_tor_1_ucover}%
\newSubFig{}{{$T(62|53)$. Valencies $(6,2|5,3)$. The order of the automorphism group: $1$.
Dual dessin   $T(53|62)$, see Figure \ref{four_v53_f62_tor}
on the page  \pageref{four_v53_f62_tor}.
Belyi function is
$\beta=\frac{1}{162}({x}^{5}+y{x}^{3}+{\frac {65}{8}}\,{x}^{4}+6\,y{x}^{2}+{x}^{3}-64\,{x}^
{2}-16\,y-16\,x+80)$
  on the curve
$X:{y}^{2}=\frac14 ( x-1 )  ( 4x+5 )  ( x^2+4x-20)$.
$n_0={\frac {1}{20736}}\, \left( x+8 \right) ^{2}{x}^{6}$,
$n_1={\frac {1}{20736}}\, ( {x}^{4}+8\,{x}^{3}-128\,x-16) ^{2}$
} \label{four_v62_f53_tor}}
\subFig{four_v62_f53_tor}%
\subFig{four_v62_f53_tor_ucover}%
\newSubFig{}{{$T(62|44)$. Valencies $(6,2|4,4)$. The order of the automorphism group: $2$.
Dual dessin   $T(44|62)$, see Figure \ref{four_v44_f62_tor}
on the page  \pageref{four_v44_f62_tor}.
Belyi function is $(X:y^2=(x-2)(4x^2+4x+3)x,
\beta = -\frac{16}{27} x^3 (x-2) )
$.
} \label{four_v62_f44_tor}}
\subFig{four_v62_f44_tor}%
\subFig{four_v62_f44_tor_ucover}%
\end{arrangedFigure}
\clearpage

\begin{align*}
\langle\langle Tr(H^6)Tr^2(H)\rangle\rangle
=6\cdot 2!\left(5N^3+\frac{5}{2}N\right)=\\
=6\cdot
2!\left(\left(4+2\left(\frac{1}{2}\right)\right)N^3+\left(2+\frac{1}{2}\right)N\right).\\
\end{align*}

\begin{arrangedFigure}{3}{1}{}{
{$S(611|611)\_2$. Valencies $(6,1,1|6,1,1)$. The order of the automorphism group: $2$.
Dual dessin   $S(611|611)\_2$, see Figure \ref{four_v611_f611_sf_2}
on the page  \pageref{four_v611_f611_sf_2}.
Belyi function is
$\beta=-4\,{\frac { \left( {z}^{2}-2 \right) {z}^{6}}
{4\,{z}^{2}+1}}
$.
} \label{four_v611_f611_sf_2}}
\subFig{four_v611_f611_sf_2}%
\newSubFig{}{{$S(611|611)\_1$. Valencies $(6,1,1|6,1,1)$. The order of the automorphism group:
$1$. Dual dessin   $S(611|611)\_1$, see Figure \ref{four_v611_f611_sf_1}
on the page  \pageref{four_v611_f611_sf_1}.
Belyi function is
$\beta=-{\frac {27}{4}}\,{\frac { \left( 3\,{z}^{2}+6\,z+7 \right) {z}^{6}}
{21\,{z}^{2}-12\,z+4}}
$.
} \label{four_v611_f611_sf_1}}
\subFig{four_v611_f611_sf_1}%
\newSubFig{}{{$S(611|521)$. Valencies $(6,1,1|5,2,1)$. The order of the automorphism group: $1$.
Dual dessin   $S(521|611)$, see Figure \ref{four_v521_f611_sf}
on the page  \pageref{four_v521_f611_sf}.
Belyi function is
$\beta=4\,{\frac {{z}^{6} \left( 9\,{z}^{2}+24\,z+70 \right) }{ \left( 4\,z-1
 \right)  \left( 14\,z-5 \right) ^{2}}}
$.
} \label{four_v611_f521_sf}}
\subFig{four_v611_f521_sf}%
\end{arrangedFigure}
\begin{arrangedFigure}{2}{1}{}{
{$(S(611|431)_{+}$. Valencies $(6,1,1|4,3,1)$. The order of the automorphism group: $1$.
Dual dessin   $(S(431|611)_{+}$, see Figure \ref{four_v431_f611_sf_plus}
on the page  \pageref{four_v431_f611_sf_plus}.
Belyi function is
$\beta=12\,{\frac { \left( -6\,{z}^{2}+20\,i\sqrt {3}z-12\,z-19\,i\sqrt {3}+
17 \right) {z}^{6}}{ \left( -747+1763\,i\sqrt {3} \right)  \left( 3\,z
+2\,i\sqrt {3}+3 \right)  \left( z-1 \right) ^{3}}}
$.
} \label{four_v611_f431_sf_plus}}
\subFig{four_v611_f431_sf_plus}%
\newSubFig{}{{$(S(611|431)_{-}$. Valencies $(6,1,1|4,3,1)$. The order of the automorphism group:
$1$. Dual dessin   $(S(431|611)_{-}$, see Figure \ref{four_v431_f611_sf_minus}
on the page  \pageref{four_v431_f611_sf_minus}.
Belyi function is
$\beta=-12\,{\frac { \left( 6\,{z}^{2}+20\,i\sqrt {3}z+12\,z-19\,i\sqrt {3}-
17 \right) {z}^{6}}{ \left( 1763\,i\sqrt {3}+747 \right)  \left( -3\,z
-3+2\,i\sqrt {3} \right)  \left( z-1 \right) ^{3}}}
$.
} \label{four_v611_f431_sf_minus}}
\subFig{four_v611_f431_sf_minus}%
\newSubFig{}{{$S(611|332)$. Valencies $(6,1,1|3,3,2)$. The order of the automorphism group: $2$.
Dual dessin   $S(332|611)$, see Figure \ref{four_v332_f611_sf}
on the page  \pageref{four_v332_f611_sf}.
Belyi function is
$\beta=-4\,{\frac {{z}^{6} \left( 9\,{z}^{2}+2 \right) }
{ \left( 4\,{z}^{2}+1 \right) ^{3}}}
$.
} \label{four_v611_f332_sf}}
\subFig{four_v611_f332_sf}%
\end{arrangedFigure}
\begin{arrangedFigure}{2}{2}{}{
{$T(611|8)\_2$. Valencies $(6,1,1|8)$. The order of the automorphism group: $2$. Dual dessin
  $T(8|611)\_2$, see Figure \ref{four_v8_f611_tor_2} on the page \pageref{four_v8_f611_tor_2}.
Belyi function is $(X:{y}^{2}=x \left( {x}^{2}+1/2\,x+3/16 \right) ,
\beta=-{\frac {256}{27}}\,{x}^{3} \left( x-1 \right) )
$.
}
\label{four_v611_f8_tor_2}}
\subFig{four_v611_f8_tor_2}%
\subFig{four_v611_f8_tor_2_ucover}%
\newSubFig{}{{$T(611|8)\_1A$. Valencies $(6,1,1|8)$. The order of the automorphism group: $1$.
Dual dessin   $T(8|611)\_1A$, see Figure \ref{four_v8_f611_tor_1_a}
on the page  \pageref{four_v8_f611_tor_1_a}.} \label{four_v611_f8_tor_1_a}}
\subFig{four_v611_f8_tor_1_a}%
\subFig{four_v611_f8_tor_1_a_ucover}%
\newSubFig{}{{$T(611|8)\_1B$. Valencies $(6,1,1|8)$. The order of the automorphism group: $1$.
Dual dessin   $T(8|611)\_1B$, see Figure \ref{four_v8_f611_tor_1_b}
on the page  \pageref{four_v8_f611_tor_1_b}.} \label{four_v611_f8_tor_1_b}}
\subFig{four_v611_f8_tor_1_b}%
\subFig{four_v611_f8_tor_1_b_ucover}%
\end{arrangedFigure}
\begin{figure}
For the dessin $T(611|8)\_1A$ Belyi function is
$\beta=-{\frac {1}{22769316864}}\, \left( -835-872\,\sqrt {2} \right)
( {x}^{4}-16\,y{x}^{2}+24\,
\sqrt {2}y{x}^{2}-648\,{x}^{3}+272\,\sqrt {2}{x}^{3}-4032\,\sqrt {2}yx
+4224\,yx-19296\,\sqrt {2}{x}^{2}+34200\,{x}^{2}+30240\,\sqrt {2}y-
50112\,y-818208\,x+706752\,\sqrt {2}x-1088640\,\sqrt {2}+1708560)
$
  on the curve
$X:{y}^{2}=-{x}^{3}+20\,{x}^{2}-8\,\sqrt {2}{x}^{2}+1104\,\sqrt {2}x-132\,x$.
\\
For the dessin $T(611|8)\_1B$ Belyi function is
$\beta=-{\frac {1}{22769316864}}\, \left( -835+872\,\sqrt {2} \right)
( {x}^{4}-16\,y{x}^{2}-24\,\sqrt {2}y{x}^{2}
-648\,{x}^{3}-272\,\sqrt {2}{x}^{3}+4032\,\sqrt {2}yx+4224\,yx
+19296\,\sqrt {2}{x}^{2}+34200\,{x}^{2}-30240\,\sqrt {2}y
-50112\,y-818208\,x-706752\, \sqrt {2}x+1088640\,\sqrt {2}+1708560 )$
  on the curve
$X:{y}^{2}=-{x}^{3}+20\,{x}^{2}+8\,\sqrt {2}{x}^{2}-1104\,\sqrt {2}x-132\,x$.
\\
│Ёш ¤Єюь
$n_0=-{\frac {\left( -2217993-1456240\,\nu \right)}{518441790453234794496}}\,
( {x}^{2}+76\,x-152\,\nu\,x+44100-19600\,\nu)
\left( x+6-12\,\nu \right) ^{6}$,
$n_1=-{\frac {\left( -2217993-1456240\,\nu \right)}{518441790453234794496}}\,
 ( {x}^{4}+56\,{x}^{3}-112\,\nu\,{x}^{3}+22680\,{x}^{2}-10080\,
\nu\,{x}^{2}-756000\,x+665280\,\nu\,x+24794640-25197696\,\nu) ^{2}
$. │фхё№ $\nu^2=2$, $\nu>0$ т ёыєўрх $A$ ш $\nu<0$ т ёыєўрх $B$.
\end{figure}
\clearpage
\begin{align*}
\langle\langle Tr(H^5)Tr(H^3)\rangle\rangle
=5\cdot 3\left(3N^4+4N^2\right).\\
\end{align*}

\begin{arrangedFigure}{3}{1}{}{
{$S(53|5111)$. Valencies $(5,3|5,1,1,1)$. The order of the automorphism group: $1$. Dual dessin
  $S(5111|53)$, see Figure \ref{four_v5111_f53_sf} on the page \pageref{four_v5111_f53_sf}.
Belyi function is
$\beta=\frac14\,{\frac { \left( z-4 \right) ^{3}{z}^{5}}
{6\,{z}^{3}-22\,{z}^{2}-12\,z-9}}
$.
}
\label{four_v53_f5111_sf}}
\subFig{four_v53_f5111_sf}%
\newSubFig{}{{$S(53|4211)$. Valencies $(5,3|4,2,1,1)$. The order of the automorphism group: $1$.
Dual dessin   $S(4211|53)$, see Figure \ref{four_v4211_f53_sf}
on the page  \pageref{four_v4211_f53_sf}.
Belyi function is
$\beta=\frac14\,{\frac { \left( -1+8\,z \right) ^{3}}{ \left( 3+z \right) ^{2}{z}
^{4} \left( 9\,{z}^{2}+42\,z-5 \right) }}
$.
} \label{four_v53_f4211_sf}}
\subFig{four_v53_f4211_sf}%
\newSubFig{}{{$S(53|3221)$. Valencies $(5,3|3,2,2,1)$. The order of the automorphism group: $1$.
Dual dessin   $S(3221|53)$, see Figure \ref{four_v3221_f53_sf}
on the page  \pageref{four_v3221_f53_sf}.
Belyi function is
$\beta=4\,{\frac { \left( 3\,z+4 \right) ^{3}{z}^{5}}{ \left( 1+7\,z \right)
 \left( 3\,z+1 \right) ^{3} \left( {z}^{2}+1 \right) ^{2}}}
$.
} \label{four_v53_f3221_sf}}
\subFig{four_v53_f3221_sf}%
\end{arrangedFigure}
\begin{arrangedFigure}{2}{2}{}{
{$T(53|71)A$. Valencies $(5,3|7,1)$. The order of the automorphism group: $1$. Dual dessin
  $T(71|53)A$, see Figure \ref{four_v71_f53_tor_a} on the page \pageref{four_v71_f53_tor_a}.}
\label{four_v53_f71_tor_a}}
\subFig{four_v53_f71_tor_a}%
\subFig{four_v53_f71_tor_a_ucover}%
\newSubFig{}{{$T(53|71)B$. Valencies $(5,3|7,1)$. The order of the automorphism group: $1$.
Dual dessin   $T(71|53)B$, see Figure \ref{four_v71_f53_tor_b}
on the page  \pageref{four_v71_f53_tor_b}.} \label{four_v53_f71_tor_b}}
\subFig{four_v53_f71_tor_b}%
\subFig{four_v53_f71_tor_b_ucover}%
\end{arrangedFigure}
\begin{figure}
For the dessins $T(53|71)A$ and $T(53|71)B$ Belyi function has the form
$\beta=\frac { \left( 135\,\nu+1379 \right)}
{2730105000({64\,x-105+45\,\nu}}
( 297675\,\nu\,y{x}^{3}-4456305\,y{x}^{3}+49546350\,{x}^{4}
-201285\,\nu\,y{x}^{2}+25806879\,y{x}^{2}-151587135\,{x}^{3}-14586075\,\nu\,{x}^{3}
-8511615\,\nu\,yx+173009165\,yx+24310125\,\nu\,{x}^{2}-1233745275\,{x}^{2}
-54618075\,\nu\,y+59676225\,y-4281377625\,x+1348801875\,\nu\,x
+972759375\,\nu-25761054375)
$
  on the curve of the form
$X:{y}^{2}={\frac {225\left( 675\,\nu+5033 \right)}{45019072}}\,
( -420\,{x}^{3}+119\,{x}^{2}+405\,\nu\,{x}^{2}-21210\,x+630\,\nu\,x-176400+18900\,\nu)
$. Here $\nu^2=105$.\\
Now $n_0={\frac {343\left( 17983+1755\,\nu \right) }{15552000}}\,
{\frac {  \left( x+5 \right) ^{3} \left( x-3 \right) ^{5}}{64\,x-105+45\,\nu}}$,\\
$n_1={\frac {7(17983+1755\,\nu )}{62208000}}
\,{\frac {\left( 14\,{x}^{4}-420\,{x}^{2}+560\,x-60825+6075\,\nu
 \right) ^{2}}{64\,x-105+45\,\nu}}$.
\end{figure}
\begin{arrangedFigure}{2}{2}{}{
{$T(53|62)$. Valencies $(5,3|6,2)$. The order of the automorphism group: $1$.
Dual dessin   $T(62|53)$, see Figure \ref{four_v62_f53_tor}
on the page  \pageref{four_v62_f53_tor}.
Belyi function is $\beta=1296\, ( 8\,y{x}^{3}+48\,y{x}^{2}-128\,y+8\,{x}^{5}+65\,{x}^{4}+8
\,{x}^{3}-512\,{x}^{2}-128\,x+640 ) ^{-1}$
  on the curve
$X:{y}^{2}=\frac14 ( x-1 )  ( 4x+5 )  ( x^2+4x-20)$.
$n_0=20736\,{\frac {1}{ \left( x+8 \right) ^{2}{x}^{6}}}$,
$n_1={\frac { \left( {x}^{4}+8\,{x}^{3}-128\,x-16 \right) ^{2}}{ \left( x+
8 \right) ^{2}{x}^{6}}}$.
} \label{four_v53_f62_tor}}
\subFig{four_v53_f62_tor}%
\subFig{four_v53_f62_tor_ucover}%
\newSubFig{}{{$T(53|53)$. Valencies $(5,3|5,3)$. The order of the automorphism group: $1$.
Dual dessin   $T(53|53)$, see Figure \ref{four_v53_f53_tor}
on the page  \pageref{four_v53_f53_tor}.
Belyi function is $\beta=-{\frac {26}{25}}\,{x}^{2}+1/2\,{x}^{5}-{\frac {24}{25}}\,x-{\frac {8}{25}}\,y+{\frac {59}{50}}\,{x}^{3}+{
\frac {25}{16}}\,{x}^{4}+{\frac {8}{25}}\,yx+y{x}^{2}+1/2\,y{x}^{3}+{
\frac {8}{25}}$
  on the curve
$X:{y}^{2}=\frac14 \left( x-1 \right) ( 13\,{x}^{2}+12\,x+4\,{x}^{3}-4)
$.
$n_0={\frac {1}{6400}}\, \left( 9\,x+16 \right) ^{3}{x}^{5}$,
$n_1={\frac {1}{6400}}\, \left( 27\,{x}^{4}+72\,{x}^{3}+32\,{x}^{2}-128
\,x-48 \right) ^{2}$.
} \label{four_v53_f53_tor}}
\subFig{four_v53_f53_tor}%
\subFig{four_v53_f53_tor_ucover}%
\end{arrangedFigure}
\clearpage
\begin{align*}
\langle\langle Tr(H^5)Tr(H^2)Tr(H)\rangle\rangle
=5\cdot 2\left(6N^3+3N\right).\\
\end{align*}
\begin{arrangedFigure}{2}{1}{}{
{$S(521|611)$. Valencies $(5,2,1|6,1,1)$. The order of the automorphism group: $1$. Dual dessin
  $S(611|521)$, see Figure \ref{four_v611_f521_sf} on the page \pageref{four_v611_f521_sf}.
Belyi function is $\beta=-1/4\,{\frac { \left( 5\,z-14 \right) ^{2} \left( z-4 \right) {z}^{5}}
{70\,{z}^{2}+24\,z+9}}
$.
}
\label{four_v521_f611_sf}}
\subFig{four_v521_f611_sf}%
\newSubFig{}{{$S(521|521)A$. Valencies $(5,2,1|5,2,1)$.The order of the automorphism group: $1$.
Dual dessin   $S(521|521)A$, see Figure \ref{four_v521_f521_sf_a}
on the page  \pageref{four_v521_f521_sf_a}.
Belyi function is
$\beta=
 \{ 16\,{\frac { \left( 391+550\,\nu+455\,{\nu}^{2} \right)
 \left( z+2\,\nu \right)  \left( z+1 \right) ^{2}{z}^{5}}{ \left( 16\,
z-\nu+7\,{\nu}^{2}-4 \right)  \left( -8\,z+3\,\nu+3\,{\nu}^{2}-4
 \right) ^{2}}},7\,{\nu}^{3}+2\,{\nu}^{2}-\nu-4=0,\nu>0 \}
$.
} \label{four_v521_f521_sf_a}}
\subFig{four_v521_f521_sf_a}%
\end{arrangedFigure}
\begin{arrangedFigure}{2}{1}{}{
{$S(521|521)B_{+}$. Valencies $(5,2,1|5,2,1)$. The order of the automorphism group:
$1$. Dual dessin   $S(521|521)B_{+}$, see Figure \ref{four_v521_f521_sf_b_plus}
on the page  \pageref{four_v521_f521_sf_b_plus}.
Belyi function is
$\beta=
 \{ 16\,{\frac { \left( 391+550\,\nu+455\,{\nu}^{2} \right)
 \left( z+2\,\nu \right)  \left( z+1 \right) ^{2}{z}^{5}}{ \left( 16\,
z-\nu+7\,{\nu}^{2}-4 \right)  \left( -8\,z+3\,\nu+3\,{\nu}^{2}-4
 \right) ^{2}}},7\,{\nu}^{3}+2\,{\nu}^{2}-\nu-4=0, \mathop{Im} \nu<0 \}
$.
} \label{four_v521_f521_sf_b_plus}}
\subFig{four_v521_f521_sf_b_plus}%
\newSubFig{}{{$S(521|521)B_{-}$. Valencies $(5,2,1|5,2,1)$. The order of the automorphism group: $1$.
Dual dessin   $S(521|521)B_{-}$, see Figure \ref{four_v521_f521_sf_b_minus}
on the page  \pageref{four_v521_f521_sf_b_minus}.
Belyi function is
$\beta=
 \{ 16\,{\frac { \left( 391+550\,\nu+455\,{\nu}^{2} \right)
 \left( z+2\,\nu \right)  \left( z+1 \right) ^{2}{z}^{5}}{ \left( 16\,
z-\nu+7\,{\nu}^{2}-4 \right)  \left( -8\,z+3\,\nu+3\,{\nu}^{2}-4
 \right) ^{2}}},7\,{\nu}^{3}+2\,{\nu}^{2}-\nu-4=0, \mathop{Im} \nu>0 \}
$.
} \label{four_v521_f521_sf_b_minus}}
\subFig{four_v521_f521_sf_b_minus}%
\end{arrangedFigure}
\begin{arrangedFigure}{2}{1}{}{
{$S(521|431)$. Valencies $(5,2,1|4,3,1)$. The order of the automorphism group: $1$.
Dual dessin   $S(431|521)$, see Figure \ref{four_v431_f521_sf}
on the page  \pageref{four_v431_f521_sf}.
Belyi function is
$\beta=-16\,{\frac {{z}^{5} \left( z+3 \right)  \left( 6\,z-7 \right) ^{2}}{
 \left( 15\,z-4 \right)  \left( 7\,z-4 \right) ^{3}}}
$.
} \label{four_v521_f431_sf}}
\subFig{four_v521_f431_sf}%
\newSubFig{}{{$S(521|332)$. Valencies $(5,2,1|3,3,2)$. The order of the automorphism group: $1$.
Dual dessin   $S(332|521)$, see Figure \ref{four_v332_f521_sf}
on the page  \pageref{four_v332_f521_sf}.
Belyi function is
$\beta={\frac {27}{4}}\,{\frac {{z}^{5} \left( 7\,z+2 \right) ^{2} \left( 11
\,z-4 \right) }{ \left( 6\,{z}^{2}-1 \right) ^{3} \left( 4\,z+1
 \right) ^{2}}}
$.
} \label{four_v521_f332_sf}}
\subFig{four_v521_f332_sf}
\end{arrangedFigure}
\begin{arrangedFigure}{2}{2}{}{
{$T(521|8)A_{+}$. Valencies $(5,2,1|8)$. The order of the automorphism group: $1$. Dual dessin
  $T(8|521)A_{+}$, see Figure \ref{four_v8_f521_tor_a_plus}
on the page  \pageref{four_v8_f521_tor_a_plus}.} \label{four_v521_f8_tor_a_plus}}
\subFig{four_v521_f8_tor_a_plus}%
\subFig{four_v521_f8_tor_a_plus_ucover}%
\newSubFig{}{{$T(521|8)A_{-}$. Valencies $(5,2,1|8)$. The order of the automorphism group: $1$.
Dual dessin   $T(8|521)A_{-}$, see Figure \ref{four_v8_f521_tor_a_minus}
on the page  \pageref{four_v8_f521_tor_a_minus}.} \label{four_v521_f8_tor_a_minus}}
\subFig{four_v521_f8_tor_a_minus}%
\subFig{four_v521_f8_tor_a_minus_ucover}%
\newSubFig{}{{$T(521|8)B$. Valencies $(5,2,1|8)$. The order of the automorphism group: $1$.
Dual dessin   $T(8|521)B$, see Figure \ref{four_v8_f521_tor_b}
on the page  \pageref{four_v8_f521_tor_b}.} \label{four_v521_f8_tor_b}}
\subFig{four_v521_f8_tor_b}%
\subFig{four_v521_f8_tor_b_ucover}%
\end{arrangedFigure}
\begin{figure}
-ы  Ёшёєэъют $T(521|8)A_{+}$, $T(521|8)A_{-}$ ш $T(521|8)B$
Belyi function is шьххЄ тшф
$\beta={\frac {1}{735306250}}\,
\left( 552\,{\nu}^{2}-617\,\nu+68 \right)
( 42875\,{x}^{4}
+1756160\,\nu\,y{x}^{2}-860160\,{\nu}^{2}y{x}^{2}-4543840\,y{x}^{2}
+3959200\,\nu\,{x}^{3}-10346175\,{x}^{3}-1926400\,{\nu}^{2}{x}^{3}
+31782912\,{\nu}^{2}yx-63438592\,\nu\,yx+168996968\,yx
+18916352\,{\nu}^{2}{x}^{2}-37781632\,\nu\,{x}^{2}+100206428\,{x}^{2}
-257512128\,y-48381952\,{\nu}^{2}y+96684032\,\nu\,y
-62101504\,{\nu}^{2}x-330259656\,x+123960064\,\nu\,x
+48381952\,{\nu}^{2}-96684032\,\nu+257512128)
$
  on the curve  тшфр
$X:y^2=-{\frac {\left( 17\,{\nu}^{2}+8-42\,\nu \right)}{960400}}\,
(19600\,{x}^{2}+55552\,\nu\,x-18432\,{\nu}^{2}x-88408\,x+338963-130592
\,\nu+65792\,{\nu}^{2}) (x-1)$.\\
Here $256 \nu^3-544 \nu^2+1427 \nu-172=0$, and the real root corresponds to the case $B$.
\\
$n_0=-{\frac {\left( -1974439 \,\nu+8411384\,{\nu}^{2}+115196 \right)}{3460321800250000000}}\,
( 1225\,x+90376-34944\,\nu+16384\,{\nu}^{2})
( 1225\,x+56519-19936\,\nu+10496\,{\nu}^{2}) ^{2}
{x}^{5}$,
$n_1=-{\frac {\left( -1974439\,\nu+8411384\,{\nu}^{2}+115196 \right)}
{185122979184640000000}}\,
( 313600\,{x}^{4}+26036992\,{x}^{3}-9576448\,\nu
\,{x}^{3}+4784128\,{\nu}^{2}{x}^{3}+307426304\,{x}^{2}-114917376\,\nu
\,{x}^{2}+57409536\,{\nu}^{2}{x}^{2}-1834522624\,x+689504256\,\nu\,x-
344457216\,{\nu}^{2}x+6757769763-2539197472\,\nu+1270132992\,{\nu}^{2}
) ^{2}$.
\end{figure}
\clearpage
\begin{align*}
\langle\langle Tr(H^5)Tr^3(H)\rangle\rangle
=5\cdot 3!\left(2N^2\right).\\
\end{align*}
\begin{figure}[h]
\begin{minipage}[b]{.45\linewidth}
\centering\epsfig{figure=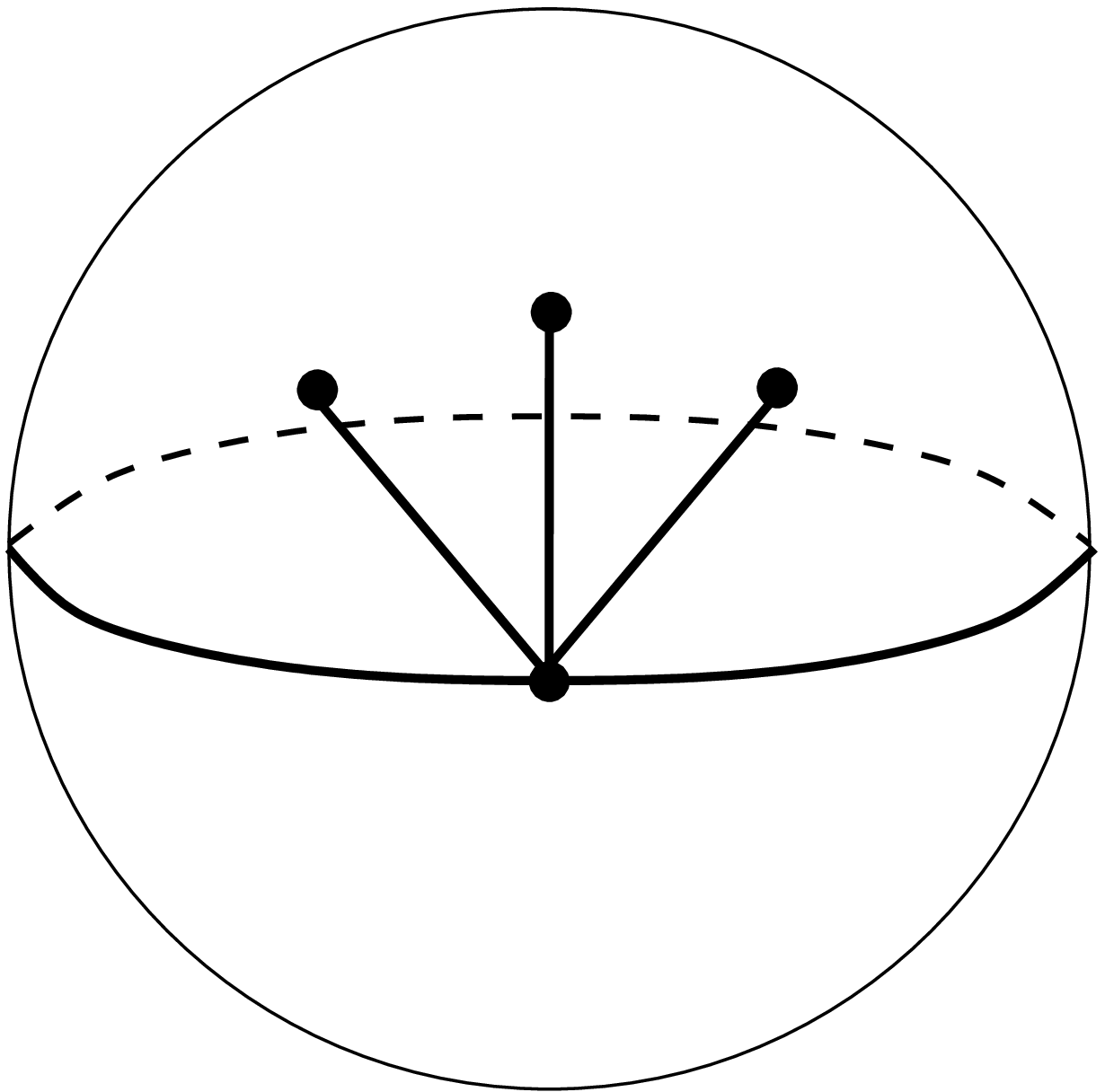,width=.85\linewidth}
\caption{$S(5111|71)$. Valencies $(5,1,1,1|7,1)$. The order of the automorphism group: $1$.
Dual dessin   $S(71|5111)$, see Figure \ref{four_v71_f5111_sf}
on the page  \pageref{four_v71_f5111_sf}.
Belyi function is $\beta=\frac{1}{16384} \frac{z^5(1225z^3-3216z^2+2912z-896)}{(z-1)^7}$.
} \label{four_v5111_f71_sf}
\end{minipage}\hfill
\begin{minipage}[b]{.45\linewidth}
\centering\epsfig{figure=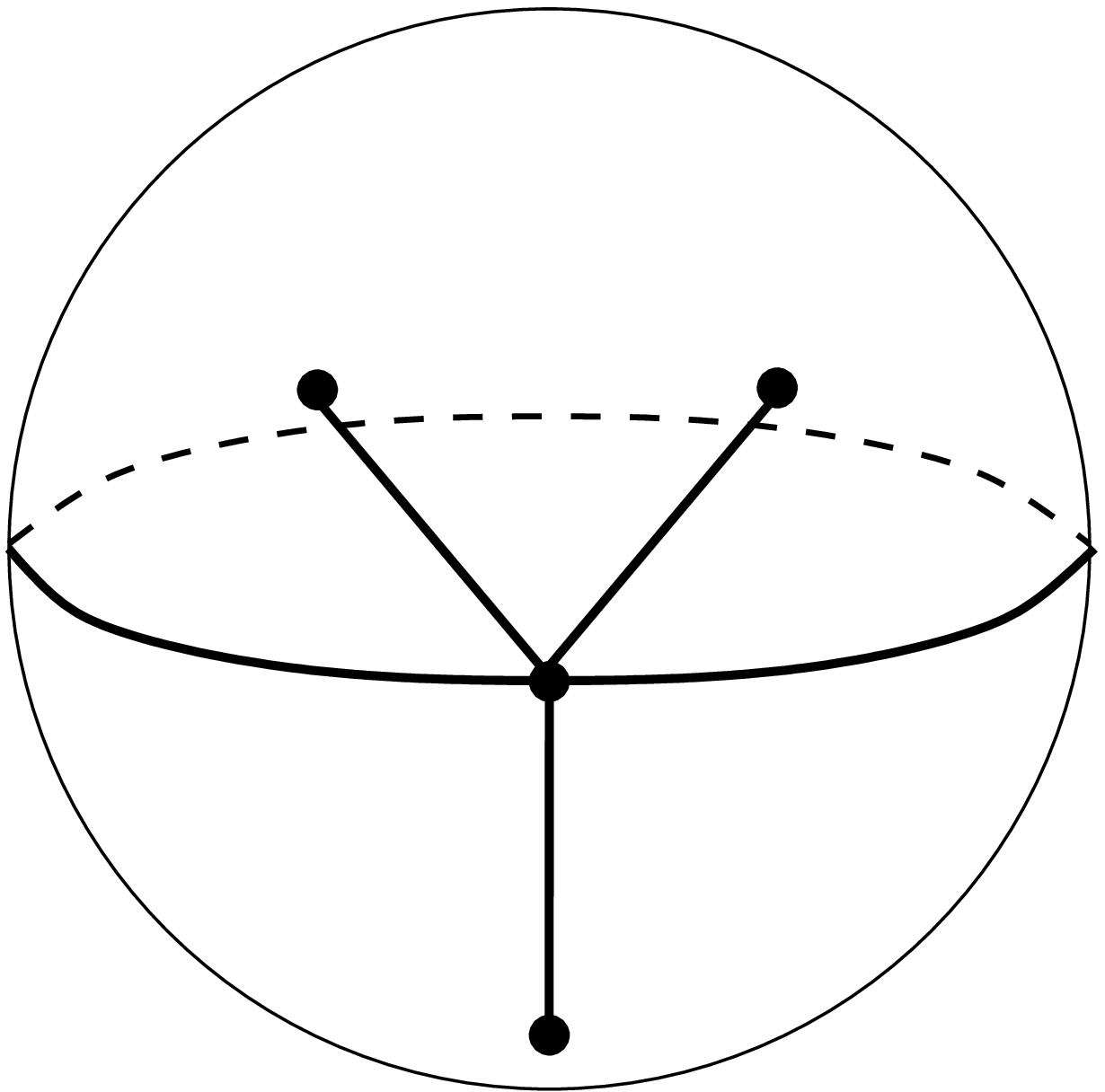,width=.85\linewidth}
\caption{$S(5111|53)$. Valencies $(5,1,1,1|5,3)$.The order of the automorphism group: $1$.
Dual dessin   $S(53|5111)$, see Figure \ref{four_v53_f5111_sf}
on the page  \pageref{four_v53_f5111_sf}.
Belyi function is
$\beta=4\,{\frac {{z}^{5} \left( 9\,{z}^{3}+12\,{z}^{2}+22\,z-6 \right) }
{\left( 4\,z-1 \right) ^{3}}}
$.
} \label{four_v5111_f53_sf}
\end{minipage}\hfill
\end{figure}
\clearpage
\begin{align*}
\langle\langle Tr^2(H^4)\rangle\rangle
=4^2\cdot 2!\left(\frac{9}{8}N^4+\frac{15}{8}N^2\right)=\\
=4^2\cdot 2!\left(\left(2\cdot\frac{1}{2}+\frac{1}{8}\right)N^4+
\left(1+\frac{1}{2}+\frac{1}{4}+\frac{1}{8}\right)N^2\right).\\
\end{align*}
\begin{arrangedFigure}{3}{1}{}{
{$S(44|4211)$. Valencies $(4,4|4,2,1,1)$. The order of the automorphism group: $2$. Dual dessin
  $S(4211|44)$, see Figure \ref{four_v4211_f44_sf} on the page \pageref{four_v4211_f44_sf}.
Belyi function is
$\beta=1/4\,{\frac { \left( {z}^{2}+1 \right) ^{4}}{{z}^{4} \left( 2\,{z}^{2}
+1 \right) }}
$.
}
\label{four_v44_f4211_sf}}
\subFig{four_v44_f4211_sf}%
\newSubFig{}{{$S(44|3311)$. Valencies $(4,4|3,3,1,1)$. The order of the automorphism group: $2$.
Dual dessin   $S(3311|44)$, see Figure \ref{four_v3311_f44_sf}
on the page  \pageref{four_v3311_f44_sf}.
Belyi function is
$\beta=-432\,{\frac {{z}^{4}}{ \left( 2\,{z}^{2}+10\,z-1 \right)  \left( 2\,{
z}^{2}+2\,z-1 \right) ^{3}}}
$.
} \label{four_v44_f3311_sf}}
\subFig{four_v44_f3311_sf}%
\newSubFig{}{{$S(44|2222)$. Valencies $(4,4|2,2,2,2)$. The order of the automorphism group: $8$.
Dual dessin   $S(2222|44)$, see Figure \ref{four_v2222_f44_sf}
on the page  \pageref{four_v2222_f44_sf}.
Belyi function is $\beta=4\frac{z^4}{(z^4+1)^2}$.
} \label{four_v44_f2222_sf}}
\subFig{four_v44_f2222_sf}%
\end{arrangedFigure}
\begin{arrangedFigure}{2}{2}{}{
{$T(44|71)$. Valencies $(4,4|7,1)$. The order of the automorphism group: $1$. Dual dessin
  $T(71|44)$, see Figure \ref{four_v71_f44_tor} on the page \pageref{four_v71_f44_tor}.
Belyi function is
$\beta=-{\frac {1}{512x}}
(-343\,y{x}^{3}+2401\,{x}^{4}-931\,y{x}^{2}+
9604\,{x}^{3}-581\,yx+11662\,{x}^{2}-y+3860\,x-7)$
  on the curve
$X:{y}^{2}=\left( 4\,x+7 \right)  ( 7\,{x}^{2}+18\,x+7)
$.
$n_0=-{\frac {343}{65536}}\,{\frac { \left( 7\,{x}^{2}+14\,x+3 \right) ^{4}}{x}}$,
$n_1=-{\frac {7}{65536}}\,
\frac { \left( 343\,{x}^{4}+1372\,{x}^{3}+1666\,{x}^{2}+588\,x-65 \right) ^{2}}{x}$.

}
\label{four_v44_f71_tor}}
\subFig{four_v44_f71_tor}%
\subFig{four_v44_f71_tor_ucover}%
\newSubFig{}{{$T(44|62)$. Valencies $(4,4|6,2)$. The order of the automorphism group: $2$.
Dual dessin   $T(62|44)$, see Figure \ref{four_v62_f44_tor}
on the page  \pageref{four_v62_f44_tor}.
Belyi function is $(X:{y}^{2}= \left( x-1 \right)  \left( 3\,{x}^{2}+8\,x+16 \right) ,
\beta={\frac {27}{256}}\,{\frac {{x}^{4}}{x-1}})
$.
} \label{four_v44_f62_tor}}
\subFig{four_v44_f62_tor}%
\subFig{four_v44_f62_tor_ucover}%
\newSubFig{}{{$T(44|44)\_8$. Valencies $(4,4|4,4)$. The order of the automorphism group: $8$.
Dual dessin   $T(44|44)\_8$, see Figure \ref{four_v44_f44_tor_8}
on the page  \pageref{four_v44_f44_tor_8}.
Belyi function is $ (X: y^2=x^4-1, \beta = x^4)$.
} \label{four_v44_f44_tor_8}}
\subFig{four_v44_f44_tor_8}%
\subFig{four_v44_f44_tor_8_ucover}%
\newSubFig{}{{$T(44|44)\_4$. Valencies $(4,4|4,4)$. The order of the automorphism group: $4$.
Dual dessin   $T(44|44)\_4$, see Figure \ref{four_v44_f44_tor_4}
on the page  \pageref{four_v44_f44_tor_4}.
Belyi function is $( X:{y}^{2}=\left( {x}^{2}-1 \right)  \left( {x}^{2}-2 \right),
\beta= \left( {x}^{2}-1 \right) ^{2} )
$.
} \label{four_v44_f44_tor_4}}
\subFig{four_v44_f44_tor_4}%
\subFig{four_v44_f44_tor_4_ucover}%
\end{arrangedFigure}
\clearpage
\begin{align*}
\langle\langle Tr(H^4)Tr(H^3)Tr(H)\rangle\rangle
=4\cdot 3\left(6N^3+2N\right).\
\end{align*}
\begin{arrangedFigure}{2}{1}{}{
{$S(431|611)_{+}$. Valencies $(4,3,1|6,1,1)$. The order of the automorphism group: $1$.
Dual dessin   $S(611|431)_{+}$, see Figure \ref{four_v611_f431_sf_plus}
on the page  \pageref{four_v611_f431_sf_plus}.
Belyi function is
$\beta=-{\frac {49}{4}}\,{\frac { \left( 87\,i\sqrt {3}+211 \right)  \left( z
-1 \right) ^{3} \left( -7\,z-3+2\,i\sqrt {3} \right) {z}^{4}}{686\,{z}
^{2}-672\,z+56\,i\sqrt {3}z-57\,i\sqrt {3}-51}}
$.
} \label{four_v431_f611_sf_plus}}
\subFig{four_v431_f611_sf_plus}%
\newSubFig{}{{$S(431|611)_{-}$. Valencies $(4,3,1|6,1,1)$. The order of the automorphism group:
$1$. Dual dessin   $S(611|431)_{-}$, see Figure \ref{four_v611_f431_sf_minus}
on the page  \pageref{four_v611_f431_sf_minus}.
Belyi function is
$\beta={\frac {49}{4}}\,{\frac { \left( 87\,i\sqrt {3}-211 \right)  \left( z-
1 \right) ^{3} \left( 7\,z+3+2\,i\sqrt {3} \right) {z}^{4}}{-686\,{z}^
{2}+56\,i\sqrt {3}z+672\,z+51-57\,i\sqrt {3}}}
$.
} \label{four_v431_f611_sf_minus}}
\subFig{four_v431_f611_sf_minus}%
\end{arrangedFigure}
\begin{arrangedFigure}{2}{1}{}{
{$S(431|521)$. Valencies $(4,3,1|5,2,1)$. The order of the automorphism group: $1$.
Dual dessin   $S(521|431)$, see Figure \ref{four_v521_f431_sf}
on the page  \pageref{four_v521_f431_sf}.
Belyi function is
$\beta=-1/16\,{\frac { \left( 4\,z-7 \right) ^{3} \left( 4\,z-15 \right) {z}^
{4}}{ \left( 3\,z+1 \right)  \left( 7\,z-6 \right) ^{2}}}
$.
} \label{four_v431_f521_sf}}
\subFig{four_v431_f521_sf}%
\newSubFig{}{{$S(431|431)A$. Valencies $(4,3,1|4,3,1)$. The order of the automorphism group: $1$.
Dual dessin   $S(431|431)A$, see Figure \ref{four_v431_f431_sf_a}
on the page  \pageref{four_v431_f431_sf_a}.
Belyi function is
$\beta={\frac {1}{3294172}}\,
{\frac { \left( 835+872\,\sqrt {2} \right)
 \left( -z+8+5\,\sqrt {2} \right)  \left( -z-8+9\,\sqrt {2} \right) ^{3}{z}^{4}}
{ \left( -z-11+8\,\sqrt {2} \right)  \left( z-1 \right) ^{3}}}
$.
} \label{four_v431_f431_sf_a}}
\subFig{four_v431_f431_sf_a}%
\newSubFig{}{{$S(431|431)B$. Valencies $(4,3,1|4,3,1)$. The order of the automorphism group:
$1$. Dual dessin   $S(431|431)B$, see Figure \ref{four_v431_f431_sf_b}
on the page  \pageref{four_v431_f431_sf_b}.
Belyi function is
$\beta={\frac {1}{3294172}}\,{\frac { \left( -835+872\,\sqrt {2} \right)
 \left( z-8+5\,\sqrt {2} \right)  \left( z+8+9\,\sqrt {2} \right) ^{3} {z}^{4}}
{ \left( z+11+8\,\sqrt {2} \right)  \left( z-1 \right) ^{3}}}
$.
} \label{four_v431_f431_sf_b}}
\subFig{four_v431_f431_sf_b}%
\newSubFig{}{{$S(431|422)$. Valencies $(4,3,1|4,2,2)$. The order of the automorphism group: $1$.
Dual dessin   $S(422|431)$, see Figure \ref{four_v422_f431_sf}
on the page  \pageref{four_v422_f431_sf}.
Belyi function is
$\beta=-4\,{\frac { \left( z+1 \right) ^{4}{z}^{3} \left( z+4 \right) }
{\left( 6\,{z}^{2}+4\,z+1 \right) ^{2}}}
$.
} \label{four_v431_f422_sf}}
\subFig{four_v431_f422_sf}%
\end{arrangedFigure}
\begin{arrangedFigure}{2}{2}{}{
{$T(431|8)A$. Valencies $(4,3,1|8)$. The order of the automorphism group: $1$. Dual dessin
  $T(8|431)A$, see Figure \ref{four_v8_f431_tor_a} on the page \pageref{four_v8_f431_tor_a}.
Belyi function is
$\beta=-{\frac {256}{85766121}} ( 5488\,{x}^{4}+14112\,y{x}^{2}-26264
\,{x}^{3}+37548\,yx-202741\,{x}^{2}+3240\,y-73368\,x-3240 )$
  on the curve
$X:{y}^{2}={\frac {1}{81}}\, ( 1-x ) ( 448\,{x}^{2}+1872\,x+81 )$.
$n_0={\frac {65536}{62523502209}}\, \left( 4\,x+45 \right)
 \left( 4\,x+21 \right) ^{3}{x}^{4}$,
$n_1={\frac {( 4096\,{x}^{4}+55296\,{x}^{3}+158976\,{x}^{2}+55296\,x-247617
) ^{2}}{62523502209}}
$.
}
\label{four_v431_f8_tor_a}}
\subFig{four_v431_f8_tor_a}%
\subFig{four_v431_f8_tor_a_ucover}%
\newSubFig{}{{$T(431|8)B$. Valencies $(4,3,1|8)$. The order of the automorphism group: $1$.
Dual dessin   $T(8|431)B$, see Figure \ref{four_v8_f431_tor_b}
on the page  \pageref{four_v8_f431_tor_b}.
Belyi function is
$\beta=-9/4\,{x}^{4}+3\,y{x}^{2}-8\,{x}^{3}+8/3\,yx
-{\frac {77}{9}}\,{x}^{2}+2/3\,y-{\frac {40}{9}}\,x-{\frac {8}{9}} $
  on the curve
$X:{y}^{2}=\frac{4}{9} \left( x+1 \right) ( 9\,{x}^{2}+4\,x+4 )$.
$n_0=\frac{3}{16} \left( x-2 \right)  ( 3x+2 ) ^{3}{x}^{4}$,
$n_1=\frac {1}{144}\, ( 27\,{x}^{4}-36\,{x}^{2}-32\,x-20) ^{2}$.
} \label{four_v431_f8_tor_b}}
\subFig{four_v431_f8_tor_b}%
\subFig{four_v431_f8_tor_b_ucover}%
\end{arrangedFigure}
\clearpage
\begin{align*}
\langle\langle Tr(H^4)Tr^2(H^2)\rangle\rangle
=4\cdot2^2\cdot 2!\left(\frac{3}{2}N^3+\frac{3}{4}N\right)=\\
=4\cdot2^2\cdot 2!\left(\left(1+\frac{1}{2}\right)N^3+
\left(\frac{1}{2}+\frac{1}{4}\right)N\right).\\
\end{align*}
\begin{figure}[h]
\begin{minipage}[b]{.45\linewidth}
\centering\epsfig{figure=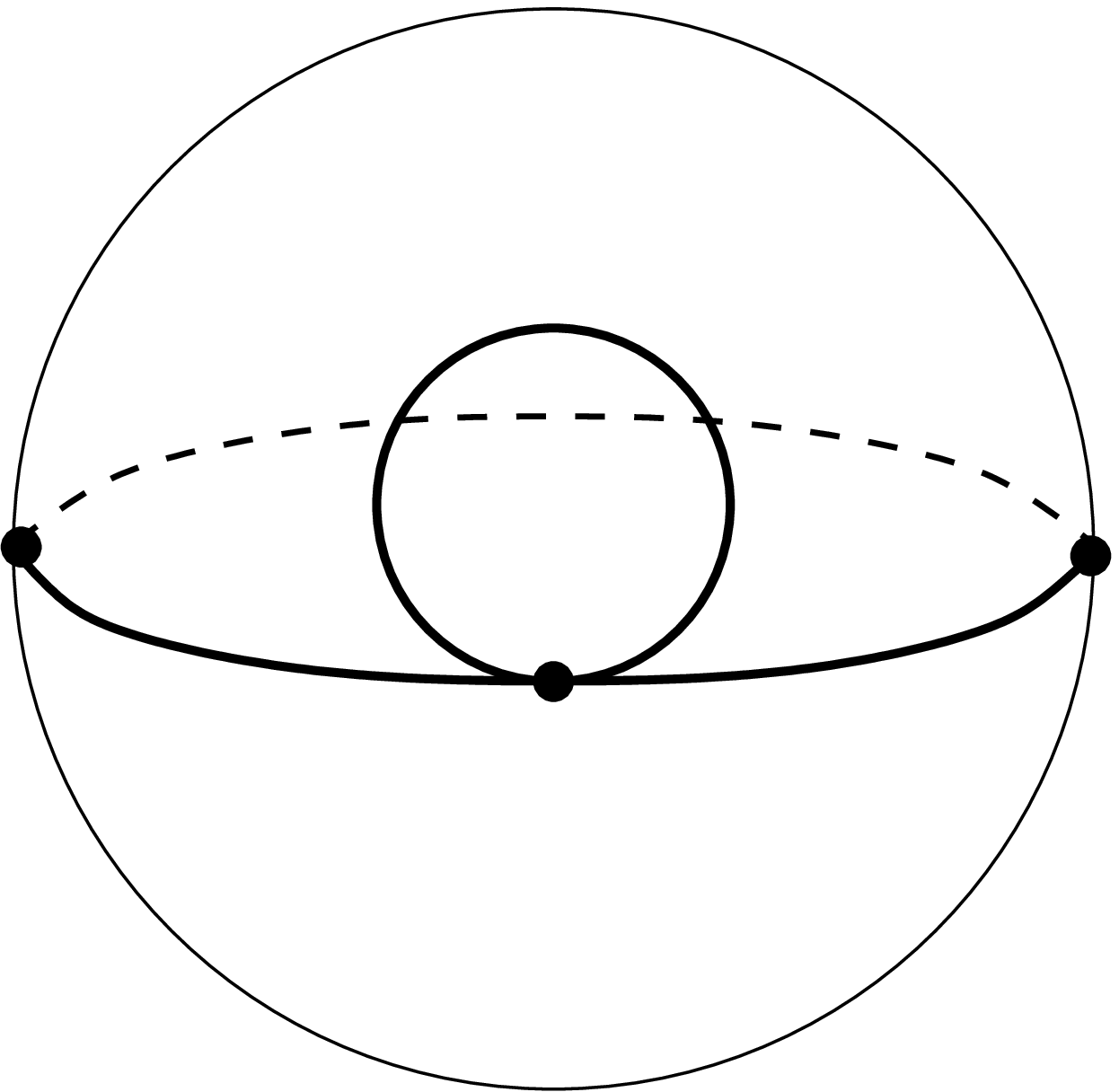,width=.85\linewidth}
\caption{$S(422|431)$. Valencies $(4,2,2|4,3,1)$. The order of the automorphism group: $1$.
Dual dessin   $S(431|422)$, see Figure \ref{four_v431_f422_sf}
on the page  \pageref{four_v431_f422_sf}.
Belyi function is
$\beta=-1/4\,{\frac { \left( {z}^{2}+4\,z+6 \right) ^{2}{z}^{4}}
{ \left( z+1 \right) ^{4} \left( 4\,z+1 \right) }}
$.
} \label{four_v422_f431_sf}
\end{minipage}\hfill
\begin{minipage}[b]{.45\linewidth}
\centering\epsfig{figure=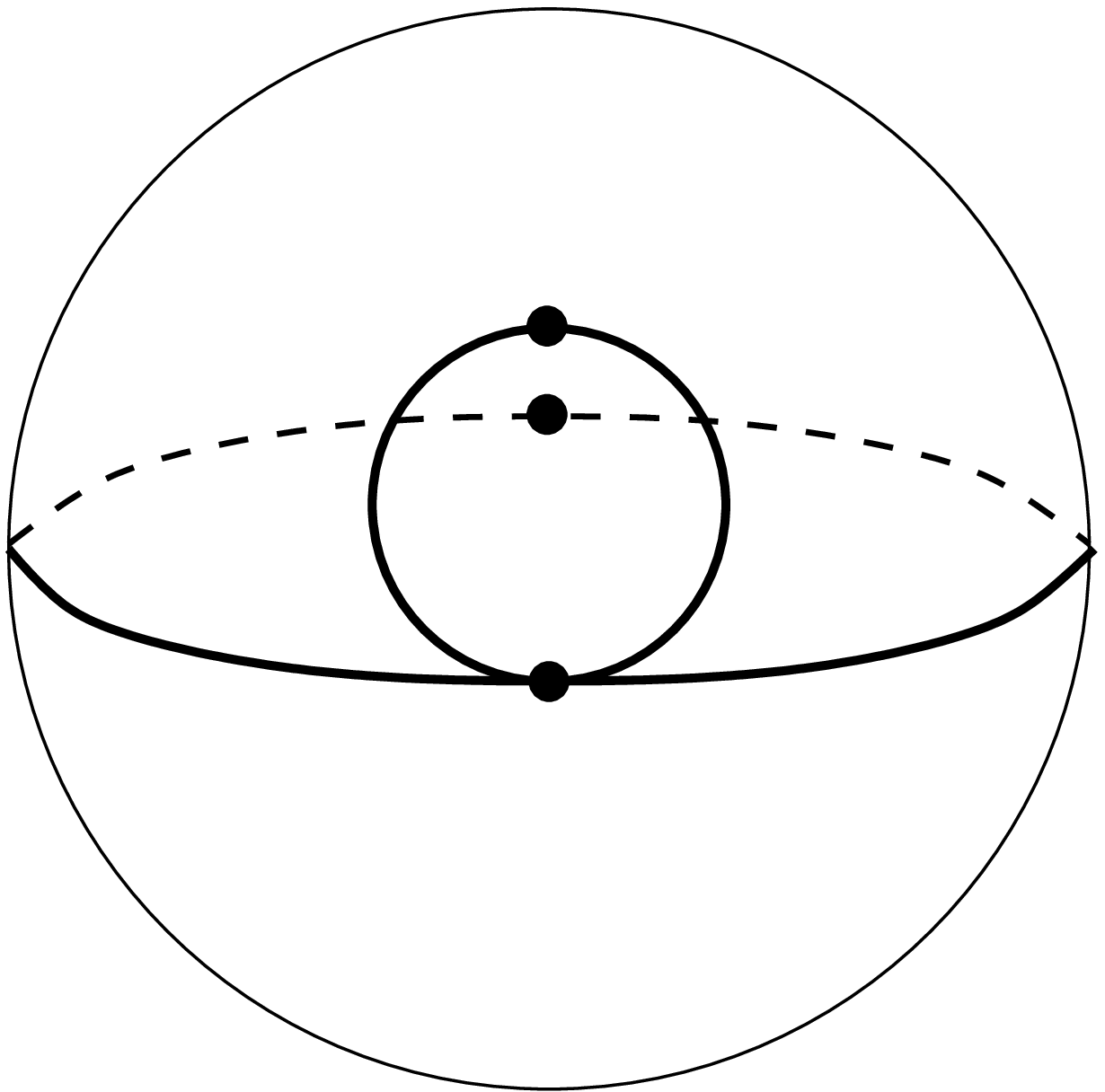,width=.85\linewidth}
\caption{$S(422|422)$. Valencies $(4,2,2|4,2,2)$. The order of the automorphism group: $2$.
Dual dessin   $S(422|422)$, see Figure \ref{four_v422_f422_sf}
on the page  \pageref{four_v422_f422_sf}.
Belyi function is
$\beta=-1/4\,{\frac { \left( {z}^{2}-2 \right) ^{2}{z}^{4}}
{ \left( z-1 \right) ^{2} \left( z+1 \right) ^{2}}}
$.
} \label{four_v422_f422_sf}
\end{minipage}\hfill
\end{figure}
\begin{arrangedFigure}{2}{2}{}{
{$T(422|8)\_4$. Valencies $(4,2,2|8)$. The order of the automorphism group: $4$. Dual dessin
  $T(8|422)\_4$, see Figure \ref{four_v8_f422_tor_4} on the page \pageref{four_v8_f422_tor_4}.
Belyi function is $(X:\,{y}^{2}=x ( x-1 )  ( x+1 ),
\beta=-4\,{x}^{2} ( x-1 )  ( x+1 ) )
$.
}
\label{four_v422_f8_tor_4}}
\subFig{four_v422_f8_tor_4}%
\subFig{four_v422_f8_tor_4_ucover}%
\newSubFig{}{{$T(422|8)\_2$. Valencies $(4,2,2|8)$. The order of the automorphism group: $2$.
Dual dessin   $T(8|422)\_2$, see Figure \ref{four_v8_f422_tor_2}
on the page  \pageref{four_v8_f422_tor_2}.
Belyi function is $(X:y^2=\left( {x}^{2}-2\,x-1 \right) x,
\beta= \left( x-2 \right) ^{2}{x}^{2})
$.
} \label{four_v422_f8_tor_2}}
\subFig{four_v422_f8_tor_2}%
\subFig{four_v422_f8_tor_2_ucover}%
\end{arrangedFigure}
\clearpage
\begin{align*}
\langle\langle Tr(H^4)Tr(H^2)Tr^2(H)\rangle\rangle
=4\cdot2\cdot 2!\left(\frac{9}{2}N^2\right)=\\
=4\cdot2\cdot 2!\left(\left(4+\frac{1}{2}\right)N^2\right).\\
\end{align*}
\begin{arrangedFigure}{2}{1}{}{
{$S(4211|71)_{+}$. Valencies $(4,2,1,1|7,1)$. The order of the automorphism group: $1$.
Dual dessin   $S(71|4211)_{+}$, see Figure \ref{four_v71_f4211_sf_plus}
on the page  \pageref{four_v71_f4211_sf_plus}.
Belyi function is
$\beta=\frac{(49 z^2-90z -2 i\sqrt 7 z+42+2i\sqrt 7) (128 z+5 i\sqrt 7 - 119)^2 z^4}
{512 (16377+181 i\sqrt 7 ) (z-1)^7}
$
} \label{four_v4211_f71_sf_plus}}
\subFig{four_v4211_f71_sf_plus}%
\newSubFig{}{{$S(4211|71)_{-}$. Valencies $(4,2,1,1|7,1)$. The order of the automorphism group:
$1$. Dual dessin   $S(71|4211)_{-}$, see Figure \ref{four_v71_f4211_sf_minus}
on the page  \pageref{four_v71_f4211_sf_minus}.
Belyi function is
$\beta=\frac{(49 z^2-90z +2 i\sqrt 7 z+42-2i\sqrt 7) (128 z-5 i\sqrt 7 - 119)^2 z^4}
{512 (16377-181 i\sqrt 7 ) (z-1)^7}
$
} \label{four_v4211_f71_sf_minus}}
\subFig{four_v4211_f71_sf_minus}%
\end{arrangedFigure}
\begin{arrangedFigure}{2}{1}{}{
{$S(4211|62)$. Valencies $(4,2,1,1|6,2)$. The order of the automorphism group: $1$.
Dual dessin   $S(62|4211)$, see Figure \ref{four_v62_f4211_sf}
on the page  \pageref{four_v62_f4211_sf}.
Belyi function is
$\beta=-108\,{\frac {{z}^{4} \left( 1+2\,z+3\,{z}^{2} \right)  \left( -1+3\,z
 \right) ^{2}}{ \left( -1+4\,z \right) ^{2}}}
$
} \label{four_v4211_f62_sf}}
\subFig{four_v4211_f62_sf}%
\newSubFig{}{{$S(4211|53)$. Valencies $(4,2,1,1|5,3)$. The order of the automorphism group: $1$. Dual dessin
  $S(53|4211)$, see Figure \ref{four_v53_f4211_sf} on the page \pageref{four_v53_f4211_sf}.
Belyi function is
$\beta=4\,{\frac { \left( 3\,z+1 \right) ^{2} \left( -9-42\,z+5\,{z}^{2}
 \right) }{ \left( z-8 \right) ^{3}{z}^{5}}}
$.
}
\label{four_v4211_f53_sf}}
\subFig{four_v4211_f53_sf}%
\newSubFig{}{{$S(4211|44)$. Valencies $(4,2,1,1|4,4)$. The order of the automorphism group: $2$.
Dual dessin   $S(44|4211)$, see Figure \ref{four_v44_f4211_sf}
on the page  \pageref{four_v44_f4211_sf}.
Belyi function is
$\beta=4\,{\frac {{z}^{2} \left( {z}^{2}+2 \right) }
{ \left( {z}^{2}+1 \right) ^{4}}}
$.
} \label{four_v4211_f44_sf}}
\subFig{four_v4211_f44_sf}%
\end{arrangedFigure}
\clearpage
\begin{align*}
\langle\langle Tr(H^4)Tr^4(H)\rangle\rangle
=4\cdot 4!\left(\frac{1}{4}N\right).\\
\end{align*}
\begin{figure}[h]
\begin{minipage}[b]{.45\linewidth}
\centering\epsfig{figure=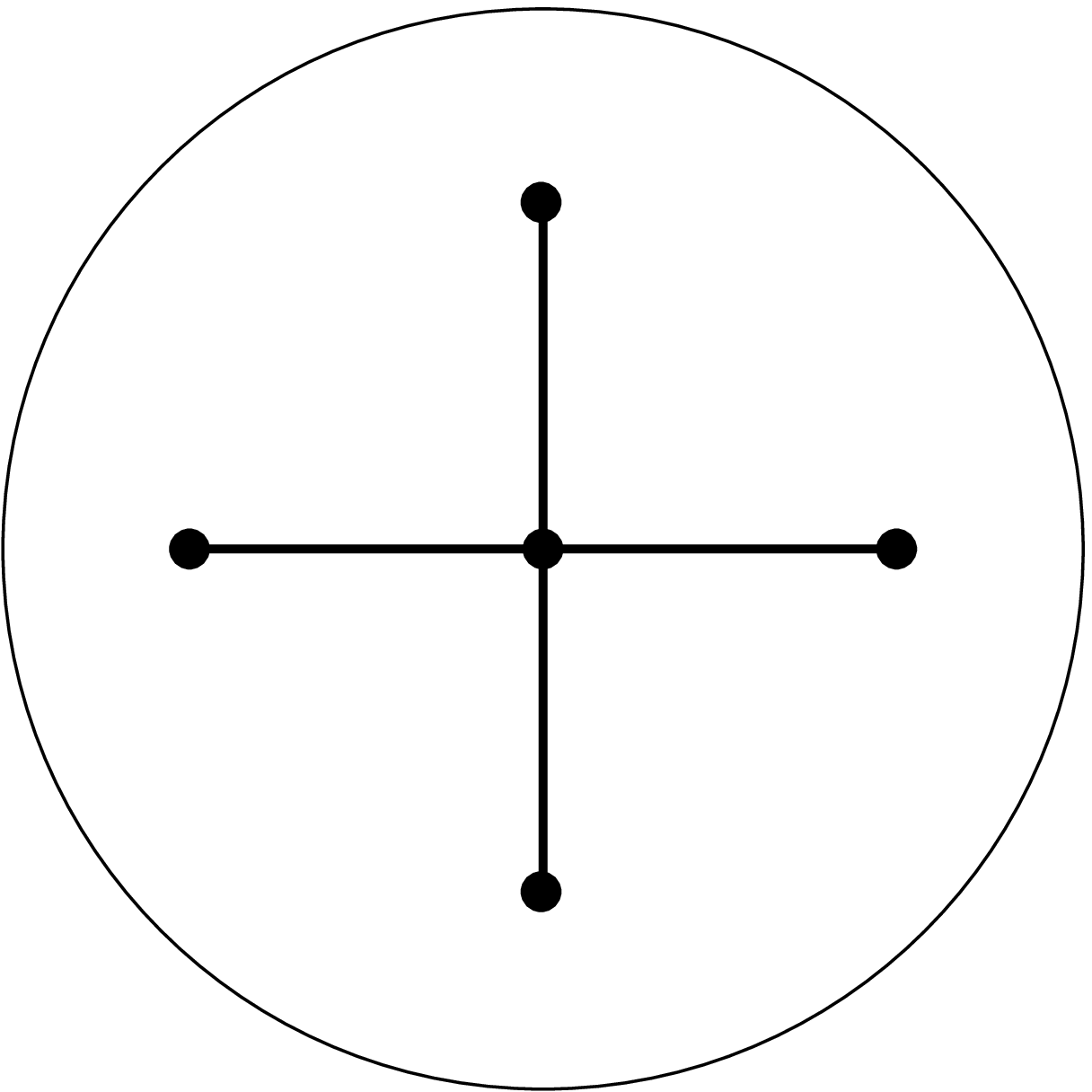,width=.85\linewidth}
\caption{$S(41111|8)$. Valencies $(4,1,1,1,1|8)$. The order of the automorphism group: $4$.
Dual dessin   $S(8|41111)$, see Figure \ref{four_v8_f41111_sf}
on the page  \pageref{four_v8_f41111_sf}.
Belyi function is $\beta=-4z^4(z-1)(z+1)(z^2+1)$.
} \label{four_v41111_f8_sf}
\end{minipage}\hfill
\end{figure}
\clearpage
\begin{align*}
\langle\langle Tr^2(H^3)Tr(H^2)\rangle\rangle
=3^2\cdot2\cdot 2!\left(2N^3+\frac{1}{2}N\right)=\\
=3^2\cdot2\cdot 2!\left(\left(1+2\cdot\frac{1}{2}\right)N^3+\frac{1}{2}N\right).\\
\end{align*}
\begin{arrangedFigure}{3}{1}{}{
{$S(332|611)$. Valencies $(3,3,2|6,1,1)$. The order of the automorphism group: $2$. Dual dessin
  $S(611|332)$, see Figure \ref{four_v611_f332_sf} on the page \pageref{four_v611_f332_sf}.
Belyi function is
$\beta=-1/4\,{\frac { \left( {z}^{2}+4 \right) ^{3}{z}^{2}}{2\,{z}^{2}+9}}$.
}
\label{four_v332_f611_sf}}
\subFig{four_v332_f611_sf}%
\newSubFig{}{{$S(332|521)$. Valencies $(3,3,2|5,2,1)$. The order of the automorphism group: $1$.
Dual dessin   $S(521|332)$, see Figure \ref{four_v521_f332_sf}
on the page  \pageref{four_v521_f332_sf}.
Belyi function is
$\beta={\frac {4}{27}}\,{\frac { \left( z+4 \right) ^{2} \left( {z}^{2}-6
 \right) ^{3}}{ \left( 2\,z+7 \right) ^{2} \left( 4\,z-11 \right) }}
$.
} \label{four_v332_f521_sf}}
\subFig{four_v332_f521_sf}%
\newSubFig{}{{$S(332|332)$. Valencies $(3,3,2|3,3,2)$. The order of the automorphism group: $2$.
Dual dessin   $S(332|332)$, see Figure \ref{four_v332_f332_sf}
on the page  \pageref{four_v332_f332_sf}.
Belyi function is
$\beta=64\,{\frac {{z}^{2} \left( {z}^{2}+1 \right) ^{3}}
{ \left( 8\,{z}^{2}-1 \right) ^{3}}}
$.
} \label{four_v332_f332_sf}}
\subFig{four_v332_f332_sf}%
\end{arrangedFigure}
\begin{arrangedFigure}{1}{2}{}{
{$T(332|8)$. Valencies $(3,3,2|8)$. The order of the automorphism group: $2$. Dual dessin
  $T(8|332)$, see Figure \ref{four_v8_f332_tor} on the page \pageref{four_v8_f332_tor}.
Belyi function is $(X:{y}^{2}=\left( x-2 \right)  \left( 4\,{x}^{2}+4\,x+3 \right) ,
\beta=-{\frac {16}{27}}\, \left( x-2 \right) {x}^{3})
$.
}
\label{four_v332_f8_tor}}
\subFig{four_v332_f8_tor}%
\subFig{four_v332_f8_tor_ucover}%
\end{arrangedFigure}
\clearpage
\begin{align*}
\langle\langle Tr^2(H^3)Tr^2(H)\rangle\rangle
=3^2\cdot2!\cdot 2!\left(2N^2\right)=\\
=3^2\cdot2!\cdot 2!\left(\left(1+2\cdot\frac{1}{2}\right)N^2\right).\\
\end{align*}
\begin{arrangedFigure}{3}{1}{}{
{$S(3311|71)$. Valencies $(3,3,1,1|7,1)$. The order of the automorphism group: $1$. Dual dessin
  $S(71|3311)$, see Figure \ref{four_v71_f3311_sf} on the page \pageref{four_v71_f3311_sf}.
Belyi function is
$\beta=-{\frac {1}{1728}}\,{\frac { \left( 1+{z}^{2}-5\,z \right) ^{3}
 \left( 49\,{z}^{2}-13\,z+1 \right) }{{z}^{7}}}
$
}
\label{four_v3311_f71_sf}}
\subFig{four_v3311_f71_sf}%
\newSubFig{}{{$S(3311|62)$. Valencies $(3,3,1,1|6,2)$. The order of the automorphism group: $2$.
Dual dessin   $S(62|3311)$, see Figure \ref{four_v62_f3311_sf}
on the page  \pageref{four_v62_f3311_sf}.
Belyi function is $\beta=-{\frac {1}{64}}\,{\frac { \left( z-3 \right)  \left( 3+z \right)
 \left( z-1 \right) ^{3} \left( z+1 \right) ^{3}}{{z}^{2}}}
$
} \label{four_v3311_f62_sf}}
\subFig{four_v3311_f62_sf}%
\newSubFig{}{{$S(3311|44)$. Valencies $(3,3,1,1|4,4)$. The order of the automorphism group: $2$.
Dual dessin   $S(44|3311)$, see Figure \ref{four_v44_f3311_sf}
on the page  \pageref{four_v44_f3311_sf}.
Belyi function is
$\beta=-{\frac {1}{432}}\,{\frac { \left( {z}^{2}-10\,z-2 \right)  \left( {z}
^{2}-2\,z-2 \right) ^{3}}{{z}^{4}}}
$.
} \label{four_v3311_f44_sf}}
\subFig{four_v3311_f44_sf}%
\end{arrangedFigure}
\clearpage
\begin{align*}
\langle\langle Tr(H^3)Tr^2(H^2)Tr(H)\rangle\rangle
=3\cdot2^2\cdot 2!\left(3N^2\right).\\
\end{align*}
\begin{arrangedFigure}{3}{1}{}{
{$S(3221|71)$. Valencies $(3,2,2,1|7,1)$. The order of the automorphism group: $1$. Dual dessin
  $S(71|3221)$, see Figure \ref{four_v71_f3221_sf} on the page \pageref{four_v71_f3221_sf}.
Belyi function is
$\beta=256\,{\frac {{z}^{3} \left( z+1 \right)  \left( 7+28\,z+24\,{z}^{2}
 \right) ^{2}}{48\,z-1}}$.

}
\label{four_v3221_f71_sf}}
\subFig{four_v3221_f71_sf}%
\newSubFig{}{{$S(3221|62)$. Valencies $(3,2,2,1|6,2)$. The order of the automorphism group: $1$.
Dual dessin   $S(62|3221)$, see Figure \ref{four_v62_f3221_sf}
on the page  \pageref{four_v62_f3221_sf}.
Belyi function is
$\beta=-\frac14\,{\frac {{z}^{3} \left( 4+z \right)  \left( {z}^{2}+2\,z-2
 \right) ^{2}}{ \left( -1+2\,z \right) ^{2}}}
$.
} \label{four_v3221_f62_sf}}
\subFig{four_v3221_f62_sf}%
\newSubFig{}{{$S(3221|53)$. Valencies $(3,2,2,1|5,3)$. The order of the automorphism group: $1$.
Dual dessin   $S(53|3221)$, see Figure \ref{four_v53_f3221_sf}
on the page  \pageref{four_v53_f3221_sf}.
Belyi function is
$\beta=\frac14\,{\frac { \left( z+7 \right)  \left( 3+z \right) ^{3} \left( {z}^{
2}+1 \right) ^{2}}{ \left( 3+4\,z \right) ^{3}}}
$.
} \label{four_v3221_f53_sf}}
\subFig{four_v3221_f53_sf}%
\end{arrangedFigure}
\clearpage
\begin{align*}
\langle\langle Tr(H^3)Tr(H^2)Tr^3(H)\rangle\rangle
=3\cdot2\cdot 3!\left(N\right).\\
\end{align*}
\begin{figure}[h]
\begin{minipage}[b]{.45\linewidth}
\centering\epsfig{figure=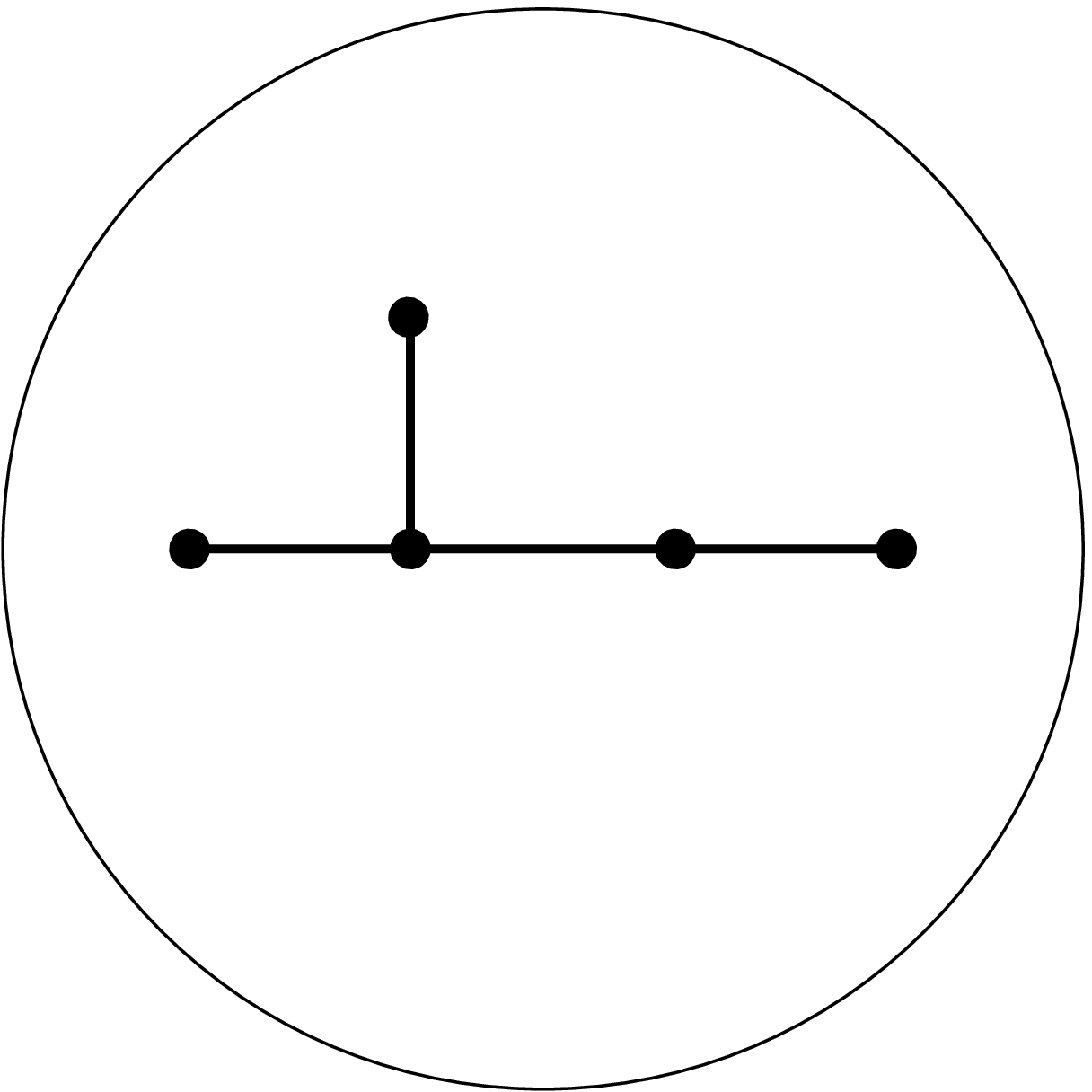,width=.85\linewidth}
\caption{$S(32111|8)$. Valencies $(3,2,1,1,1|8)$. The order of the automorphism group: $1$.
Dual dessin   $S(8|32111)$, see Figure \ref{four_v8_f32111_sf}
on the page  \pageref{four_v8_f32111_sf}.
Belyi function is $\beta=-\frac{1024}{729}z^3(z-1)(16z^2+8z+3)(4z-3)^2$
} \label{four_v32111_f8_sf}
\end{minipage}\hfill
\end{figure}
\clearpage
\begin{align*}
\langle\langle Tr^4(H^2)\rangle\rangle
=2^4\cdot 4!\left(\frac{1}{8}N^2\right).\\
\end{align*}
\begin{figure}[h]
\begin{minipage}[b]{.45\linewidth}
\centering\epsfig{figure=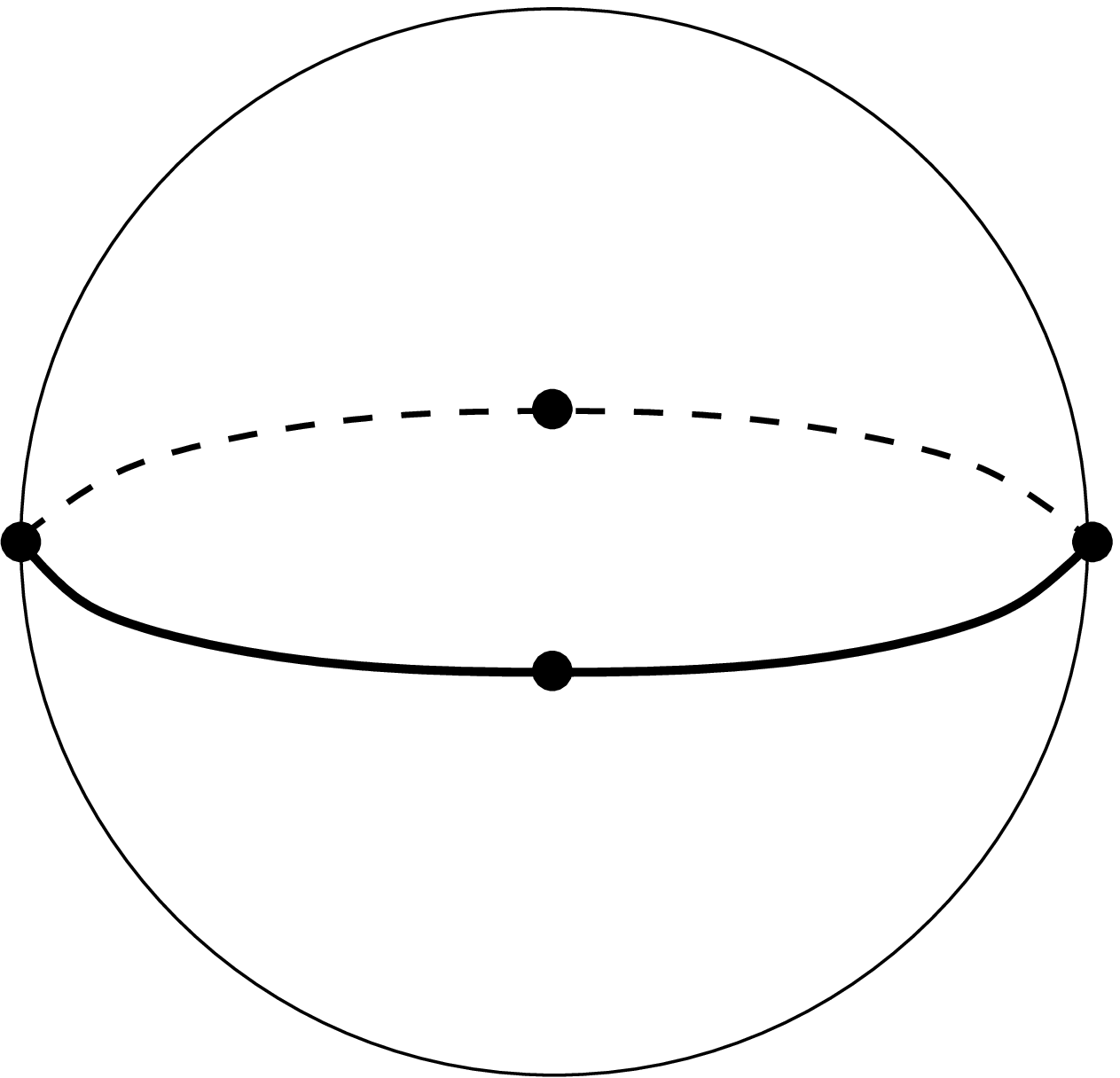,width=.85\linewidth}
\caption{$S(2222|44)$. Valencies $(2,2,2,2|4,4)$. The order of the automorphism group: $8$.
Dual dessin   $S(44|2222)$, see Figure \ref{four_v44_f2222_sf}
on the page  \pageref{four_v44_f2222_sf}.
Belyi function is $\beta=\frac{(z^4+1)^2}{4z^4}$.
} \label{four_v2222_f44_sf}
\end{minipage}\hfill
\end{figure}
\clearpage
\begin{align*}
\langle\langle Tr^3(H^2)Tr^2(H)\rangle\rangle
=2^3\cdot 3!\cdot 2!\left(\frac{1}{2}N\right).\\
\end{align*}
\begin{figure}[h]
\begin{minipage}[b]{.45\linewidth}
\centering\epsfig{figure=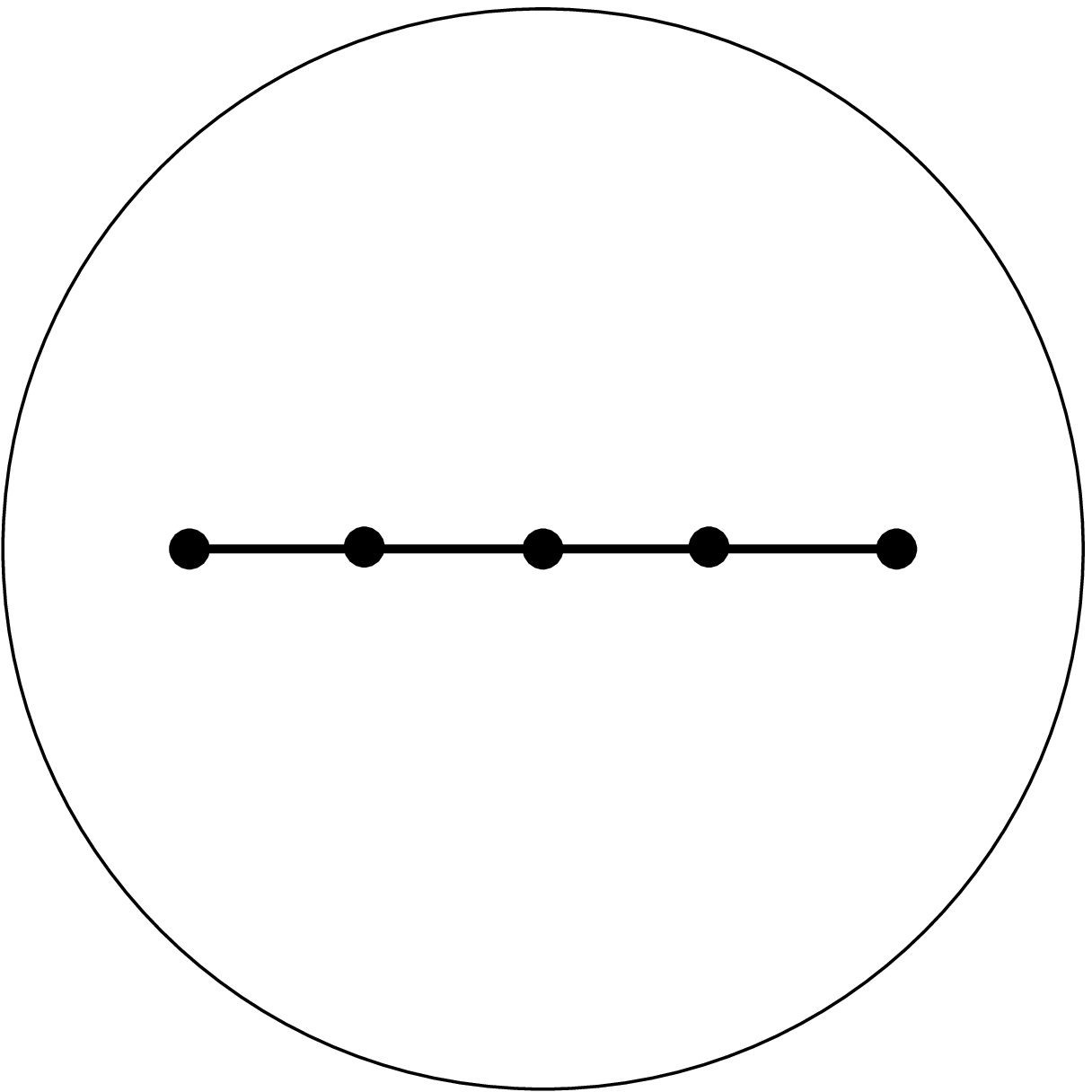,width=.85\linewidth}
\caption{$S(22211|8)$. Valencies $(2,2,2,1,1|8)$. The order of the automorphism group: $2$.
Dual dessin   $S(8|22211)$, see Figure \ref{four_v8_f22211_sf}
on the page  \pageref{four_v8_f22211_sf}.
Belyi function is $\beta=-4 z^2(z^2-2) (z-1)^2 (z+1)^2$
} \label{four_v22211_f8_sf}
\end{minipage}\hfill
\end{figure}

\end{document}